\documentclass[12pt]{amsart}
\usepackage{amsmath,amsthm,amsfonts,amssymb,latexsym,mathrsfs}
\usepackage[top=1in, bottom=1in, left=1in, right=1in]{geometry}
\usepackage{graphics, color}
\usepackage{bm}
\usepackage{times}

\usepackage{xspace}
\usepackage{hyperref}
\usepackage{hhline}

\definecolor{myblue}{rgb}{0.09,0.32,0.44} 
\hypersetup{pdfborder={0 0 0},pdfborderstyle={},colorlinks=true,linkcolor=myblue,citecolor=myblue,
urlcolor=black}

\usepackage[shortlabels]{enumitem}

\makeatletter
\renewcommand\part{\@startsection{part}{2}%
	\z@{0.5\linespacing\@plus2\linespacing}{\linespacing}%
	{\normalfont\large\scshape\bfseries\centering}}

\usepackage[utf8]{inputenc}

\usepackage{upgreek}
\usepackage{enumitem}

\usepackage{graphicx,tikz}
\usetikzlibrary{shapes.geometric, arrows}

\usepackage [autostyle, english = american]{csquotes}
\MakeOuterQuote{"}
\usepackage{amsthm}

\numberwithin{equation}{section}
\newtheorem{theorem}{Theorem}[section]
\newtheorem{lemma}[theorem]{Lemma}
\newtheorem{cor}[theorem]{Corollary}
\newtheorem{prop}[theorem]{Proposition}

\newtheorem{thm}{Theorem}

\theoremstyle{remark}
\newtheorem*{remark}{Remark}
\newtheorem{rem}{Remark}[section]

\theoremstyle{definition}
\newtheorem{definition}[theorem]{Definition}

\newcommand{\imod}[1]{\ (\mathrm{mod}\ #1)}

\usepackage[utf8]{inputenc}
\usepackage[T1]{fontenc}

\begin{document}

\title[Primes with a missing digit in arithmetic progressions]{Primes with a missing digit: distribution in arithmetic progressions and an application in sieve theory}

\author[Kunjakanan Nath]{Kunjakanan Nath}

\address{Department of Mathematics,
University of Illinois Urbana-Champaign,
Altgeld Hall, 1409 W. Green Street,
Urbana, IL, 61801, USA}
\email{\href{kunjakanan@gmail.com}{\nolinkurl{kunjakanan@gmail.com}}}

\date{\today}

\begin{abstract}
We prove Bombieri-Vinogradov type theorems for primes with a missing digit in their $b$-adic expansion for some large positive integer $b$. The proof is based on the circle method, which relies on the Fourier structure of the integers with a missing digit and the exponential sums over primes in arithmetic progressions. 

Combining our results with the semi-linear sieve, we obtain an upper bound and a lower bound of the correct order of magnitude for the number of primes of the form $p=1+m^2+n^2$ with a missing digit in a large odd base $b$.
\end{abstract}

\maketitle

\setcounter{tocdepth}{1}

\tableofcontents

\part{Main results and outline of the proof}

\section{Introduction}\label{sec: Intro}
Let $b\geq 3$ be an integer and let $a_0\in\{0, 1, \dotsc, b-1\}$. Consider 
\[\mathcal{A}:=\bigg\{\sum_{j\geq 0}n_jb^j: n_j\in \{0, \dotsc, b-1\}\setminus \{a_0\}\: \text{for all $j$}\bigg\},\]
the set of non-negative integers without the digit $a_0$ in their $b$-adic expansion. For any $k\in \mathbb{N}$, the cardinality of the set $\mathcal{A}\cap [1, b^k)$ is $\approx (b-1)^k$. 

For the rest of the paper, we set $X=b^k$ and note that there are $\approx X^{\zeta}$ elements in $\mathcal{A}$ less than $X$, where
\[\zeta :=\frac{\log(b-1)}{\log b}<1.\] 
This reveals that $\mathcal{A}$ is a ``sparse set''. It is often the case that sparseness is one of the obstacles in analytic number theory. However, the set $\mathcal{A}$ admits some interesting structure in the sense that its Fourier transform has an explicit description, which is often small in size. There has been a considerable amount of work (see Dartyge-Mauduit \cite{darmau, darmau2}, Erd\H{o}s-Mauduit-S\'ark\"ozy \cite{erdmausar, erdmausar2}, Konyagin \cite{kon}, Maynard \cite{May, May1, May2}, Pratt \cite{Pratt}) in this direction by exploiting the Fourier structure of the set $\mathcal{A}$.

\begin{remark}
Note that $\zeta \rightarrow 1$ as $b\rightarrow \infty$.  We shall have many occasions to use
this fact in the paper, and we do so without further comment.
\end{remark}

It is a natural question to ask if the set $\mathcal{A}$ contains infinitely many primes. We expect the answer to be affirmative. In his celebrated paper \cite{May1}, Maynard showed that the set $\mathcal{A}$ contains infinitely many primes for any base $b\geq 10$. Moreover, for a large base, say $b\geq 2\times 10^6$, he established an asymptotic formula (see \cite[Theorem 2.5]{May} or \cite[Theorem 1.1]{May2}).

Prior to Maynard's work, Dartyge-Mauduit \cite{darmau, darmau2} showed the existence of infinitely many almost-primes (integers with at most $2$ prime factors) in $\mathcal{A}$ for any base $b\geq 3$. They used crucially the fact that $\mathcal{A}$ is well-distributed in arithmetic progressions (see, for example, \cite{darmau}, \cite{erdmausar}). In that spirit, we are interested in understanding how the primes of $\mathcal{A}$ are distributed in arithmetic progressions.  For $(c, d)=1$ and $(d,b)=1$, one expects that as $X\to \infty$,
\begin{align}
\notag \#\{p<X:p\equiv c\imod d, p\in \mathcal{A}\} \sim \dfrac{1}{\varphi(d)}\#\{p<X:p\in \mathcal{A}\} 
\end{align}
holds uniformly for $d\leq X^{\zeta(1-\varepsilon)}$ with any fixed $\varepsilon>0$. This seems to be a difficult question at present. Instead, we aim for a Bombieri-Vinogradov Theorem of the following type:
\begin{align}
\notag \sum_{\substack{d\leq D\\(d,b)=1}}\max_{(c, d)=1}\bigg|\#\{p<X:p\equiv c\imod d, p\in \mathcal{A}\} - \dfrac{1}{\varphi(d)}\#\{p<X:p\in \mathcal{A}\}\bigg| \ll_{A,b} \dfrac{X^{\zeta}}{(\log X)^A},
\end{align}
where $D\leq X^{1/2-\varepsilon}$, for any fixed $\varepsilon>0$, provided that $b$ is large enough in terms of $\varepsilon$ (so that $\zeta$ is close enough to $1$). However, using the current techniques, we are not able to prove that the above estimate holds for $D\leq X^{1/2-\varepsilon}$. Nevertheless, we can make some progress in this direction. 

For technical convenience, we will work with the von Mangoldt function $\Lambda$ (recall that $\Lambda(n)=\log p$ if $n=p^m$, and $0$ otherwise). For $X=b^k$ with $k\in \mathbb{N}$ and for $(c, d)=(r, b)=1$, we set
\begin{align}
\notag \mathcal{E}(X; d, c; b, r):=\sum_{\substack{n < X\\n\equiv c \imod d\\n\equiv r\imod b}}\Lambda(n)1_{\mathcal{A}}(n) - \dfrac{1}{\varphi(d)}\dfrac{b}{\varphi(b)}\sum_{\substack{n<X\\n\equiv r\imod b}}1_{\mathcal{A}}(n).
\end{align}
Note that
\[\sum_{\substack{n<X\\ n\equiv r\imod b}}1_{\mathcal{A}}(n) = \dfrac{X^{\zeta}}{b-1}\]
whenever $r\not\equiv a_0\imod b$; otherwise, both sums in the definition of $\mathcal{E}(X; d, c;b,r)$ is $0$. Moreover, the condition $n\equiv r\imod b$ is equivalent to $n$ having $r$ as its last digit in its $b$-adic expansion. We add this condition in order to simplify some technical details later on.

\begin{thm}\label{Main Theorem}
Let $\delta>0$ and let $b$ be an integer that is sufficiently large in terms of $\delta$. Let $D\in [1, X^{1/3-\delta}]$ and let $r\in \mathcal{A}\cap[0, b)$ be an integer such that $(r, b)=1$. Then for any $A>0$, we have
\begin{align}
\label{MT1eq1} \sum_{\substack{d\leq D\\(d,b)=1}}\max_{(c,d)=1}\big|\mathcal{E}(X;d,c;b,r)\big| \ll_{A, b, \delta} \dfrac{X^{\zeta}}{(\log X)^A}.
\end{align}
\end{thm}

We can do a little better if we allow the moduli to be the product of two integers. However, the parameter $c$ is now fixed, so we must drop the expression $\max_{(c, d)=1}$ from \eqref{MT1eq1}.

\begin{thm}\label{Main Theorem2}
Let $\delta>0$,  let $b$ be an integer that is sufficiently large in terms of $\delta$, and let 
\[D_1\in [1, X^{1/3-\delta}] \quad \text{and} \quad D_2\in [1, X^{1/9}].\]
Let $c$ be a non-zero integer and let $r\in \mathcal{A}\cap[0, b)$ with $(r, b)=1$. Then for any $A>0$, we have
\begin{align}
\notag \sum_{\substack{d_1\leq D_1}} \sideset{}{^*}\sum_{\substack{d_2\leq D_2}}\big|\mathcal{E}(X; d_1d_2, c; b, r)\big| \ll_{A, b, \delta} \dfrac{X^{\zeta}}{(\log X)^A},
\end{align}
where $^*$ in the sum denotes the conditions $(c,d_1d_2)=(b,d_1d_2)=(d_1, d_2)=1$.
\end{thm}

We can further have better result in this direction when we replace the absolute value inside the sum over $d$ by a \emph{well-factorable function}. Before proceeding to state our result, we formally define the ``well-factorable'' function.

\begin{definition}[Well-factorable]\label{def: well fac}
Let $D\geq 1$ be a real number. We say an arithmetic function $\xi:\mathbb{N}\rightarrow\mathbb{R}$ \emph{well-factorable} of level $D$ if, for any choice of factorization $D = D_1D_2$ with $D_1, D_2 \geq 1$, there exist two arithmetic functions $\xi_1, \xi_2:\mathbb{N}\rightarrow\mathbb{R}$ such that
\begin{enumerate}[(i)]
\item $|\xi_1|, |\xi_2|\leq 1$.
\item $\xi_1$ is supported on $[1, D_1]$ and $\xi_2$ is supported on $[1, D_2]$.
\item We have\footnote{In general, we do not require the co-primality condition in the definition of $\xi$. However, in order to avoid some technical issues, we impose this condition here.}
\[\xi(d)=\sum_{\substack{d_1d_2=d\\(d_1, d_2)=1}}\xi_1(d_1)\xi_2(d_2).\]
\end{enumerate}
\end{definition}
With this definition, we are now ready to state the following result.

\begin{thm}\label{Main Theorem3}
Let $\delta>0$ and let $b$ be an integer that is sufficiently large in terms of $\delta$. Let $\xi:\mathbb{N}\rightarrow\mathbb{R}$ be a well-factorable arithmetic function of level $D\in [1, X^{1/2-\delta}]$. Let $c$ be a non-zero integer and let $r\in \mathcal{A}\cap[0,b)$ be such that $(r, b)=1$. Then, for any $A>0$, we have
\begin{align}
\notag \sum_{\substack{d\leq D\\ (d,bc)=1}}\xi(d)\mathcal{E}(X; d, c; b, r) \ll_{A,b,\delta}\dfrac{X^{\zeta}}{(\log X)^A}.
\end{align}
\end{thm}
\begin{remark}
We have not explicitly mentioned in the above three theorems the size of $b$. In fact, it will be evident from the proof that $\delta$ and $b$ are \emph{inversely related} to each other. A simple calculation suggests that the size of $b$ is approximately of order $10^{632}$ if we take $\delta=1/100$. Therefore, we will refrain from explicitly calculating $\delta$ and $b$ in the above three theorems.
\end{remark}

The key point of Theorem \ref{Main Theorem2} and Theorem \ref{Main Theorem3} is the quantitative improvement over Theorem \ref{Main Theorem} allowing us to handle moduli as large as $X^{4/9-\delta}$ and $X^{1/2-\delta}$, respectively (instead of $X^{1/3-\delta}$). However, Theorem \ref{Main Theorem2} has the disadvantage that it has a stronger requirement that the moduli need to be composite and Theorem \ref{Main Theorem3} requires moduli weighted by a well-factorable function. But, in some of the applications in sieve theory, we do have well-factorable moduli. In fact, we give such an application in this paper: we prove the existence of infinitely many primes of the form $p=1+m^2+n^2$ with a missing digit in a large odd base $b$. The following theorem gives a precise statement.

\begin{thm}\label{Main TheoremApp}
Let $b$ be an odd integer that is sufficiently large, and let 
\[\mathbb{B}=\{n:n=n_1^2+n_2^2~~\text{for some}~(n_1, n_2)=1\}.\]
Let $r\in \mathcal{A}\cap[0, b)$ with $\big(r(r-1), b\big)=1$. Then, we have
\begin{equation}
\notag \sum_{\substack{p < X\\p\equiv r\imod b}}1_{\mathcal{A}}(p)1_{\mathbb{B}}(p-1)\asymp_{b} \dfrac{X^{\zeta}}{(\log X)^{3/2}}.
\end{equation}
\end{thm}

\begin{remark}
The implicit upper bound in Theorem \ref{Main TheoremApp} follows from Theorem \ref{Main Theorem} and a standard upper bound sieve estimate (for example, see Lemma \ref{fundsemi}). However, for the lower bound, we need to be more careful and use an argument due to Iwaniec \cite{iwa, iwa2} that allows sieving for primes of the form $1+m^2+n^2$ using \emph{level of distribution} slightly less than $X^{1/2}$. Additionally, in order to use the sieve estimates efficiently, we need two technical results, namely, Theorem \ref{Hyp1} and Theorem \ref{Hyp2} (similar in nature to Theorems \ref{Main Theorem2}, \ref{Main Theorem3}).
\end{remark}

\subsection*{Acknowledgements}
The author would like to thank his Ph.D. advisor, Dimitris Koukoulopoulos, for many helpful conversations and support during the course of this project, for his valuable comments on the earlier versions of the manuscript and for suggesting some key inputs. The author would also like to thank Kevin Ford, Andrew Granville, Joni Ter\"{a}v\"{a}inen and the anonymous referee for their useful comments and suggestions.

The author was supported by bourse de doctorat en recherche (B2X) of  Fonds de recherche du Qu\'ebec - Nature et technologies (FRQNT); bourse de fin d’\'etudes doctorales of \'Etudes sup\'erieures et postdoctorales (ESP), and bourse Ars\`ene David of Universit\'e de Montr\'eal while carrying out this work.

\subsection*{Notations} We employ some standard notation that will be used throughout the paper.
\begin{itemize}
\item Expressions of the form $\mathfrak{f}(X)=O(\mathfrak{g}(X))$, $\mathfrak{f}(X) \ll \mathfrak{g}(X)$ and $\mathfrak{g}(X) \gg \mathfrak{f}(X)$ signify that $|\mathfrak{f}(X)| \leq C|\mathfrak{g}(X)|$ for all sufficiently large $X$, where $C>0$ is an absolute constant. A subscript of the form $\ll_A$ means the implied constant may depend on the parameter $A$. The notation $\mathfrak{f}(X) \asymp \mathfrak{g}(X)$ indicates that $\mathfrak{f}(X) \ll \mathfrak{g}(X) \ll \mathfrak{f}(X)$. Here all the quantities should be thought of $X=b^k$ with $k$ an integer and $k\rightarrow \infty$. 
\item All sums, products and maxima will be taken over $\mathbb{N}=\{1, 2, \dotsc\}$ unless specified otherwise.
 \item We reserve the letters $p, p^\prime, p_1, p_2$ to denote primes. 
\item The letter $\gamma$ will always denote the Euler-Mascheroni constant.
\item As usual, $\mathbb{R}$ will denote the set of real numbers, $\mathbb{P}$ the set of primes and $\mathbb{Z}$ the set of integers. Furthermore, $p^v\|m$ means that $p^v|m$ and $p^{v+1}\nmid m$.
\item Throughout the paper, $\varphi$ will denote the totient function, $\mu$ the M\"{o}bius function, and $\uptau_h(n)$ the number of ways of writing $n$ as a product of $h$ natural numbers. 
\item As it is customary, we denote $e(y)= e^{2\pi i y}$ for any real number $y$. We write $n\sim N$ to denote $ N<n\leq 2N$. We use $\|y\|$ to denote $\min_{n\in \mathbb{Z}}|y-n|$. 
\item Unless otherwise specified, $\chi$ will always denote a Dirichlet character modulo some positive integer. The symbol $\chi_0$ will always denote a principal character.
\item We will set $(a, b)$ to be the greatest common divisor of integers $a$ and $b$ and by abuse of notation it will also denote the open interval on the real line. On the other hand, $[a, b]$ will denote the closed interval on the real line, and sometimes it will denote the least common multiple of integers $a$ and $b$. Its exact meaning will always be clear from the context.
\item For co-prime integers $m$ and $n$, we set $\overline{n}$ to denote the inverse of $n$ modulo $m$, that is, $n\overline{n}\equiv 1\imod m$.
\item We let $1_\mathcal{E}$ to be the characteristic function of the set $\mathcal{E}$ (so $1_{\mathcal{E}}(x)=1$ if $x\in \mathcal{E}$, and $0$, otherwise).
\item For any set $\mathcal{E}$, $\#\mathcal{E}$ denotes the cardinality of the set $\mathcal{E}$.
\item For any two arithmetic functions $\mathfrak{f}$ and $\mathfrak{g}$, we write $(\mathfrak{f}*\mathfrak{g})(n):= \sum_{ab=n}\mathfrak{f}(a)\mathfrak{g}(b)$ for their Dirichlet convolution.
\item For any arithmetic function $\mathfrak{f}:\mathbb{N}\rightarrow \mathbb{C}$, we set $\lVert \mathfrak{f}\rVert_2: = \big(\sum_{n}|\mathfrak{f}(n)|^2\big)^{1/2}$.
\item For any arithmetic function $F$, we also set $F_{\leq U}(n):=F(n)\cdot 1_{n\leq U}$ and $F_{>U}(n):=F(n)\cdot 1_{n>U}$.
\item We set $\mathbb{B}=\{n\in \mathbb{Z}: n=n_1^2+n_2^2~~\text{for some}~(n_1, n_2)=1\}$ and $\mathcal{B}=\{n\geq 1: p|n\Rightarrow p\equiv 1\imod 4\}$.
\end{itemize}

We set $X=b^k$ with $k\in \mathbb{N}$ and $k\rightarrow \infty$ for the rest of the paper except Part \ref{part: exp sum}. Throughout, we fix a choice of an integer $b\geq 3$ and $a_0\in\{0,1,\dotsc,b-1\}$, and we set
\[\mathcal{A}:=\bigg\{\sum_{j\geq 0}n_jb^j: n_j\in \{0, \dotsc, b-1\}\setminus \{a_0\}\: \forall j \bigg\}.\]
In addition, given an integer $r\in \mathcal{A}\cap [0, b)$, we let
\[\mathcal{A}_r=\{n\in \mathcal{A}: n\equiv r\imod b\}.\]
For $(c,d)=(r,b)=1$, we set
\begin{align}
\notag \mathcal{E}(X;d,c;b,r)=\sum_{\substack{n < X\\n\equiv c \imod d\\n\equiv r\imod b}}\Lambda(n)1_{\mathcal{A}}(n) - \dfrac{1}{\varphi(d)}\dfrac{b}{\varphi(b)}\sum_{\substack{n<X\\n\equiv r\imod b}}1_{\mathcal{A}}(n).
\end{align}
Furthermore, we set $\zeta := \frac{\log (b-1)}{\log b}$ for the rest of the paper.

\subsection*{Organization of the paper}
We will give a proof outline in Section \ref{Setup} following Maynard \cite{May2}, which is based on the circle method.

We devote Part \ref{part: sieve} to establish Theorem \ref{Main TheoremApp}. 


The graphical structure for the Sections \ref{prelim: sieve} and \ref{ThmappProof} can be described below:
\medskip

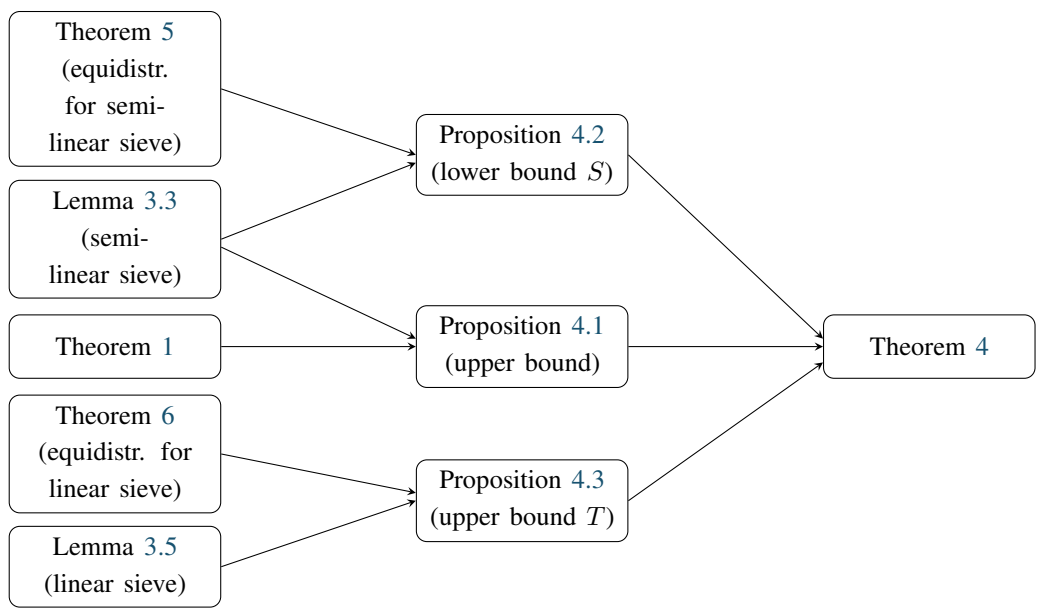
\begin{center}

\makebox[\textwidth]{\parbox{1.5\textwidth}{
			\begin{center}
				\tikzstyle{interface}=[draw, text width=6em,
				text centered, minimum height=2.0em]
				\tikzstyle{daemon}=[draw, text width=6em,
				minimum height=2em, text centered, rounded corners]
				\tikzstyle{lemma}=[draw, text width=6em,
				minimum height=2em, text centered, rounded corners]
				\tikzstyle{dots} = [above, text width=6em, text centered]
				\tikzstyle{wa} = [daemon, text width=6em,
				minimum height=2em, rounded corners]
				\tikzstyle{ur}=[draw, text centered, minimum height=0.01em]
				\def\blockdist{1.3}
				\def\edgedist{0.}
				\begin{tikzpicture}
		\node (thm6)[daemon]  {\footnotesize Proposition \ref{PropUpper} (upper bound)};
		\path (thm6.west)+(-4,0) node (thm2)[daemon] {\footnotesize Theorem \ref{Main Theorem}};
		\path (thm6.north) + (0,2) node (hyp2)[daemon] {\footnotesize Proposition \ref{PropAppS} (lower bound $S$)};
		\path (thm6.south) + (0,-1.5) node (appT)[daemon] {\footnotesize Proposition \ref{PropAppT} (upper bound $T$)};
	\path (thm6.east)+(4,0) node (thm3)[daemon] {\footnotesize Theorem \ref{Main TheoremApp}};	
\path (thm2.north)+(0,3) node (d1)[daemon] {\footnotesize Theorem \ref{Hyp1} (equidistr. for semi-linear sieve)};
\path (thm2.south)+(0, -1) node (d2)[daemon] {\footnotesize Theorem \ref{Hyp2} \\(equidistr. for linear sieve)};	
\path (thm2.north)+(0, 1) node (d3)[daemon] {\footnotesize Lemma \ref{fundsemi}\\ (semi-linear sieve)};	
\path (thm2.south)+(0, -2.5) node (d4)[daemon] {\footnotesize Lemma \ref{fundlin} (linear sieve)};
\path [draw, ->,>=stealth] (thm6.east) -- node [above] {} (thm3.west) ;
\path [draw, ->,>=stealth] (hyp2.east) -- node [above] {} ([yshift =3]thm3.west); 
\path [draw, ->,>=stealth] (appT.east) -- node [above] {} ([yshift=-6]thm3.west);
\path [draw, ->,>=stealth] (thm2.east) -- node [above] {} (thm6.west);
\path [draw, ->,>=stealth] (d1.east) -- node [above] {} (hyp2.west);
\path [draw, ->,>=stealth] (d2.east) -- node [above] {} ([yshift=3]appT.west);
\path [draw, ->,>=stealth] ([yshift=-3]d3.east) -- node [above] {} ([yshift=3]thm6.west);
\path [draw, ->,>=stealth] (d3.east) -- node [above] {} ([yshift =-3]hyp2.west);
\path [draw, ->,>=stealth] (d4.east) -- node [above] {} (appT.west);

	\end{tikzpicture}
	\end{center}}}
\end{center}

\medskip

In Part \ref{part: exp sum} we will establish exponential sums estimates over primes in arithmetic progressions, which is one of the key ingredients to prove our main results.

Finally, in Part \ref{part: circle} we will employ the circle method to establish our main theorems. In particular, we will deduce Theorems \ref{Main Theorem}, \ref{Main Theorem2}, \ref{Main Theorem3}, \ref{Hyp1} and \ref{Hyp2} from a more general theorem, Theorem \ref{Thm G} in Section \ref{mainproof}.

The dependency graph for Part \ref{part: circle} leading to the proofs of Theorems \ref{Main Theorem}, \ref{Main Theorem2}, \ref{Main Theorem3}, \ref{Hyp1}, and \ref{Hyp2} is given below:

\begin{center}

\makebox[\textwidth]{\parbox{1.5\textwidth}{
			\begin{center}
				\tikzstyle{interface}=[draw, text width=6em,
				text centered, minimum height=2.0em]
				\tikzstyle{daemon}=[draw, text width=6em,
				minimum height=2em, text centered, rounded corners]
				\tikzstyle{lemma}=[draw, text width=6em,
				minimum height=2em, text centered, rounded corners]
				\tikzstyle{dots} = [above, text width=6em, text centered]
				\tikzstyle{wa} = [daemon, text width=6em,
				minimum height=2em, rounded corners]
				\tikzstyle{ur}=[draw, text centered, minimum height=0.01em]
				\def\blockdist{1.3}
				\def\edgedist{0.}
				\begin{tikzpicture}
		\node (thm6)[daemon]  {\footnotesize Theorem \ref{Main Theorem3}};
		\path (thm6.north)+(0,1.6) node (thm1)[daemon] {\footnotesize Theorem \ref{Main Theorem}};
		\path (thm6.north)+(0,0.6) node (thm2)[daemon] {\footnotesize Theorem \ref{Main Theorem2}};

		\path (thm6.south) + (0,-0.6) node (hyp1)[daemon] {\footnotesize Theorem \ref{Hyp1}};
		\path (thm6.south) + (0, -1.6) node (hyp2)[daemon] {\footnotesize Theorem \ref{Hyp2}};

	\path (thm6.west)+(-7,1) node (propM)[daemon] {\footnotesize Proposition \ref{PropMajArc}\\ (major arcs estimate)};
	\path (thm6.west)+(-7,-1) node (propm)[daemon] {\footnotesize Proposition \ref{PropMinorArc1}\\ (minor arcs estimate)};
	\path (thm6.west)+(-11,-1) node (hyb)[daemon] {\footnotesize Lemma \ref{lemd2} (hybrid bound)};
	\path (thm6.west)+(-11,1) node (linf)[daemon] {\footnotesize Lemma \ref{lemd3} ($L^\infty$ bound)};

\path (thm6.west)+(-3,0) node (thmG)[daemon] {\footnotesize Theorem \ref{Thm G}\\ (gen thm)};		
\path (thmG.north)+(0,1.7) node (d1)[daemon] {\footnotesize Proposition \ref{PropminorPrimeI} };
\path (thmG.north)+(0,0.7) node (d2)[daemon] {\footnotesize Proposition \ref{PropminorPrimeII}};	
\path (thmG.south)+(0,-0.7) node (d3)[daemon] {\footnotesize Proposition \ref{expwellfac}};	
\path (thmG.south)+(0,-1.7) node (d4)[daemon] {\footnotesize Proposition \ref{expsemi}};
\path (thmG.south)+(0,-2.7) node (d5)[daemon] {\footnotesize Proposition \ref{explin}};
\path [draw, ->,>=stealth] (d1.east) -- node [above] {} (thm1.west) ;
\path [draw, ->,>=stealth] (d2.east) -- node [above] {} (thm2.west) ;
\path [draw, ->,>=stealth] (d3.east) -- node [above] {} ([yshift=-3]thm6.west) ;
\path [draw, ->,>=stealth] (d4.east) -- node [above] {} ([yshift=-3]hyp1.west) ;
\path [draw, ->,>=stealth] (d5.east) -- node [above] {} (hyp2.west) ;
\path [draw, ->,>=stealth] ([yshift=6]thmG.east) -- node [above] {} ([yshift=-2]thm1.west) ;
\path [draw, ->,>=stealth] ([yshift=3]thmG.east) -- node [above] {} ([yshift=-2]thm2.west) ;
\path [draw, ->,>=stealth] (thmG.east) -- node [above] {} (thm6.west) ;
\path [draw, ->,>=stealth] ([yshift=-3]thmG.east) -- node [above] {} ([yshift=3]hyp1.west) ;
\path [draw, ->,>=stealth] ([yshift=-6]thmG.east) -- node [above] {} ([yshift=3]hyp2.west) ;
\path [draw, ->,>=stealth] (propM.east) -- node [above] {} (thmG.west) ;
\path [draw, ->,>=stealth] (propm.east) -- node [above] {} ([yshift=-2]thmG.west) ;
\path [draw, ->,>=stealth] (hyb.east) -- node [above] {} (propm.west);
\path [draw, ->,>=stealth] (linf.east) -- node [above] {} (propM.west);

	\end{tikzpicture}
	\end{center}}}
\end{center}

\section{Set-up and outline of the proof}\label{Setup}
The strategy to prove Theorem \ref{Main Theorem}, Theorem \ref{Main Theorem2}, and Theorem \ref{Main Theorem3} is to apply the circle method. For the sake of exposition, we will outline the proof of Theorem \ref{Main Theorem} following the set-up from Maynard \cite{May2}.

Let $\widehat{{1}}_{\mathcal{A}}$ be the Fourier transform of the set $\mathcal{A}$ restricted to $\{1, \dotsc, X\}$ with $X=b^k$. Then, for any real number $\theta\in [0, 1)$, we have
\begin{align}\label{Afourier}
\widehat{{1}}_{\mathcal{A}}(\theta): =\sum_{n<X}{1}_{\mathcal{A}}(n)e(n\theta) = \prod_{j=0}^{k-1}\bigg(\sum_{0\leq n_j<b}{1}_{\mathcal{A}}(n_i)e(n_jb^j \theta)\bigg),
\end{align}
where $n=\sum_{j=0}^{k-1}n_jb^j$. Next, for $r\in \mathcal{A}$, we set
\begin{align}
\notag \mathcal{A}_r=\{n\in \mathcal{A}:n\equiv r\imod b\}.
\end{align}
We then define 
\begin{align}
\label{Arfourier}\widehat{1}_{\mathcal{A}_r}(\theta)&:=\sum_{n<X}{1}_{\mathcal{A}_r}(n)e(n\theta)\\
\notag &=e(r\theta)\prod_{j=1}^{k-1}\bigg(\sum_{0\leq n_j<b}{1}_{\mathcal{A}}(n_j)e(n_jb^j \theta)\bigg).
\end{align}
Note that for $r\in \mathcal{A} \cap[0, b)$ and for any real number $\theta\in [0, 1)$, we have the trivial bound:
\[\big|\widehat{1}_{\mathcal{A}_r}(\theta)\big| \leq \dfrac{X^\zeta}{b-1} \leq X^\zeta,\]
which we will often use in the paper.

Next, by Fourier inversion on $\mathbb{Z}/X\mathbb{Z}$, for $n<X$, we have
\begin{align}
\label{fourinv} {1}_{\mathcal{A}_r}(n) = \dfrac{1}{X}\sum_{0\leq t < X}\widehat{{1}}_{\mathcal{A}_r}\bigg(\frac{t}{X}\bigg)e\bigg(\frac{-nt}{X}\bigg).
\end{align}

In order to prove Theorem \ref{Main Theorem}, we consider the following setup. For $(c, d)=1$ and for any real number $\theta\in [0, 1)$, we set  
\begin{align}\label{lambdafour}
\widehat{\Lambda}_{d,\: c}(\theta)=\sum_{\substack{n<X\\n\equiv c\imod {d}}}\Lambda(n)e(n\theta).
\end{align}
Then, by the relations \eqref{fourinv} and \eqref{lambdafour}, we have
\begin{align}
\sum_{\substack{n<X\\n\equiv c\imod {d}\\n\equiv r\imod b}}\Lambda(n)1_{\mathcal{A}}(n) = \dfrac{1}{X}\sum_{0\leq t<X}\widehat{1}_{\mathcal{A}_r}\bigg(\dfrac{t}{X}\bigg)\widehat{\Lambda}_{d, c}\bigg(\dfrac{-t}{X}\bigg).
\end{align}
Therefore, our task in \eqref{MT1eq1} reduces to showing that
\begin{align}
\label{Mainsum}\sum_{\substack{d\leq D\\(d, b)=1}}\max_{\substack{1\leq c<d\\(c, d)=1}}\bigg|\dfrac{1}{X}\sum_{0\leq t<X}\widehat{1}_{\mathcal{A}_r}\bigg(\dfrac{t}{X}\bigg)\widehat{\Lambda}_{d, c}\bigg(\dfrac{-t}{X}\bigg) - \dfrac{1}{\varphi(d)}\dfrac{b}{\varphi(b)}\sum_{n<X}1_{\mathcal{A}_r}(n)\bigg|\ll_{A, b, \delta}\dfrac{X^\zeta}{(\log X)^A}.
\end{align}
We then consider two cases according to whether $t/X$ is close to a rational number with a small denominator or not, namely, major arcs and minor arcs, respectively.

\medskip

\noindent {\it Major arcs:} The \emph{major arcs} $\mathfrak{M}$ are those $t$'s in $[0, X)\cap \mathbb{Z}$ such that 
\[\bigg|\dfrac{t}{X}-\dfrac{a}{q}\bigg|\leq \dfrac{(\log X)^C}{X}\]
for some $(a, q)=1$, $0\leq a<q$, $1\leq q\leq (\log X)^C$ with $C>0$ to be chosen later in terms of $A$. It will be convenient to divide the major arcs $\mathfrak{M}$ into three disjoint subsets:
\begin{align}
\label{maj}\mathfrak{M}=\mathfrak{M}_1\cup\mathfrak{M}_2\cup \mathfrak{M}_3,
\end{align}
where 
\begin{align}
\notag \mathfrak{M}_1 =\bigg\{t& \in [0, X)\cap \mathbb{Z}: \bigg\lvert\dfrac{t}{X}-\dfrac{a}{q}\bigg\rvert \leq \dfrac{(\log X)^C}{X}~\text{for some}\, (a, q)=1, 1\leq a <q\leq (\log X)^C, q\nmid X\bigg\},\\
\notag \mathfrak{M}_2=\bigg\{t& \in [0, X)\cap \mathbb{Z}: \dfrac{t}{X}=\dfrac{a}{q}+ \dfrac{\eta}{X}\\
\notag & \text{for some}\, (a, q)=1, 0\leq a< q\leq (\log X)^C, q\geq 1, q\vert X, 0<\lvert \eta \rvert \leq (\log X)^C\bigg\},\\
\notag \mathfrak{M}_3 =\bigg\{t& \in [0, X)\cap \mathbb{Z}: \dfrac{t}{X}=\dfrac{a}{q}\: \text{for some}\, (a, q)=1, 0\leq a< q\leq (\log X)^C, q\geq 1,  q\lvert X\bigg\}.
\end{align}
We now briefly explain how we will estimate the sum \eqref{Mainsum} when $t$ is in one of the above-defined three sets of the \emph{major arcs}.
\begin{enumerate}[(a)]
\item  We use the $L^\infty$ bound for the Fourier transform of the set $\mathcal{A}_r$ and the trivial bound for $\widehat{\Lambda}_{d,c}(t/X)$ to estimate the sum \eqref{Mainsum} when $t\in \mathfrak{M}_1$.
\item It turns out that when $t\in \mathfrak{M}_2$, we can use the Bombieri-Vinogradov Theorem to handle the exponential sum $\widehat{\Lambda}_{d,c}(t/X)$ and the trivial bound for $\widehat{1}_{\mathcal{A}_r}(t/X)$ in \eqref{Mainsum}. 
\item  When $t\in \mathfrak{M}_3$, we get the main term in $\eqref{Mainsum}$ and the error term is again controlled by using the Bombieri-Vinogradov Theorem.
We note that $\widehat{1}_{\mathcal{A}_r}(t/X)$ is large if $t$ is close to a number with few non-zero base-$b$ digits. 
\end{enumerate}
 
This will establish our estimate in \eqref{Mainsum} when $t/X$ is in \emph{major arcs}.

\medskip

\noindent{\it Minor arcs:} The `\emph{minor arcs}' $\mathfrak{m}$ are those $t \in [0,X)\cap\mathbb{Z}$ such that $t\not\in\mathfrak{M}$. We use a $L^\infty-L^1$ bound to handle minor arcs as follows:
\begin{align}
\label{minorD}
\begin{split}
\sum_{\substack{d\leq D\\(d, b)=1}}\max_{(c, d)=1} \bigg| \frac{1}{X}\sum_{t\in\mathfrak{m}} \widehat{1}_{\mathcal{A}_r}&\bigg(\dfrac{t}{X}\bigg)\widehat{\Lambda}_{d,c}\bigg(\dfrac{-t}{X}\bigg)\bigg|\\
& \le \bigg(\sup_{t \in \mathfrak{m}}\sum_{\substack{d\leq D\\(d, b)=1}}\max_{(c, d)=1}\bigg|\widehat{\Lambda}_{d,c}\bigg(\frac{-t}{X}\bigg)\bigg|\bigg)\sum_{t\in \mathfrak{m}}\dfrac{1}{X}\bigg|\widehat{1}_{\mathcal{A}_r}\bigg(\frac{t}{X}\bigg)\bigg|.
\end{split}
\end{align}
As in Maynard \cite{May2}, we use a large-sieve type argument to control the $L^1$ sum of $\widehat{1}_{\mathcal{A}_r}$, which is shown to be small in Lemma \ref{lemd2}. Next, our goal is to save over the trivial bound on $\sum_{d\leq D}\max_{(c, d)=1}\Big\lvert\widehat{\Lambda}_{d,c}(t/X)\Big\rvert$ when $t \in \mathfrak{m}$ and $D$ as large as possible. We use estimates from exponential sum over primes in arithmetic progressions from the works of Matom\"aki \cite{Mat}, Mikawa \cite{Mik}, and Ter\"{a}v\"{a}inen \cite{Ter} to handle those sums over primes in Part \ref{part: exp sum}. Combining these $L^1$ and $L^\infty$ bounds, we will show that
\begin{align}
\notag \bigg(\sup_{t \in \mathfrak{m}}\sum_{\substack{d\leq D\\(d, b)=1}}\max_{(c, d)=1}\bigg|\widehat{\Lambda}_{d,c}\bigg(\frac{-t}{X}\bigg)\bigg|\bigg)\sum_{t\in \mathfrak{m}}\dfrac{1}{X}\bigg|\widehat{1}_{\mathcal{A}_r}\bigg(\frac{t}{X}\bigg)\bigg| \ll_{A,b, \delta}\dfrac{X^{\zeta}}{(\log X)^{A}}.
\end{align} 
This completes the rough outline of the proof of Theorem \ref{Main Theorem}.

The key difference in the proofs of Theorems \ref{Main Theorem2} and \ref{Main Theorem3} compared to Theorem \ref{Main Theorem} is better exponential sums estimate over primes in arithmetic progressions, which allows us to take a bigger range of the moduli $d\leq D$.

\begin{remark}
Note that we will establish a much more general theorem, Theorem \ref{Thm G}, for an arithmetic function $\mathfrak{f}$ satisfying some appropriate conditions in Part \ref{part: circle}. In particular, Theorem \ref{Thm G} will incorporate Theorems \ref{Main Theorem}, \ref{Main Theorem2}, and \ref{Main Theorem3} by choosing $\mathfrak{f}$ and other parameters appropriately.
\end{remark}

\part{Sieve methods and their applications}\label{part: sieve}

\section{Preliminaries from sieve methods}\label{prelim: sieve}
In this section, we collect some technical results from sieve methods that will be needed to prove Theorem \ref{Main TheoremApp}. 


Given a sequence of weights $\mathcal{C}=\big(c(n)\big)_{n=1}^\infty \subset \mathbb{R}_{\geq 0}$ with $\sum_{n=1}^\infty c(n) <\infty$ and a set of primes $\mathcal{P}$, we consider the {\it sifting function},
\[S(\mathcal{C},\mathcal{P}, z) := \sum_{(n, P(z))=1}c(n),\]
where for some real number $z>1$,
\[P(z): = \prod_{\substack{p \leq z\\p \in \mathcal{P}}}p.\]
Here $z$ is often called the \emph{sifting parameter} in the sieve setting.

In order to proceed further, for any $x\geq 1$ and for each integer $d\geq 1$, we set 
\[C_d(x) := \sum_{\substack{n \leq x\\d|n}}c(n), \]
and we impose the following {\it axioms of sieve theory}:
\begin{enumerate}[label=(A\arabic*)]
\item \label{A1} For some multiplicative function $\mathfrak{g}$, we have \[C_d(x)= \dfrac{\mathfrak{g}(d)}{d}C(x) + {E}(d),\]
where $C(x)$  can be interpreted as an approximation to $\sum_{n\leq x}c(n)$ and ${E}(d)$ is a real number which we think of as an error term.
\item \label{A2} We assume that the multiplicative function $\mathfrak{g}$ satisfies $\mathfrak{g}(p) \leq \min\{2, p-1\}$ for all primes $p\in \mathcal{P}$.
\item \label{A3} There is a constant $A>0$, and a quantity $D\geq 1$ such that
\[\sum_{d\leq D}\mu^2(d)|{E}(d)| \ll_{A}\dfrac{C(x)}{(\log x)^A}.\]
If such an estimate holds, then we say $\mathcal{C}$ has {\it level of distribution} $D$.
\item \label{A4} We have
\[\sum_{\substack{p \leq x\\ p \in \mathcal{P}}}\dfrac{\mathfrak{g}(p)\log p}{p} = \varkappa\log x + O(1) \quad \text{for all $x$}.\]
Here we say $\varkappa$ as the {\it dimension of the sieve}.
\end{enumerate}

Next, we state the definition of what is an {\it upper bound sieve} and a {\it lower bound sieve}. 

\begin{definition}[Upper bound sieve]
An arithmetic function $\lambda^+:\mathbb{N}\rightarrow\mathbb{R}$ that is supported on the set $\{d|P(z): \: d \leq D\}$ and satisfies the relation $(\lambda^+*1)(n) \geq 1_{(n,\: P(z))=1}$ is called an \emph{upper bound sieve of level $D$} for the set of primes $\mathcal{P}$.
\end{definition}

\begin{definition}[Lower bound sieve]
An arithmetic function $\lambda^-:\mathbb{N}\rightarrow\mathbb{R}$ that is supported on the set $\{d|P(z): \: d \leq D\}$ and satisfies the relation $1_{(n,\:P(z))=1} \geq (\lambda^-*1)(n)$ is called a \emph{lower bound sieve of level $D$} for the set of primes $\mathcal{P}$.
\end{definition}

\begin{remark}
We will refer to $\lambda^{\pm}$ as the \emph{sieve weights} or \emph{sifting weights} in this paper.
\end{remark}

Now we are ready to state the Fundamental Lemma of Sieve Theory in the special case when the dimension $\varkappa$ equals $1/2$, often referred to as the {\it semi-linear sieve} or the {\it half-dimensional sieve}.

\begin{lemma}[Fundamental Lemma for the Semi-linear Sieve]\label{fundsemi}
Consider a sequence $\mathcal{C}=\big(c(n)\big)_{n=1}^\infty$ of non-negative real numbers and a set of primes $\mathcal{P}$ satisfying axioms \ref{A1}, \ref{A2}, and \ref{A4} with $\varkappa = 1/2$. If $u_1>0$ and $D=z_1^{u_1}$, then there exist two arithmetic functions $\lambda^\pm_{\textup{sem}}: \mathbb{N}\rightarrow [-1, 1]$ supported on the set $\{d|P(z_1): d\leq D\}$, and we have 
\begin{equation}\label{lower}
S(\mathcal{C}, \mathcal{P}, z_1) \geq C(x)\bigg\{f_{\textup{sem}}(u_1) + o(1)\bigg\}\prod_{p\leq z_1,\: p\in \mathcal{P}}\bigg(1-\dfrac{\mathfrak{g}(p)}{p}\bigg) - \sum_{\substack{p|d \Rightarrow p \in \mathcal{P}\\ d\leq D}}\lambda^-_{\textup{sem}}(d){E}(d),
\end{equation}
and
\begin{equation}
\label{upper} S(\mathcal{C}, \mathcal{P}, z_1) \leq C(x)\bigg\{F_{\textup{sem}}(u_1) + o(1)\bigg\}\prod_{p\leq z_1,\: p\in \mathcal{P}}\bigg(1-\dfrac{\mathfrak{g}(p)}{p}\bigg)+\sum_{\substack{p|d \Rightarrow p \in \mathcal{P}\\ d\leq D}}\lambda^+_{\textup{sem}}(d){E}(d),
\end{equation}
where $f_{\textup{sem}}, F_{\textup{sem}}$ are continuous functions in $u_1 = \log D/\log z_1$ such that 
\begin{equation}
\begin{cases} \label{s}
\sqrt{u_1}F_{\textup{sem}}(u_1)=2\sqrt{e^\gamma/\pi} & \text{if}\quad 0 <u_1 \leq 2, \\ 
f_{\textup{sem}}(u_1) = 0 & \text{if}\quad 0<u_1\leq 1, 
\end{cases}
\end{equation}
where $\gamma$ is the Euler-Mascheroni constant, and 
for $1\leq u_1\leq 3$ we have
\begin{equation}\label{fsem}
\dfrac{\sqrt{u_1}f_{\textup{sem}}(u_1)}{\sqrt{e^\gamma/\pi}}=\int_1^{u_1} \dfrac{\textup{d}y}{\sqrt{y(y-1)}}=\log \bigg(1+2(u_1-1)+2\sqrt{u_1(u_1-1)}\bigg).
\end{equation}
\end{lemma}
\begin{proof}
The proof follows from \cite[Theorem 11.12--Theorem 11.13 ]{opera}) with $\beta=1$ and \cite[Chapter 14 (pp. 275--276)]{opera}.
\end{proof}

We also state the partial well-factorability (see Definition \ref{def: well fac}) of the semi-linear sieve in the next lemma.

\begin{lemma}[Partial well-factorability of semi-linear sieve]\label{wellfacsem}
Let $\varepsilon>0$ be small. Let $\delta\in (0, 10^{-3}]$ and let $\rho_{\textup{sem}}=\frac{3}{7}(1-4\delta)-\varepsilon$. Then the lower bound semi-linear sieve weights $\lambda_{\textup{sem}}^-$ as given in Lemma \ref{fundsemi} with level $X^{\rho_{\textup{sem}}}$ and sifting parameter $z_1\leq X^{1/3-2\delta-2\varepsilon^2}$ is supported in the set 
\begin{align}
\mathfrak{D}^{-,\: \textup{sem}}=\{p_1\cdots p_r\leq X^{\rho_{\textup{sem}}}: z_1\geq p_1> \dotsc >p_r, p_1\cdots p_{2m-1}p_{2m}^2\leq X^{\rho_{\textup{sem}}}\: \forall \: m \geq 1\},
\end{align}
where $p_1, \dotsc, p_r$ denote primes. In addition, for any $D_0\in [X^{1/3-2\delta-2\varepsilon^2}, X^{\rho_{\textup{sem}}}]$, every $d\in \mathfrak{D}^{-, \textup{sem}}\cap[X^{1/10}, X^{\rho_{\textup{sem}}}]$ can be factorized as $d=d_1d_2$ such that $d_1\in [X^{1/10}, D_0]$ and $d_1d_2^2\leq X^{1-4\delta-2\varepsilon^2}/D_0$.
\end{lemma}
\begin{proof}
This is \cite[Lemma 9.2]{Ter} with $\theta=\delta$ and $D=D_0$. 
\end{proof}

Next, we state the Fundamental Lemma for the {\it linear sieve}, that is, for dimension $\varkappa=1$.
\begin{lemma}[Fundamental Lemma for the Linear Sieve]\label{fundlin}
Consider a sequence $\mathcal{C} = \big(c(n)\big)_{n=1}^\infty\subset\mathbb{R}_{\geq 0}$ and a set of primes $\mathcal{P}$ satisfying axioms \ref{A1}, \ref{A2}, and \ref{A4} with $\varkappa = 1$. If $u_2>0$ and $D=z_2^{u_2}$, then there exist two arithmetic functions $\lambda^\pm_{\textup{lin}}: \mathbb{N}\rightarrow [-1, 1]$ supported on the set $\{d|P(z_2): d\leq D\}$, and we have
\begin{equation}\label{lowerl}
 S(\mathcal{C}, \mathcal{P}, z_2) \geq C(x)\Big\{f_{\textup{lin}}(u_2) + o(1)\Big\}\prod_{p\leq z_2,\: p\in \mathcal{P}}\bigg(1-\dfrac{\mathfrak{g}(p)}{p}\bigg) - \sum_{\substack{p|d \Rightarrow p \in \mathcal{P}\\ d\leq D}}\lambda^-_{\textup{lin}}(d){E}(d),
\end{equation}
and
\begin{equation}\label{upperl}
S(\mathcal{C}, \mathcal{P}, z_2) \leq C(x)\Big\{F_{\textup{lin}}(u_2) + o(1)\Big\}\prod_{p\leq z_2,\: p\in \mathcal{P}}\bigg(1-\dfrac{\mathfrak{g}(p)}{p}\bigg) +\sum_{\substack{p|d \Rightarrow p \in \mathcal{P}\\ d\leq D}}\lambda^+_{\textup{lin}}(d){E}(d),
\end{equation}
where $f_{\textup{lin}}, F_{\textup{lin}}$ are continuous functions in $u_2 = \log D/\log z_2$ such that 
\begin{equation}
\begin{cases} \label{sl}
u_2 F_{\textup{lin}}(u_2)=2e^\gamma & \text{if}\quad 1 \leq u_2 \leq 3, \\ 
f_{\textup{lin}}(u_2) = 0 & \text{if}\quad 0<u_2\leq 2,
\end{cases}
\end{equation}
where $\gamma$ is the Euler-Mascheroni constant.

\end{lemma}

\begin{proof}
The proof follows from \cite[Theorem 11.12--Theorem 11.13 ]{opera}) with $\beta=2$ and \cite[Chapter 12 (pp. 235--236)]{opera}.
\end{proof}

In order to deal with the linear sieve in Theorem \ref{Hyp2}, we need the following well-factorability lemma.

\begin{lemma}[Well-factorability of linear sieve]\label{wellfacclin}
Let $\varepsilon>0$ be small. Let $\delta\in (0, 10^{-3}]$  and let $\rho_{\textup{lin}}=\frac{1}{2}-2\delta-\varepsilon$. Then the upper bound linear sieve weights $\lambda_{\textup{lin}}^+$ as given in Lemma \ref{fundlin} with level $X^{\rho_{\textup{lin}}}$ and sifting parameter $z_2\leq X^{1/2}$ is supported in the set
\begin{align}
\mathfrak{D}^{+,\: \textup{lin}}=\{p_1\cdots p_r\leq X^{\rho_{\textup{lin}}}: z_2\geq p_1> \dotsc >p_r, p_1\cdots p_{2m-2}p_{2m-1}^3\leq X^{\rho_{\textup{lin}}}\: \forall \: m\geq 1\},
\end{align}
where $p_1, \dotsc, p_r$ denote primes. 
In addition, for any $D_0 \in [X^{1/5}, X^{\rho_{\textup{lin}}}]$,  every $d \in \mathfrak{D}^{+, \textup{lin}}\cap[X^{1/10}, X^{\rho_{\textup{lin}}}]$ can be factorized as $d=d_1d_2$ such that $d_1\in [X^{1/10}, D_0]$ and $d_1d_2^2\leq X^{1-4\delta-2\varepsilon^2}/D_0$. 
\end{lemma}
\begin{proof}
See \cite[Lemma 9.1]{Ter} or \cite[Lemma 12.16]{opera}.
\end{proof}

\section{Proof of Theorem \ref{Main TheoremApp}}\label{ThmappProof}

\subsection{Upper bound in Theorem \ref{Main TheoremApp}}
We first establish the upper bound in Theorem \ref{Main TheoremApp} by using Lemma \ref{fundsemi} and assuming Theorem \ref{Main Theorem}.
\begin{prop}\label{PropUpper}
Let $b$  be a sufficiently large odd integer and $r \in \mathcal{A}\cap[0, b)$ be such that $(r, b)=(r-1, b)=1$. Then, we have
\begin{align}
\notag \sum_{\substack{p<X\\p\equiv r\imod b}}1_{\mathcal{A}}(p)1_{\mathbb{B}}(p-1)\ll_{b} \dfrac{X^{\zeta}}{(\log X)^{3/2}}.
\end{align}
\end{prop}

\begin{proof}
Let $z\in [2, X]$ be a parameter to be chosen later. We let
\begin{align}
\notag \mathcal{P}_3=\{p\equiv 3\imod 4, p \nmid b\} \quad \text{and} \quad P_3(z) = \prod_{\substack{p\leq z\\p\in \mathcal{P}_3}}p.
\end{align}
Then, we have
\begin{align}
\label{upper1}\sum_{\substack{p<X\\p\equiv r\imod b}}1_{\mathcal{A}}(p) 1_{\mathbb{B}}(p-1) &\leq \sum_{\substack{p<X\\(p-1, P_3(z))=1}}1_{\mathcal{A}_r}(p)\\
\notag  &< \dfrac{10}{9\log X}\sum_{\substack{X^{9/10}<p<X\\(p-1, P_3(z))=1}}1_{\mathcal{A}_r}(p)\log p + X^{9/10}\\ 
\label{upper2} &< \dfrac{10}{9\log X}\sum_{\substack{n<X\\(n-1, P_3(z))=1}}\Lambda(n)1_{\mathcal{A}_r}(n) + X^{9/10}.
\end{align}
Next, for $d|P_3(z)$, we set
\begin{align}
\notag {E}(d)=\sum_{\substack{n<X\\n\equiv 1\imod d\\}}\Lambda(n)1_{\mathcal{A}_r}(n)- \dfrac{1}{\varphi(d)}\dfrac{b}{\varphi(b)}\sum_{\substack{n<X}}1_{\mathcal{A}_r}(n).
\end{align}
Therefore, by Theorem \ref{Main Theorem} with $D=X^{3/10}$, for any large real number $A>0$, we find that
\begin{align}
\notag \sum_{\substack{d\leq X^{3/10}\\d|P_3(z)}}\lvert {E}(d)\rvert \ll \dfrac{X^{\zeta}}{(\log X)^A}.
\end{align}
Now, we choose  $c(n)=\Lambda(n)1_{\mathcal{A}_r}(n)$  for $n<X$ and $z_1=D=z=X^{3/10}$ in Lemma \ref{fundsemi}. Clearly, the sequence $c(n)$ satisfies the axioms of sieve theory with $\mathfrak{g}(d)=d/\varphi(d)$ and $\varkappa =1/2$. Therefore, by the upper bound semi-linear sieve \eqref{upper} with $u_1=1$, we have 
\begin{align}
\notag \sum_{\substack{n<X\\(n-1, P_3(z))=1}}\Lambda(n)1_{\mathcal{A}_r}(n)\leq&~\bigg(\dfrac{2e^{\gamma/2}}{\pi^{1/2}} + o(1)\bigg) \dfrac{b}{\varphi(b)}\sum_{\substack{n<X}}1_{\mathcal{A}_r}(n)\prod_{\substack{p<z\\p\in \mathcal{P}_3}}\bigg(1-\dfrac{1}{p-1}\bigg)+ O\bigg(\dfrac{X^{\zeta}}{(\log X)^A}\bigg),
\end{align}
where $\gamma$ is the Euler-Mascheroni constant. The above estimate together with the estimates from \eqref{upper1} and \eqref{upper2} allows us to obtain
\begin{align}
\notag \sum_{\substack{p<X\\p\equiv r\imod b}}1_{\mathcal{A}}(p)1_{\mathbb{B}}(p-1) \leq &~\bigg(\dfrac{2e^{\gamma/2}}{\pi^{1/2}} + o(1)\bigg) \dfrac{10b}{9\varphi(b)\log X}\sum_{\substack{n<X}}1_{\mathcal{A}_r}(n)\prod_{\substack{p<z\\p\in \mathcal{P}_3}}\bigg(1-\dfrac{1}{p-1}\bigg)\\
\notag &+ O\bigg(\dfrac{X^{\zeta}}{(\log X)^A} + X^{9/10}\bigg).
\end{align}
We can now use Mertens' estimate \cite[Theorem 3.4(c)]{kou} to the product over the primes (for example, see \cite[p.278]{opera} for a detailed estimate) and the fact that $\sum_{\substack{n<X}}1_{\mathcal{A}_r}(n)=X^\zeta/(b-1)$ to deduce that 
\begin{align}
\notag \sum_{\substack{p<X\\p\equiv r\imod b}}1_{\mathcal{A}}(p) 1_{\mathbb{B}}(p-1) \ll_{b} \dfrac{X^{\zeta}}{(\log X)^{3/2}}
\end{align}
as desired.
\end{proof}

\subsection{Lower bound in Theorem \ref{Main TheoremApp}}
The lower bound in Theorem \ref{Main TheoremApp} can also be obtained from \cite[Theorem 6.5]{Ter} by choosing $\omega_n=1_{\mathcal{A}}(n)\cdot 1_{n\equiv r\imod b}$, where Hypothesis 6.4 holds by considering variants of Theorems \ref{Hyp1} and \ref{Hyp2}. For the sake of completeness, we will establish the lower bound from scratch in this paper. In order to do so, we consider the following sieve setup.

\subsubsection{Sieve set-up for the lower bound.} For $r\in \mathcal{A}\cap[0,b)$ with $\big(r(r-1), b\big)=1$, we set
 \begin{align}\label{setup F}
\begin{split}
&\mathcal{F} ={\{p-1:\: p < X, \: p\in \mathcal{A}_r,\: p\equiv 3\imod 8\}},\\
\mathcal{P}_3 & =\{p\equiv 3\imod 4, p\nmid b\},\quad \text{and} \quad  P_3(z)=\prod_{\substack{p<z\\p\in \mathcal{P}_3}}p.
\end{split}
\end{align}
Note that, since $p\equiv r\imod b$ for the primes we are considering here, and we have assumed that $(r-1,b)=1$, so there are no primes that divide both $p-1$ and $b$. So, we have that 
\begin{align}
\label{sieve lower} \sum_{\substack{p<X\\p\equiv r\imod b}}1_{\mathcal{A}}(p)1_{\mathbb{B}}(p-1)\geq S(\mathcal{F}, \mathcal{P}_3, X^{1/2})=\sum_{\substack{p<X\\(p-1,\: P_3(X^{1/2}))=1\\p\equiv 3\imod 8}}1_{\mathcal{A}_r}(p).
\end{align}
For notational convenience, we set $z=X^{1/\alpha}$ for some $\alpha\in [2, 4)$. Later, we will choose $\alpha\approx3$.

By the Buchstab identity (see \cite[eqn (6.4)]{opera}), we have
\begin{align}
\label{Buch} S(\mathcal{F}, \mathcal{P}_3, \sqrt{X}) = S(\mathcal{F}, \mathcal{P}_3, z)-\sum_{\substack{z<p_1\leq \sqrt{X}\\p_1\equiv 3\imod 4}}S(\mathcal{F}_{p_1}, \mathcal{P}_3, p_1)=:S-T.
\end{align}
We will give a lower bound for $S$ using the semi-linear sieve and Theorem \ref{Hyp1}. On the other hand, an upper bound for $T$ is given using the linear sieve and Theorem \ref{Hyp2}. 

Since $p-1$ has an even number of prime factors in the class $3\imod 4$ and by our choice of $z$, we can write the sum $T$ as 
\begin{align}
\notag T=\sum_{p<X}\sum_{\substack{p-1=2n_1p_1p_2\\p_1,\: p_2\in \mathcal{P}_3\\p_2\geq p_1\geq X^{1/\alpha}\\n_1\in \mathcal{B}}}1_{\mathcal{A}_r}(p),
\end{align}
where $\mathcal{B}=\{n: p|n \Rightarrow p\equiv 1\imod 4\}$.
Following Matom\"aki \cite{Mat}, we define
\begin{align}
\label{mathcalL} \mathcal{L}&=\{\ell  =n_1p_1: n_1\leq X^{1-2/\alpha}, n_1\in \mathcal{B}, X^{1/\alpha}\leq p_1 <(X/n_1)^{1/2}, p_1\in \mathcal{P}_3\},\\
\intertext{and for each $\ell \in \mathcal{L}$,}
\label{mathcalMl} \mathcal{M}(\ell)& ={\{m=2\ell p_2+1:\: m \in \mathcal{A}_r,\: p_2< X/2\ell, \: p_2\in \mathcal{P}_3,\: p_2\geq X^{1/\alpha}\}}.
\end{align}
Note that for each $m\in \mathcal{M}(\ell)$, we have $m\equiv r\imod b$. Since, by our assumption $(r-1, b)=1$, we have that $(\ell, b)=1$.
This allows us to bound the sum $T$ as
\begin{align}
\label{Tsimple} T\leq \sum_{\substack{\ell \in \mathcal{L}\\(\ell, b)=1}}\bigg(S(\mathcal{M}(\ell), \mathcal{P}(\ell), X^{1/\nu}) + O(X^{1/\nu})\bigg),
\end{align}
where $\mathcal{P}(\ell)=\{p: p\nmid 2b\ell\}$ and we will choose $\nu$ appropriately later. In fact, we will choose $\nu \approx 5$.

\begin{remark}
Note that if $m\in \mathcal{M}(\ell)$ in \eqref{mathcalMl}, we have $2\ell p_2+1\equiv r\imod b$. Since $(r-1, b)=1$ this implies that $(2\ell p_2, b)=1$, which in turn restricts the base $b$ to be odd.
\end{remark}

Now we are ready to bound the sums $S$ from below and $T$ from above separately in the following two propositions.

\begin{prop}\label{PropAppS} Assume the above sieve set-up. Let $\varepsilon>0$ be small. Let $\delta\in (0, 10^{-3}]$ and let $b$ be an odd integer that is sufficiently large in terms of $\delta$. Let $\alpha=(1/3-2\delta)^{-1}+\varepsilon$ be such that $\alpha\in [2, 4)$ and let $\rho_{\textup{sem}}\leq \frac{3}{7}(1-4\delta)-\varepsilon$. Then we have 
\begin{align}
\label{EqS1} S\geq \dfrac{\mathfrak{S}+ o(1)}{(\log X)^{3/2}}\dfrac{b}{\varphi(b)}\prod_{\substack{p|b\\p\equiv 3\imod 4}}\bigg(1-\dfrac{1}{p-1}\bigg)^{-1}I_{\textup{sem}}(\rho_{\textup{sem}}, \alpha)\sum_{\substack{n<X}}1_{\mathcal{A}_r}(n),
\end{align}
where
\begin{align}
\label{ac} \mathfrak{S} &=\dfrac{1}{4\sqrt{2}}\prod_{p\equiv 3\imod 4}\bigg(1-\dfrac{1}{p^2}\bigg)^{1/2}\prod_{p\equiv 3\imod 4}\bigg(1-\dfrac{1}{(p-1)^2}\bigg)
\end{align}
and
\begin{align}
\label{eq: isem} I_{\textup{sem}}(\rho_{\textup{sem}}, \alpha)&=\dfrac{1}{\sqrt{\rho_{\textup{sem}}}}\int_{1}^{\alpha \rho_{\textup{sem}}} \dfrac{\textup{d}y}{\sqrt{y(y-1)}}.
\end{align}

\end{prop}

\begin{prop}\label{PropAppT} Assume the above sieve set-up. Let $\varepsilon>0$ be small. Let $\delta\in (0, 10^{-3}]$ and let $b$ be an odd integer that is sufficiently large in terms of $\delta$. Let $\alpha=(1/3-2\delta)^{-1} + \varepsilon$ be such that $\alpha\in [2, 4)$ and let $\rho_{\textup{lin}}\leq \frac{1}{2}-2\delta-\varepsilon$. Then we have
\begin{align}
\label{EqT1} T\leq \dfrac{10\mathfrak{S}+o(1)}{9(\log X)^{3/2}} \dfrac{b}{\varphi(b)}\prod_{\substack{p|b\\p\equiv 3\imod 4}}\bigg(1-\dfrac{1}{p-1}\bigg)^{-1}I_{\textup{lin}}(\rho_{\textup{lin}}, \alpha)\sum_{\substack{n<X}}1_{\mathcal{A}_r}(n),
\end{align}
where $\mathfrak{S}$ is given by the relation \eqref{ac} and
\begin{align}
\label{eq: ilin} I_{\textup{lin}}(\rho_{\textup{lin}}, \alpha)=\dfrac{1}{\rho_{\textup{lin}}}\int_{2}^\alpha\dfrac{\log (y-1)}{y(1-y/\alpha)^{1/2}}\textup{d}y.
\end{align}

\end{prop}

Now we can give the proof of Theorem \ref{Main TheoremApp} from Propositions \ref{PropUpper}, \ref{PropAppS}, and \ref{PropAppT}.

\begin{proof}[Proof of Theorem \ref{Main TheoremApp} assuming Propositions \ref{PropAppS} and \ref{PropAppT}]
From \eqref{sieve lower}, \eqref{Buch}, \eqref{EqS1} and \eqref{EqT1}, we have 
\begin{align}
\notag \sum_{\substack{p<X\\p\equiv r\imod b}}1_{\mathcal{A}}(p)1_{\mathbb{B}}(p-1)\geq &~S(\mathcal{F}, \mathcal{P}_3, \sqrt{X})\\
\notag =&~S(\mathcal{F}, \mathcal{P}_3, X^{1/\alpha})-T\\
\notag \geq&~\dfrac{\mathfrak{S} + o(1)}{(\log X)^{3/2}}\dfrac{b}{\varphi(b)}\prod_{\substack{p|b\\p\equiv 3\imod 4}}\bigg(1-\dfrac{1}{p-1}\bigg)^{-1}\\
\notag &\times \Big(I_{\textup{sem}}(\rho_{\textup{sem}}, \alpha)-\dfrac{10}{9}\cdot I_{\textup{lin}}(\rho_{\textup{lin}}, \alpha)+o(1)\Big)\sum_{\substack{n<X}}1_{\mathcal{A}_r}(n).
\end{align}
A simple numerical computation yields that
\begin{align}
\notag I_{\textup{sem}}(\rho_{\textup{sem}}, \alpha)-\dfrac{10}{9}\cdot I_{\textup{lin}}(\rho_{\textup{lin}}, \alpha)>1.60492-1.4566=0.1482 > 0
\end{align}
for $\rho_{\textup{sem}}=3(1-4\delta)/7-\varepsilon$, $\rho_{\textup{lin}}=1/2-2\delta-\varepsilon$, $\alpha=(1/3-2\delta)^{-1} +\varepsilon$, $\delta=1/1000$ with $\varepsilon>0$ small. Hence, we obtain
\begin{align}
\notag \sum_{\substack{p<X\\p\equiv r\imod b}}1_{\mathcal{A}}(p)1_{\mathbb{B}}(p-1)\gg \dfrac{X^{\zeta}}{(\log X)^{3/2}}.
\end{align}
This establishes the lower bound in Theorem \ref{Main TheoremApp}. Along with Proposition \ref{PropUpper}, this completes the proof of Theorem \ref{Main TheoremApp}.
\end{proof}

\subsection{Auxiliary results} We collect two key estimates essential for us while computing the lower bound.
\begin{lemma}\label{Aux1}
We have
\begin{align}
\notag \prod_{\substack{p\leq y\\p\equiv 3\imod 4}}\bigg(1-\dfrac{1}{\varphi(p)}\bigg)=2C_2C_3\Big(1+o(1)\Big)\bigg(\dfrac{\pi e^{-\gamma}}{\log y}\bigg)^{1/2} \quad \text{as} \quad y\rightarrow\infty,
\end{align}
where $\gamma$ is the Euler-Mascheroni constant,
\begin{align}
\notag C_2=\dfrac{1}{2\sqrt{2}}\prod_{p\equiv 3\imod 4}\bigg(1-\dfrac{1}{p^2}\bigg)^{1/2}\quad \text{and}\quad C_3=\prod_{p\equiv 3\imod 4}\bigg(1-\dfrac{1}{(p-1)^2}\bigg).
\end{align}
\end{lemma}
\begin{proof}
The proof is standard and can be easily derived following \cite[pp.~277--278]{opera}.
\end{proof}

\begin{lemma}\label{Aux2}
Let $\mathcal{L}$ be as in \eqref{mathcalL} and let $\alpha \in [2, 4)$. For any positive integer $n\geq 3$, let
\begin{align}\label{eq: mathfrakt}
    \mathfrak{t}(n) = \prod_{\substack{p|n\\ p>2}}\dfrac{p-1}{p-2}.
\end{align} Then, we have
\begin{align}
\notag \sum_{\substack{\ell \in \mathcal{L}\\ (\ell, 2b)=1}}\dfrac{\mathfrak{t}(\ell)}{\ell \log (X/\ell)}=\dfrac{1+o(1)}{(\log X)^{1/2}}\dfrac{C_2}{2C_1}\prod_{\substack{p|b\\p\equiv 1\imod 4}}\bigg(1+\dfrac{1}{p-2}\bigg)^{-1}\int_{2}^\alpha \dfrac{\log (y-1)}{y(1-y/\alpha)^{1/2}}\textup{d}y,
\end{align}
where 
\begin{align}
\notag C_2=\dfrac{1}{2\sqrt{2}}\prod_{p\equiv 3\imod 4}\bigg(1-\dfrac{1}{p^2}\bigg)^{1/2}\quad \text{and}\quad C_1=\prod_{p\equiv 1\imod 4}\bigg(1-\dfrac{1}{(p-1)^2}\bigg).
\end{align}
\end{lemma}
\begin{proof}
The proof follows from the proof of \cite[Lemma 5]{Mat2} in conjunction with \cite[Satz 1]{Wir} to incorporate the extra condition $(\ell, 2b)=1$.
\end{proof}

\subsection{Proof of Proposition \ref{PropAppS}}

We establish Proposition \ref{PropAppS} assuming Theorem \ref{Hyp1}, given below.

\begin{thm}[Semi-linear sieve equidistribution estimate]\label{Hyp1}
Let $\varepsilon>0$ be small. Let $\delta \in (0, 10^{-3}]$ and let $b$ be an odd integer that is sufficiently large in terms of $\delta$. Let $r\in \mathcal{A}\cap[0, b)$ with $(r(r-1), b)=1$. Let $\lambda^-_{\textup{sem}}$ be as in Lemma \ref{fundsemi} and Lemma \ref{wellfacsem} with $z_1\leq X^{1/3-2\delta-2\varepsilon^2}$ and $D=X^{\rho_{\textup{sem}}}$, where $\rho_{\textup{sem}}=3(1-4\delta)/7-\varepsilon$. Then for any $A>0$, we have
\begin{align}\notag 
\sum_{\substack{d\leq D\\(d,2b)=1}}\lambda^-_{\textup{sem}}(d)\Bigg(\sum_{\substack{n < X\\n\equiv 1 \imod d\\n\equiv 3\imod 8}}\Lambda(n)1_{\mathcal{A}_r}(n) - \dfrac{1}{4\varphi(d)}\dfrac{b}{\varphi(b)}\sum_{\substack{n<X}}1_{\mathcal{A}_r}(n)\Bigg) \ll_{A, b, \delta, \varepsilon} \dfrac{X^{\zeta}}{(\log X)^A}.
\end{align}
\end{thm}

\begin{proof}[Proof of Proposition \ref{PropAppS} assuming Theorem \ref{Hyp1}]
We have
\begin{align}
\label{semi0} S\geq \dfrac{1}{\log X}\sum_{\substack{p<X\\(p-1,\: P_3(X^{1/\alpha}))=1\\p\equiv 3\imod 8}}1_{\mathcal{A}_r}(p)\log p.
\end{align}
Next, for $d|P_3(X^{1/\alpha})=\prod_{p<X^{1/\alpha},\: p\in \mathcal{P}_3}p$, where $\mathcal{P}_3=\{p\equiv 3\imod 4: p\nmid b\}$, let
\begin{align}
\notag {E}_1(d)=\sum_{\substack{p<X\\p\equiv 1\imod {d}\\p\equiv 3\imod 8}}1_{\mathcal{A}_r}(p)\log p-\dfrac{1}{4\varphi(d)}\dfrac{b}{\varphi(b)}\sum_{n<X}1_{\mathcal{A}_r}(n).
\end{align}
Now we choose $c(n)=1_{\mathcal{A}_r\cap \mathbb{P}}(n)\log n$ for $n<X$ and $n\equiv 3\imod 8$ in Lemma \ref{fundsemi}. Then, for $1\leq u_1\leq 3$, the lower bound semi-linear sieve \eqref{lower} yields
\begin{align}
\notag \sum_{\substack{p<X\\(p-1,\: P_3(X^{1/\alpha}))=1\\p\equiv 3\imod 8}}1_{\mathcal{A}_r}(p)\log p \geq \dfrac{1}{4}\Big(f_{\textup{sem}}(u)+o(1)\Big)V_{\textup{sem}}(X^{1/\alpha})\dfrac{b}{\varphi(b)}\sum_{n<X}1_{\mathcal{A}_r}(n)\\
\label{semi1} + \sum_{\substack{d\leq X^{u_1/\alpha}\\(d, 2b)=1}}\lambda^-_{\textup{sem}}(d){E}_1(d),
\end{align}
where $\lambda^-_{\textup{sem}}$ are the lower bound semi-linear sieve weights with sifting parameter $z_1=X^{1/\alpha}$, $f_{\textup{sem}}(u_1)$ is given by \eqref{fsem}, and 
\begin{align}
 V_{\textup{sem}}(X^{1/\alpha})=\prod_{\substack{p<X^{1/\alpha}\\p\equiv 3\imod 4\\(p,b)=1}}\bigg(1-\dfrac{1}{\varphi(p)}\bigg).
\end{align}
We have $z_1=X^{1/\alpha}\leq X^{1/3-2\delta-2\varepsilon^2}$, so that we can take $u_1=\rho_{\textup{sem}} \alpha$, where $\rho_{\textup{sem}}=\frac{3}{7}(1-4\delta)-\varepsilon$ in Theorem \ref{Hyp1}. We can then use Theorem \ref{Hyp1} and the fact that the contribution of prime powers is negligible to bound the error term $E_1(d)$. In fact, using Chebyshev's estimate \cite[Theorem 2.4]{kou}, the contribution of prime powers can be bounded by 
\[\ll \sum_{d\leq X^{\frac{3}{7}(1-4\delta)-\varepsilon}}\sum_{\substack{p^m<X\\ p\equiv 1\imod d\\ p\equiv 3\imod 8\\ m\geq 2}}1_{\mathcal{A}_r}(p)\log p \ll (\log X)\sum_{d\leq X^{\frac{3}{7}(1-4\delta)-\varepsilon}}\sum_{p\leq X^{1/2}}1 \ll X^{13/14-12\delta/7},\]
which is admissible. Hence, the error term in \eqref{semi1} can be bounded as
\begin{align}\label{semierr}
\sum_{\substack{d\leq X^{\rho_{\textup{sem}}}\\(d, 2b)=1}}\lambda^-_{\textup{sem}}(d){E}_1(d) \ll_{A, b, \delta, \varepsilon} \dfrac{X^{\zeta}}{(\log X)^A}.
\end{align}
Next, we simplify the main term in \eqref{semi1} using Lemma \ref{Aux1}, so that
\begin{align}\label{Vx}
V_{\textup{sem}}(X^{1/\alpha}) = (1+o(1))\prod_{\substack{p|b\\p\equiv 3\imod 4}}\bigg(1-\dfrac{1}{p-1}\bigg)^{-1}\cdot 2C_2C_3\cdot \bigg(\dfrac{\alpha \pi e^{-\gamma}}{\log X}\bigg)^{1/2},
\end{align}
where 
\begin{align}
\label{A} C_2 =\dfrac{1}{2\sqrt{2}}\prod_{p\equiv 3\imod 4}\bigg(1-\dfrac{1}{p^2}\bigg)^{1/2}\quad 
\text{and} \quad C_3 = \prod_{p\equiv 3\imod 4}\bigg(1-\dfrac{1}{(p-1)^2}\bigg).
\end{align}
Putting the estimates from \eqref{semi1}, \eqref{semierr} and \eqref{Vx} in \eqref{semi0}, and noting that $u_1=\alpha\rho_{\textup{sem}}$, we have 
\begin{align}
\notag S \geq &~\dfrac{2C_2C_3(1+o(1))}{4(\log X)^{3/2}} \prod_{\substack{p|b\\p\equiv 3\imod 4}}\bigg(1-\dfrac{1}{p-1}\bigg)^{-1}\bigg(\dfrac{\alpha}{u_1}\bigg)^{1/2} \int_{1}^{u_1}\dfrac{\text{d}y}{\sqrt{y(y-1)}}\\ \notag &\times \dfrac{b}{\varphi(b)}\sum_{n<X}1_{\mathcal{A}_r}(n)\\
\notag = &~\dfrac{C_2C_3(1 + o(1))}{2(\log X)^{3/2}}\prod_{\substack{p|b\\p\equiv 3\imod 4}}\bigg(1-\dfrac{1}{p-1}\bigg)^{-1}\dfrac{b}{\varphi(b)}\sum_{n<X}1_{\mathcal{A}_r}(n)I_{\textup{sem}}(\rho_{\textup{sem}}, \alpha),
\end{align}
where $I_{\textup{sem}}(\rho_{\textup{sem}}, \alpha)$ is given by \eqref{eq: isem} 
Therefore,
\begin{align}
\label{semfin} S\geq \dfrac{\mathfrak{S} + o(1)}{(\log X)^{3/2}}\dfrac{b}{\varphi(b)}\prod_{\substack{p|b\\p\equiv 3\imod 4}}\bigg(1-\dfrac{1}{p-1}\bigg)^{-1}I_{\textup{sem}}(\rho_{\textup{sem}}, \alpha)\sum_{\substack{n<X}}1_{\mathcal{A}_r}(n),
\end{align} 
where
\begin{align}
\label{mathfrakS} \mathfrak{S}=\dfrac{C_2C_3}{2}=\dfrac{1}{4\sqrt{2}}\prod_{p\equiv 3\imod 4}\bigg(1-\dfrac{1}{p^2}\bigg)^{1/2}\prod_{p\equiv 3\imod 4}\bigg(1-\dfrac{1}{(p-1)^2}\bigg).
\end{align}
This completes the proof of Proposition \ref{PropAppS}.
\end{proof}

Thus, we are left to establish Theorem \ref{Hyp1}, which we do in Part \ref{part: circle}.

\subsection{Proof of Proposition \ref{PropAppT}}

Finally, we give the proof of Proposition \ref{PropAppT} assuming Theorem \ref{Hyp2}, given below.

\begin{thm}[Linear sieve equidistribution estimate]\label{Hyp2}
Let $\varepsilon>0$ be small. Let $\delta\in (0, 10^{-3}]$ and let $b$ be an odd integer that is sufficiently large in terms of $\delta$. Let $r\in \mathcal{A}\cap[0,b)$ with $(r, b)=(r-1, b)=1$. Let $L$ be a real number such that $L\in[X^{1/3-2\delta-\varepsilon}, X^{2/3+2\delta + \varepsilon}]$. Suppose $\mathfrak{h}$ is a bounded arithmetic real-valued function, and $\lambda^+_{\textup{lin}}$ is as in Lemma \ref{fundlin} and Lemma \ref{wellfacclin} with $z_2=X^{1/5}$, and $D=X^{\rho_{\textup{lin}}}$ for $\rho_{\textup{lin}}=1/2-2\delta-\varepsilon$. Then for any $A>0$, we have
\begin{align}
\notag \sum_{\substack{d \leq X^{1/2-2\delta}\\(d,2b)=1}}\lambda^+_{\textup{lin}}(d)\bigg(\sum_{\substack{\ell\sim L\\(\ell,2b)=1}}\mathfrak{h}(\ell)&\sum_{\substack{n< X/2\ell\\2\ell n+1\equiv 0\imod {d}\\ \ell n\equiv 1\imod 4}}1_{\mathcal{A}_r}(2\ell n+1)\Lambda(n)\\
\label{hyp2eqn} &-\dfrac{1}{4\varphi(d)}\dfrac{b}{\varphi(b)}\sum_{\substack{\ell \sim L\\(\ell, 2bd)=1}}\dfrac{\mathfrak{h}(\ell)}{\ell}\sum_{n<X}1_{\mathcal{A}_r}(n)\bigg)\ll_{A, b, \delta, \varepsilon} \dfrac{X^{\zeta}}{(\log X)^A}.
\end{align}
\end{thm}

\begin{proof}[Proof of Proposition \ref{PropAppT} assuming Theorem \ref{Hyp2}]
 By the inequality \eqref{Tsimple}, for some parameter $\nu$ (to be chosen later), we find that
\begin{align}
\label{T00} T\leq \sum_{\substack{\ell \in \mathcal{L}\\(\ell, b)=1}}\bigg( S(\mathcal{M}(\ell), \mathcal{P}(\ell), X^{1/\nu}) + O(X^{1/\nu})\bigg),
\end{align}
where $\mathcal{L}$ and $\mathcal{M}(\ell)$ are given by \eqref{mathcalL} and \eqref{mathcalMl}, respectively. Furthermore, $\mathcal{P}(\ell) = \{p: p\nmid 2b\ell\}$.  

Next, we set $P_\ell (X^{1/\nu}): = \prod_{p < X^{1/\nu},\: p\in \mathcal{P}(\ell)}p$ and note that
\begin{align}
\notag \sum_{\substack{\ell \in \mathcal{L}\\(\ell, b)=1}}\bigg( S(\mathcal{M}(\ell), \mathcal{P}(\ell), X^{1/\nu}) + O(X^{1/\nu}\bigg)\leq \sum_{\substack{\ell \in \mathcal{L}\\(\ell, b)=1}}\sideset{}{^\flat}\sum_{\substack{p_2< X/2\ell}}1_{\mathcal{A}_r}(2\ell p_2 +1)+ O(X^{1/\nu}\#\mathcal{L}),
\end{align}
where $\sum^\flat$ denotes a sum over values of $p_2$ satisfying  
\[\ell p_2\equiv 1\imod 4 \quad \text{and} \quad \big(2\ell p_2+1, P_{\ell}(X^{1/\nu})\big)=1.\]
As in the proof of Proposition \ref{PropUpper}, we first split the range of $p_2$ to obtain
\begin{align}
\notag \sideset{}{^\flat}\sum_{\substack{p_2< X/2\ell}}1_{\mathcal{A}_r}(2\ell p_2 +1) \leq &~\dfrac{10}{9\log (X/\ell)}\sum_{\substack{n< X/2\ell\\\ell n\equiv 1\imod 4\\ (2\ell n + 1, P_{\ell}(X^{1/\nu}))=1}}\Lambda(n)1_{\mathcal{A}_r}(2\ell n +1)\\
\notag  &+ \sideset{}{^\flat}\sum_{\substack{p_2\leq (X/\ell)^{9/10}}}1_{\mathcal{A}_r}(2\ell p_2 +1).
\end{align}
 Next, we use Chebyshev's bound \cite[Theorem 2.4]{kou} for the sum over primes $p_2$. Note that, since $\alpha=(1/3-2\delta)^{-1} + \varepsilon$, by \eqref{mathcalL}, we have 
\begin{align}\label{Range: L}
\mathcal{L} \subset [X^{1/\alpha}, X^{1-1/\alpha}] \subset[X^{1/3-2\delta-\varepsilon}, X^{2/3+2\delta +\varepsilon}].
\end{align}
This allows us to bound the second sum as
\begin{align}
\notag \sum_{\substack{\ell \in \mathcal{L}\\(\ell, b)=1}}\sideset{}{^\flat}\sum_{\substack{p_2\leq (X/\ell)^{9/10}}}1_{\mathcal{A}_r}(2\ell p_2 +1)\ll \sum_{\ell\in [X^{1/3-2\delta-\varepsilon}, X^{2/3+2\delta+\varepsilon}]}\dfrac{(X/\ell)^{9/10}}{\log (X/\ell)} \ll X^{29/30+\delta/5+\varepsilon}.
\end{align}
The above estimates yield
\begin{align}
\label{T0}T\leq \dfrac{10}{9}\sum_{\substack{\ell \in \mathcal{L}\\(\ell, b)=1}}\dfrac{1}{\log (X/\ell)}\sideset{}{^{\flat\flat}}\sum_{\substack{n< X/2\ell}}\Lambda(n)1_{\mathcal{A}_r}(2\ell n +1) + O\bigg(X^{29/30+\delta/5+\varepsilon} + \#\mathcal{L}X^{1/\nu}\bigg),
\end{align}
where $\sum^{\flat\flat}$ denotes a sum over values of $n$ satisfying  
\[\ell n\equiv 1\imod 4 \quad \text{and} \quad \big(2\ell n + 1, P_{\ell}(X^{1/\nu})\big)=1.\]
Next, for $d|\prod_{p<z,\: p\in \mathcal{P}(\ell)}p$, where $\mathcal{P}(\ell)=\{p:p\nmid 2b\ell\}$, we let
\begin{align}
\notag {E}_2(d)=\sum_{\substack{n< X/2\ell\\  2 \ell n+1\equiv 0\imod {d}\\ \ell n\equiv 1\imod 4}}\Lambda(n)1_{\mathcal{A}_r}(2\ell n+1)-\dfrac{1}{4\varphi(d)}\dfrac{b}{\varphi(b)}\dfrac{1}{\ell}\sum_{n<X}1_{\mathcal{A}_r}(n).
\end{align}
We now apply Lemma \ref{fundlin} with the sequence $c(n) = \Lambda(n)1_{\mathcal{A}_r}(2\ell n +1)$ for $n<X/2\ell$ and $\ell n\equiv 1\imod 4$. Then given a parameter $u_2\in[1, 3]$ to be chosen later, the upper bound linear sieve \eqref{upperl} yields
\begin{align}\label{Sl}
\begin{aligned}
 \sideset{}{^{\flat\flat}}\sum_{\substack{n< X/2\ell}}\Lambda(n)1_{\mathcal{A}_r}(2\ell n +1)\leq  \dfrac{b}{4\varphi(b)}\big(F_{\textup{lin}}(u_2)+o(1)\big) &V_{\textup{lin}}(X^{1/\nu}) \dfrac{1}{\ell}\sum_{n<X}1_{\mathcal{A}_r}(n)\\
 &+ \sum_{\substack{d\leq X^{u_2/\nu}\\(d,2\ell b)=1}}\lambda^+_{\textup{lin}}(d){E}_2(d),
\end{aligned}
\end{align}
where $\lambda^+_{\textup{lin}}$ are the upper bound linear sieve weights with sifting parameter $z_2=X^{1/\nu}, F_{\textup{lin}}(u_2)=2e^\gamma/u_2$ and 
\begin{align}
\notag V_{\textup{lin}}(X^{1/\nu}) = \prod_{\substack{p< X^{1/\nu}\\ (p, 2\ell b)=1}}\bigg(1-\dfrac{1}{\varphi(p)}\bigg).
\end{align}
Since $(\ell, b)=1$ in \eqref{T0}, we may use Mertens' theorem \cite[Theorem 3.4(c)]{kou} to obtain
\begin{align}
\label{Eq: Vlin} V_{\textup{lin}}(X^{1/\nu}) = \prod_{\substack{p< X^{1/\nu}\\ (p, 2\ell b)=1}}\bigg(1-\dfrac{1}{\varphi(p)}\bigg)=\big(1+o(1)\big)\dfrac{2\nu C_1C_3e^{-\gamma}\mathfrak{t}(\ell)\mathfrak{t}(b)}{\log X},
\end{align}
where $\mathfrak{t}(n)$ is given by \eqref{eq: mathfrakt},
\begin{align}
\label{C_1C_3} C_1 &= \prod_{p\equiv 1 \imod 4}\bigg(1-\dfrac{1}{(p-1)^2}\bigg), \quad \text{and} \quad C_3 = \prod_{p\equiv 3\imod 4}\bigg(1-\dfrac{1}{(p-1)^2}\bigg).
\end{align}

Now we take $u_2=\rho_{\textup{lin}} \nu$ in the linear sieve, where $\rho_{\textup{lin}}$ corresponds to the level of the upper bound sieve in Theorem \ref{Hyp2}. Next, using \eqref{Range: L}, we write 
\begin{align}
\notag \sum_{d\leq X^{\rho_{\textup{lin}}}}\lambda^+_{\textup{lin}}(d)\sum_{\substack{\ell \in \mathcal{L}\\(\ell, 2bd)=1}}\dfrac{1}{\log (X/\ell)}{E}_2(d)&=\sum_{d\leq X^{\rho_{\textup{lin}}}}\lambda^+_{\textup{lin}}(d)\sum_{\substack{X^{1/3-2\delta-\varepsilon}<\ell\leq X^{2/3+2\delta+\varepsilon}\\(\ell, 2bd)=1}}\dfrac{1_{\mathcal{L}}(\ell)}{\log (X/\ell)}{E}_2(d).
\end{align}
We do a dyadic decomposition on the range of $\ell$, say $\ell \sim L$ with $L\in [X^{1/3-2\delta-\varepsilon}, X^{2/3+2\delta+\varepsilon}]$. Since the number of such dyadic intervals are at most $\log X$, we use Theorem \ref{Hyp2} with $\mathfrak{h}(\ell) = 1_{\mathcal{L}}(\ell)/\log (X/\ell)$ for $\ell\sim L$ to bound the above expression as
\begin{align}
\notag \sum_{d\leq X^{\rho_{\textup{lin}}}}\lambda^+_{\textup{lin}}(d)\sum_{\substack{X^{1/3-2\delta-\varepsilon}<\ell\leq X^{2/3+2\delta+\varepsilon}\\(\ell, 2bd)=1}}\mathfrak{h}(\ell){E}_2(d)&\ll (\log X)\bigg|\sum_{d\leq X^{\rho_{\textup{lin}}}}\lambda^+_{\textup{lin}}(d)\sum_{\substack{\ell \sim L\\(\ell, 2bd)=1}}\mathfrak{h}(\ell){E}_2(d)\bigg|\\
\notag &\ll_{A, b, \delta, \varepsilon} \dfrac{X^{\zeta}}{(\log X)^A},
\end{align}
which is admissible.

From \eqref{T0}, \eqref{Sl} and using the above bound for the error term, we have that
\begin{align}
\notag T \leq &~\dfrac{10}{9}\cdot \dfrac{b}{4\varphi(b)}\big(F_{\textup{lin}}(u_2)+o(1)\big)\cdot \sum_{n<X}1_{\mathcal{A}_r}(n)\sum_{\substack{\ell \in \mathcal{L}\\(\ell, b)=1}}\dfrac{V_{\textup{lin}}(X^{1/\nu})}{\ell \log(X/\ell)}\\
\notag & + O_{A, b, \delta, \varepsilon}\bigg(\dfrac{X^{\zeta}}{(\log X)^A} + X^{29/30+\delta/5+\varepsilon} + \#\mathcal{L}X^{1/\nu}\bigg)\\
 \notag\leq &~\dfrac{10}{9}\cdot \dfrac{\nu\cdot C_1C_3e^{-\gamma}}{2\log X}\cdot \dfrac{b\cdot \mathfrak{t}(b)}{\varphi(b)}\big(F_{\textup{lin}}(u_2)+o(1)\big)\sum_{\substack{n<X}}1_{\mathcal{A}_r}(n)\sum_{\substack{\ell \in \mathcal{L}\\(\ell, 2b)=1}}\dfrac
{\mathfrak{t}(\ell)}{\ell \log (X/\ell)}\\
\label{Tefin} &+  O_{A, b, \delta, \varepsilon}\bigg(\dfrac{X^{\zeta}}{(\log X)^A} + X^{29/30+\delta/5+\varepsilon} + \#\mathcal{L}X^{1/\nu}\bigg),
\end{align}
where we have used the asymptotic formula for $V_{\textup{lin}}(X^{1/\nu})$ from the relation \eqref{Eq: Vlin} in the last line. Next, by Lemma \ref{Aux2}, we have
\begin{align}
\label{tlem}\sum_{\substack{\ell \in \mathcal{L}\\(\ell, 2b)=1}}\dfrac
{\mathfrak{t}(\ell)}{\ell \log (X/\ell)}=\dfrac{(1+o(1))C_2}{2C_1(\log X)^{1/2}}\prod_{\substack{p|b\\p\equiv 1\imod 4}}\bigg(1+\dfrac{1}{p-2}\bigg)^{-1}\int_{2}^\alpha\dfrac{\log (y-1)}{y(1-y/\alpha)^{1/2}}\text{d}y, 
\end{align}
where $C_2$ is as in the relation \eqref{A} and $C_1$ is given by \eqref{C_1C_3}.

We choose $u_2=5/2$, so that $F_{\textup{lin}}(u_2)=4e^\gamma/5$. As by our choice, $\rho_{\textup{lin}}=1/2-2\delta-\varepsilon$, we can choose $\nu=5$.  Note that since $\nu=5, \#\mathcal{L} \leq X^{2/3+2\delta+\varepsilon}$, $\varepsilon>0$ is small enough, $\delta\in (0, 10^{-3}]$, and $\zeta$ tends to $1$ as $b \rightarrow \infty$, we have that
\begin{align}\label{tefin2}
\#\mathcal{L}X^{1/\nu}, X^{29/30+\delta/5+\varepsilon}  \ll_{A, \delta, \varepsilon} \dfrac{X^\zeta}{(\log X)^A}.
\end{align}
Therefore, we substitute \eqref{tlem} and \eqref{tefin2} in \eqref{Tefin} to obtain
\begin{equation}
\begin{aligned}
\label{tfin} T\leq &~\dfrac{10\mathfrak{S}+o(1)}{9(\log X)^{3/2}} \dfrac{b}{\varphi(b)}\prod_{\substack{p|b\\p>2}}\bigg(1-\dfrac{1}{p-1}\bigg)^{-1}\prod_{\substack{p|b\\p\equiv 1\imod 4}}\bigg(1+\dfrac{1}{p-2}\bigg)^{-1}I_{\textup{lin}}(\rho_{\textup{lin}}, \alpha)\\
&\times \sum_{\substack{n<X}}1_{\mathcal{A}_r}(n),
\end{aligned}
\end{equation}
where $I_{\textup{lin}}(\rho_{\textup{lin}}, \alpha)$ is given by \eqref{eq: ilin} and $\mathfrak{S}=C_2C_3/2$ is given by the relation \eqref{mathfrakS}. Hence, the estimate \eqref{tfin} along with the fact that 
\[\prod_{\substack{p|b\\p>2}}\bigg(1-\dfrac{1}{p-1}\bigg)^{-1}\prod_{\substack{p|b\\p\equiv 1\imod 4}}\bigg(1+\dfrac{1}{p-2}\bigg)^{-1}=\prod_{\substack{p|b\\p\equiv 3\imod 4}}\bigg(1-\dfrac{1}{p-1}\bigg)^{-1}\]
yields the required bound for the sum $T$. 
\end{proof}

We have therefore established Proposition \ref{PropAppT} assuming Theorem \ref{Hyp2}. So, we are left to establish Theorems \ref{Hyp1} and \ref{Hyp2}, which we do in Part \ref{part: circle}.

\part{Exponential sums}\label{part: exp sum}
In this part, we estimate the exponential sums over primes in arithmetic progressions using Vinogradov's method (see \cite[Chapter 23]{kou} for an introduction to the method), which we will employ in Part \ref{part: circle} to deduce our main results. Note that some of the estimates in this part are well-known. See, for example, \cite{Mat}, \cite{Mik}.

Recall that we set $X=b^k$ with $k\in \mathbb{Z}$ and $k\rightarrow \infty$ throughout this paper. We remark that the results in this part of the paper hold for any large real number $X$.

\section{Preliminary estimates and Type I estimate}\label{sec: type I}
We begin with the following estimate.

\begin{lemma}\label{lemma1}
Let $\theta = a/q + \beta$ with $(a, q)=1$ and $0< |\beta| < 1/q^2$. Then for any $M, N \geq 2$, we have
\begin{equation}
\notag \sum_{m=1}^M\min \bigg(N, \dfrac{1}{\lVert m\theta \rVert}\bigg) \ll \bigg(M + MNq\lvert \beta\rvert + \dfrac{1}{q\lvert \beta \rvert}\bigg)(\log 2qM).
\end{equation} 
\end{lemma}
\begin{proof}
The proof of the lemma is a standard one. However, we need a variant of it to take advantage of $\beta$ in the sum. For a detailed proof, see \cite[Lemma 4.1]{May2}.
\end{proof}

Let us now deduce the following corollary from the above lemma.
\begin{cor}\label{Corollary1}
Let $\theta = a/q + \beta$ with $(a, q)=1$ and $|\beta| < 1/q^2$. Then for any $M\geq 1$, we have
\begin{align}
\label{cor1eq}\sum_{m\leq M}\min\bigg(\dfrac{X}{m}+1,\dfrac{1}{\|m\theta\|}\bigg) \ll X\bigg(\dfrac{M}{X} + \dfrac{qH}{X}+ \dfrac{1}{qH}\bigg)(\log 2qM)^2,
\end{align}
where $H=1+|\beta|X$.
\end{cor}

\begin{proof}
If $\beta\neq 0$, we perform a dyadic decomposition and then apply Lemma \ref{lemma1} to obtain
\begin{align}
\label{cor5.2i} \sum_{m\leq M}\min \bigg(\dfrac{X}{m}+1,\dfrac{1}{\|m\theta\|}\bigg)
 &\ll  X\bigg(\dfrac{M}{X} + q\lvert\beta\rvert + \dfrac{1}{Xq\lvert\beta\rvert}\bigg)(\log 2qM)^2.
\end{align}
Next, for all $\beta$, we apply \cite[Lemma 13.7]{iwakow} to obtain
\begin{align}
\label{cor5.2ii} \sum_{m\leq M}\bigg(\dfrac{X}{m} + 1, \dfrac{1}{\|m\theta\|}\bigg) \ll X\bigg(\dfrac{M}{X} + \dfrac{1}{q} + \dfrac{q}{X}\bigg)(\log 2qM).
\end{align}
Combining the estimates from inequalities \eqref{cor5.2i} and\eqref{cor5.2ii}, we have
\begin{align}\label{cor5.2iii}
\sum_{m\leq M}\bigg(\dfrac{X}{m} + 1, \dfrac{1}{\|m\theta\|}\bigg) \ll X\bigg\{\dfrac{M}{X} + \min\bigg(q|\beta| + \dfrac{1}{Xq|\beta|}, \dfrac{1}{q} + \dfrac{q}{X}\bigg)\bigg\}(\log 2qM)^2.
\end{align}
Next, we note that
\begin{align}\label{cor5.2iv}
\min\bigg(\dfrac{1}{q}, \dfrac{1}{Xq|\beta|}\bigg)\leq \dfrac{2}{q(1+|\beta|X)}.
\end{align}
Therefore, using \eqref{cor5.2iv} and recalling that $|\beta|<1/q^2$, we obtain
\begin{align}\notag 
\min\bigg(q|\beta| + \dfrac{1}{Xq|\beta|}, \dfrac{1}{q} + \dfrac{q}{X}\bigg)\ll q|\beta| + \dfrac{q}{X} + \dfrac{1}{q(1+|\beta|X)} = \dfrac{qH}{X} + \dfrac{1}{qH},
\end{align}
where $H=1+|\beta|X$. Substituting the above inequality in \eqref{cor5.2iii} completes the proof of the corollary.
\end{proof}

We now state the following bilinear sum estimates for the exponential sum.
\begin{lemma}[Bilinear estimate]\label{bilinear lemma}
Let $M, N \geq 1$ be such that $MN\leq X$. Let $\theta = a/q + \beta$ for some $(a, q)=1$ and $|\beta| < 1/q^2$. Suppose ${\alpha}_1$ and ${\alpha}_2$ are two arithmetic functions supported in $[1, M]$ and $[1, N]$, respectively. Then we have
\begin{align}\notag
\sum_{n\leq X}({\alpha}
_1*{\alpha}_2)(n)e(n \theta) \ll X^{1/2}\lVert {\alpha}_1 \rVert_2 \lVert {\alpha}_2 \rVert_2\bigg(\dfrac{M}{X} + \dfrac{N}{X} + \dfrac{qH}{X} + \dfrac{1}{qH} \bigg)^{1/2}(\log 2qX),
\end{align}
where $H=1+|\beta|X$.
\end{lemma}
\begin{proof}
The proof follows by combining the argument of \cite[Theorem 23.6]{kou} with Lemma \ref{lemma1} and Corollary \ref{Corollary1}.
\end{proof}

Next, we will need an auxiliary lemma due to Matom\"aki \cite[Lemma 8]{Mat}, who improved on the earlier work of Mikawa \cite{Mik}.

\begin{lemma}[Matom\"aki]\label{Mikawa}
Let $M, N\geq 1$ be such that $M, N\leq X$. Let $\theta= a/q + \beta$ with $(a, q)=1$, $|\beta\rvert <1/q^2$ and $q<X$. Then for any $\psi>0$, one has
\begin{align}
\notag M\sum_{m\sim M}\sum_{n\sim N}\uptau_3(n)\min\bigg(\dfrac{X}{m^2n} + 1, \dfrac{1}{\lVert m^2n \theta \rVert}\bigg)\ll &~M^2N(\log X)^3\\
\notag & + X\bigg(\dfrac{1}{M} + \dfrac{qH}{X} + \dfrac{1}{qH}\bigg)^{1/2-\psi}(\log X)^8,
\end{align}
where $H=1+|\beta|X$.
\end{lemma}

\begin{proof}
The proof follows from the argument of \cite[Lemma 8]{Mat} in conjunction with Lemma \ref{lemma1} and Corollary \ref{Corollary1}.
\end{proof}

\subsection{Type I estimate}
We will estimate the so-called Type I sum in the following lemma.
\begin{lemma}[Type I estimate]\label{lemmatype1}
Let $v>0$. Let $D, M\geq 1$ be such that $DM<X$. Let $\theta=a/q+\beta$ with $(a, q)=1$, $\lvert \beta \rvert <1/q^2$ and $q<X$. Suppose ${\alpha}$ is an arithmetic function supported in $[1, M]$ and satisfies $\lvert {\alpha}\rvert \leq \uptau_{h_1}\cdot \log^{h_2}$ for some fixed integers $h_1\geq 1, h_2\geq 0$. Furthermore, let $h_3\geq 1$ be a fixed integer. Then, we have
\begin{align}
\notag \sum_{\substack{d\leq D}}\uptau_{h_3}(d)\cdot\max_{(c,d)=1}&\bigg\lvert \sum_{\substack{mn <X\\1\leq m\leq M\\mn\equiv c\imod d}}{\alpha}(m)(\log n)^ve(mn\theta)\bigg\rvert\\
\notag &\ll X\bigg(\dfrac{DM}{X} + \dfrac{qH}{X} + \dfrac{1}{qH}\bigg)^{1/2}(\log X)^{(h_1+h_3)^2/2 + h_2+v + 1},
\end{align}
where $H=1+|\beta|X$.
\end{lemma}

\begin{proof}
Let $\mathcal{S}_{\textup{Type I}}$ be the sum that we wish to estimate. Applying partial summation and then using the fact that $\sum_{n\leq y}e(nt) \ll \min(y, \|t\|^{-1})$ for any real numbers $y>1$ and $t$, we have
\begin{align}
\notag \sum_{\substack{n< X/m\\n\equiv c\overline{m}\imod {d}}}(\log n)^ve(mn\theta)  &\ll (\log X/m)^v \min\bigg(\dfrac{X}{dm} + 1, \dfrac{1}{\|dm\theta\|}\bigg).
\end{align}
This implies that
\begin{align}
\notag |\mathcal{S}_{\textup{Type I}}|&\ll(\log X)^{v} \sum_{d\leq D}\uptau_{h_3}(d)\sum_{\substack{m\leq M}}|{\alpha}(m)|\min\bigg(\dfrac{X}{dm} + 1, \dfrac{1}{\lVert dm\theta\rVert}\bigg),
\end{align}
Next, we write $d^\prime =dm$, so that $d^\prime \leq DM$. Then, by the Cauchy-Schwarz inequality and Corollary \ref{Corollary1} along with the fact that $|{\alpha}| \leq \uptau_{h_1}\cdot \log^{h_2}$, we have
\begin{align}
\notag |\mathcal{S}_{\textup{Type I}}|   & \ll(\log X)^{v+h_2}\bigg(X\sum_{d^\prime \leq DM}\dfrac{\uptau_{h_1+h_3}(d^\prime)^2}{d^\prime}\bigg)^{1/2}\cdot \bigg(\sum_{d^\prime \leq DM}\min\bigg(\dfrac{X}{d^\prime} + 1, \dfrac{1}{\lVert d^\prime \theta\rVert}\bigg)\bigg)^{1/2}\\
\notag & \ll (\log X)^{v+h_2} X(\log X)^{(h_1+h_3)^2/2}\bigg(\dfrac{DM}{X}+ \dfrac{qH}{X} + \dfrac{1}{qH}\bigg)^{1/2}(\log X),
\end{align}
where we have used the fact that $\sum_{n\leq y}\uptau_h(n)^2/n \ll (\log y)^{h^2}$ for any real number $y\geq 2$ and for any integer $h\geq 1$. The above estimate on simplification yields the desired result.
\end{proof}

\section{Type II estimates}\label{sec: type II}

We use Vinogradov's method to estimate the Type II sums in the following lemma.

\begin{lemma}[Point wise Type II estimate]\label{lemmatypeIIpointwise}
Let $M, N\geq 1$ be such that $MN\leq X$. Let $\theta=a/q+\beta$ with $(a, q)=1$, $\lvert \beta \rvert <1/q^2$ and $q<X$. Suppose ${\alpha}_1$ and ${\alpha}_2$ are two arithmetic functions supported in $[M, 2M]$ and $[N, 2N]$, respectively, and satisfy $|{\alpha}_1|, |{\alpha}_2| \leq \uptau_h\cdot \log $  for some fixed integer $h\geq 1$. Let $c$ and $d$ be non-zero positive integers such that $(c, d)=1$. Then, we have
\begin{align}\label{eq: pointwise}
\bigg\lvert\sum_{\substack{mn <X\\m\sim M, n\sim N\\mn\equiv c\imod d}}{\alpha}_1(m){\alpha}_2(n)e(mn\theta)\bigg\rvert\ll X\bigg(\dfrac{M}{X}+\dfrac{N}{X}+\dfrac{qH}{X} + \dfrac{1}{qH}\bigg)^{1/2}(\log X)^{h^2 + 2},
\end{align}
where $H=1+|\beta|X$.

\end{lemma}

\begin{proof}
Let $\chi$ be Dirichlet character modulo $d$. Then, by the orthogonality of Dirichlet characters, we bound the sum in the left-hand side of \eqref{eq: pointwise} as
\begin{align}
\notag &\leq \dfrac{1}{\varphi(d)}\sum_{\chi\imod d}\bigg\lvert\sum_{\substack{mn <X\\m\sim M, n\sim N}}{\alpha}_1(m)\chi(m){\alpha}_2(n)\chi(n)e(mn\theta)\bigg \rvert.
\end{align}
Now we use Lemma \ref{bilinear lemma} with ${\alpha}_1\cdot \chi$ and ${\alpha}_2\cdot \chi$ to estimate the sum over $mn <X$  and the trivial bound to sum over $\varphi(d)$ characters modulo $\chi$ to show that the above sum is
\begin{align}
\label{lemmapeq1} &\ll X^{1/2}\lVert{\alpha}_1\rVert_2\lVert {\alpha}_2\rVert_2\bigg(\dfrac{M}{X} + \dfrac{N}{X} + \dfrac{qH}{X} + \dfrac{1}{qH}\bigg)^{1/2}(\log qX).
\end{align}
Next, we recall that $|{\alpha}_1|, |{\alpha}_2|\leq \uptau_h \cdot \log $ to obtain
\begin{align}
\label{lemmapeq2}\lVert{\alpha}_1\rVert_2\lVert{\alpha}_2\rVert_2 \leq \bigg(\sum_{m\sim M}\uptau_h(m)^2\bigg)^{1/2}\bigg(\sum_{n\sim N}\uptau_h(n)^2\bigg)^{1/2}(\log X)^2 \ll (MN)^{1/2}(\log X)^{h^2+1},
\end{align}
using the fact that $\sum_{n< y}\uptau_h(n)^2 \ll y(\log y)^{h^2-1}$ for any real number $y\geq 2$.
Substituting the estimate from \eqref{lemmapeq2} in \eqref{lemmapeq1} and using the fact that $MN\leq X$ and $q< X$ completes the proof of the lemma.
\end{proof}

In the next lemma, we improve the bounds of the previous lemma by taking advantage of averaging.

\begin{lemma}\label{lemmatypeII}
Let $c$ be a fixed non-zero integer. Let $D_1, D_2, M, N \geq 1$ be such that
\[MN< X, \quad D_1M< X \quad  \text{and} \quad D_1D_2^{2}N< X.\]
 Let $\theta=a/q+\beta$ with $(a, q)=1$, $\lvert \beta \rvert <1/q^2$ and $q<X$. Suppose ${\alpha}_1$ and ${\alpha}_2$ are two arithmetic functions with support $[M, 2M]$ and $[N, 2N]$, respectively, and satisfy $\lvert {\alpha}_1|, |{\alpha}_2| \leq \uptau_h\cdot \log $ for some fixed integer $h\geq 1$. Then, for any integer $h_1\geq 1$, we have
\begin{align}
\notag &\mathop{\sum\sum}_{\substack{d_1\sim D_1\\d_2\sim D_2\\(cd_1,d_2)=1}}\uptau_{h_1}(d_1)\max_{(c', d_1)=1}\Bigg\lvert\sum_{\substack{mn <X\\m\sim M, n\sim N\\mn\equiv c'\imod {d_1}\\mn\equiv c\imod {d_2}}}{\alpha}_1(m){\alpha}_2(n)e(mn\theta)\Bigg\rvert\\
\notag &\ll  X\bigg(\dfrac{D_1M}{X} + \dfrac{(D_1D_2)^2}{X}+\dfrac{D_1D_2^2N}{X} + \dfrac{1}{D_1^{1/4}} + \dfrac{(qH)^{1/4}}{X^{1/4}} + \dfrac{1}{(qH)^{1/4}}\bigg)^{1/2}(\log X)^{h^2+h_1^2/2 + 5},
\end{align}
where $H=1+|\beta|X$.
\end{lemma}

\begin{proof}
The proof of the lemma is closely related to the proofs of \cite[Proposition 9]{Mat} and \cite[Theorem (p.~352)]{Mik}, but for the convenience of the reader we include the proof here.
We will estimate the sum:
\begin{align}\notag 
\mathcal{S}_{\text{Type II}}:=\mathop{\sum\sum}_{\substack{d_1\sim D_1\\d_2\sim D_2\\(cd_1,d_2)=1}}\uptau_{h_1}(d_1)\max_{(c', d_1)=1}&\Bigg\lvert\sum_{\substack{mn <X\\m\sim M, n\sim N\\mn\equiv c'\imod {d_1}\\mn\equiv c\imod {d_2}}}{\alpha}_1(m){\alpha}_2(n)e(mn\theta)\Bigg\rvert.
\end{align}
Let us assume that the maximum over $c'$ is attained at $c_{d_1}$. Let ${\lambda}(d_1, d_2)\in \mathbb{C}$ be of absolute value $1$  whenever $c'=c_{d_1}$ and $(d_1, c_{d_1})=(d_1, d_2)=(d_2, c)=1$  for $d_1\sim D_1$ and $d_2\sim D_2$. Then, we have
\begin{align}
\notag \mathcal{S}_{\text{Type II}}=\sum_{d_1\sim D_1}\uptau_{h_1}(d_1)\sum_{m\sim M}{\alpha}_1(m)\sum_{d_2\sim D_2}{\lambda}(d_1, d_2)\sum_{\substack{mn <X\\n\sim N\\mn\equiv c_{d_1}\imod {d_1}\\mn\equiv c\imod {d_2}}}{\alpha}_2(n)e(mn\theta).
\end{align}
We apply the Cauchy-Schwarz inequality to obtain
\begin{align}
\notag |\mathcal{S}_{\text{Type II}}|^2\leq &~D_1 (\log X)^{h_1^2-1}\|{\alpha}_1\|_2^2\sum_{\substack{d_1\sim D_1}}\sum_{\substack{m\sim M}}\bigg|\sum_{\substack{d_2\sim D_2}}{\lambda}(d_1,d_2)\sum_{\substack{n<X/m\\ n\sim N\\ mn\equiv c_{d_1} \imod {d_1}\\mn\equiv c\imod {d_2}}}{\alpha}_2(n)e(mn\theta)\bigg|^2\\
\notag \ll &~ D_1(\log X)^{h_1^2-1}\|{\alpha}_1\|_2^2\sum_{d_1\sim D_1}\sum_{\substack{d_2, d_2'\sim D_2\\(d_2d_2', d_1)=1}}\sum_{\substack{n_1, n_2\sim N\\(n_1, d_1d_2)=(n_2, d_1d_2')=1 }}\big|{\alpha}_2(n_1)\overline{{\alpha}_2(n_2)}\big|\\
\notag &\times \Bigg|\sum_{\substack{m\sim M\\m< \min (X/n_1, X/n_2)\\mn_1\equiv mn_2\equiv c_{d_1}\imod {d_1}\\mn_1\equiv c\imod {d_2}\\mn_2\equiv c\imod {d_2'}}}e\big(m(n_1-n_2)\theta\big)\Bigg|\\
\notag \ll &~ D_1(\log X)^{h_1^2-1}\|{\alpha}_1\|_2^2\sum_{d_1\sim D_1}\sum_{\substack{d_2, d_2'\sim D_2\\(d_2d_2', d_1)=1}}\sum_{j\in \{1, 2\}}\sum_{\substack{n_1, n_2\sim N\\(n_1, d_1d_2)=(n_2, d_1d_2')=1 }}\big|{\alpha}_2(n_j)|^2\\
\notag &\times \Bigg|\sum_{\substack{m\sim M\\m< \min (X/n_1, X/n_2)\\mn_1\equiv mn_2\equiv c_{d_1}\imod {d_1}\\mn_1\equiv c\imod {d_2}\\mn_2\equiv c\imod {d_2'}}}e\big(m(n_1-n_2)\theta\big)\Bigg|
\end{align}
using the fact $|{\alpha}_2(n_1)\overline{{\alpha}_2(n_2)}|\leq |{\alpha}_2(n_1)|^2 + |{\alpha}_2(n_2)|^2$. 

The above congruences $mn_1 \equiv mn_2 \equiv c_{d_1} \imod {d_1}$, $mn_1\equiv c \imod {d_2}$, and $mn_1\equiv c\imod {d_2'}$ have a solution in $m$ if and only if $(n_1, d_1d_2)=(n_2, d_1d_2')=1$ and $n_1 \equiv n_2\imod {d_1 (d_2, d_2')}$. Then, we have a unique solution $m\equiv h'\imod {d_1[d_2, d_2']}$ for some $h'\in \{0, 1, \dotsc, d_1[d_2, d_2']-1\} $. Next, we write 
\[n_1-n_2=n'd_1(d_2, d_2') \quad \text{and}\quad  m=h'+ m'd_1[d_2, d_2']\]
so that $|n'|< 4N/d_1(d_2, d_2')$ and $m'\ll 1+M/d_1[d_2, d_2']$. This implies that
\begin{align}
\notag m(n_1-n_2)=h'n'd_1(d_2, d_2') + d_1^2d_2d_2'n'm'. 
\end{align}
Then, we have
 \begin{align}
 \notag \lvert\mathcal{S}_{\textup{Type II}}\rvert^2 \ll &~D_1(\log X)^{h_1^2-1}\lVert {\alpha}_1\rVert_2^2 \sum_{\substack{d_1\sim D_1}}\sum_{\substack{d_2, d_2^\prime \sim D_2\\(d_2d_2^\prime, d_1)=1}} \sum_{n_1\sim N}|{\alpha}_2(n_1)|^2\\\notag &\times \sum_{|n^\prime|< 4N/d_1(d_2, d_2^\prime)}\bigg|\sum_{m^\prime}e(m^\prime n^\prime d_1^2d_2d_2^\prime \theta)\bigg|\\
\notag \ll &~D_1(\log X)^{h_1^2-1}\lVert {\alpha}_1\rVert_2^2\lVert {\alpha}_2 \rVert_2^2\\
\label{meth2i}  &\times \sum_{\substack{d_1\sim D_1}}\sum_{\substack{d_2, d_2^\prime \sim D_2\\(d_2d_2^\prime,d_1)=1}}\sum_{|n^\prime|< 4N/d_1(d_2, d_2^\prime)}\min\bigg(\dfrac{M}{d_1[d_2, d_2^\prime]}+1, \dfrac{1}{\|n^\prime d_1^2d_2d_2^\prime\theta\|}\bigg).
 \end{align}
The terms with $n^\prime =0$ in \eqref{meth2i} contribute
 \begin{align}
\notag  &\leq D_1(\log X)^{h_1^2-1}\lVert {\alpha}_1\rVert_2^2\lVert {\alpha}_2 \rVert_2^2 \sum_{\substack{d_1\sim D_1}}\sum_{\substack{d_2, d_2^\prime \sim D_2\\(d_2d_2^\prime, d_1)=1}}\bigg(\dfrac{M}{d_1[d_2,d_2^\prime]} + 1\bigg)\\
\notag & \ll D_1M\lVert {\alpha}_1 \rVert_2^2\lVert {\alpha}_2\rVert_2^2(\log X)^{h_1^2+2} + (D_1D_2)^2\|\alpha_1\|_2^2 \|\alpha_2\|_2^2(\log X)^{h_1^2-1},
 \end{align}
using the fact that $\sum_{h_1, h_2\leq y}1/[h_1, h_2] \ll (\log y)^3$ for any $y\geq 2$. Therefore,
 \begin{align}
\notag \lvert \mathcal{S}_{\textup{Type II}}\rvert^2\ll &~\lVert {\alpha}_1\rVert_2^2\lVert {\alpha}_2 \rVert_2^2(\log X)^{h_1^2-1} \bigg\{MD_1(\log X)^{3} + D_1^2D_2^2\\
\label{meth2iii}  & + D_1\sum_{\substack{d_1\sim D_1}}\sum_{\substack{d_2, d_2^\prime \sim D_2\\(d_2d_2^\prime,d_1)=1}}\sum_{1\leq |n^\prime| < 4N/d_1(d_2, d_2^\prime)}\min\bigg(\dfrac{M}{d_1[d_2, d_2^\prime]} + 1, \dfrac{1}{\lVert n^\prime d_1^2d_2d_2^\prime\theta\rVert}\bigg)\bigg\}.
 \end{align}
 Next, we write $n^\prime d_2d_2^\prime =d''$, so that 
\[0<|d''|=|n^\prime| d_2d_2^\prime = |n^\prime| (d_2,d_2^\prime)[d_2, d_2^\prime] < \dfrac{4D_2^2N}{D_1},\] 
since $0<|n^\prime| < 4N/d_1(d_2, d_2^\prime)$ and $d_1\sim D_1$. Moreover,
\[\dfrac{M}{d_1[d_2, d_2^\prime]} = \dfrac{Md_1(d_2, d_2^\prime)|n^\prime|}{d_1^2 |n^\prime| d_2d_2^\prime} \ll  \dfrac{MN}{d_1^2|d''|}.\]
 The above reduction yields
 \begin{align}
\notag \lvert \mathcal{S}_{\textup{Type II}}\rvert^2\ll &~\|\alpha_1\|_2^2\|\alpha_2\|_2^2\bigg\{(D_1M+D_1^2D_2^2)(\log X)^{h_1^2+2}+ \\
\label{meth2ii} & + D_1(\log X)^{h_1^2-1}\sum_{d_1\sim D_1}\sum_{1\leq |d''|\ll D_2^2N/D_1}\uptau_3(d'')\min\bigg(\dfrac{MN}{d_1^2|d''|} + 1,\dfrac{1}{\lVert d_1^2d''\theta\rVert}\bigg)\bigg\}.
 \end{align}
We observe that if $D_2^2N/D_1\ll 1$, then we can bound the sum 
\[\sum_{d_1\sim D_1}\sum_{1\leq |d''|\ll D_2^2N/D_1}\uptau_3(d'')\min\bigg(\dfrac{MN}{d_1^2|d''|} + 1,\dfrac{1}{\lVert d_1^2d''\theta\rVert}\bigg) \ll \sum_{d_1\sim D_1}\bigg(\dfrac{MN}{d_1^2} + 1\bigg) \ll M + D_1. \]
Therefore, we can assume that $D_2^2N/D_1\gg1$, otherwise the sum over $d''$  in \eqref{meth2ii} can be bounded trivially as above. Without loss of generality, we can assume that $d''>0$ in the above sum.

Next, we apply Lemma \ref{Mikawa} with $\psi=1/4$ and recalling that $MN<X$ to obtain
\begin{align}
\notag D_1\sum_{d_1\sim D_1}&\sum_{0<d''\ll D_2^2N/D_1}\uptau_3(d'')\min\bigg(\dfrac{MN}{d_1^2d''} + 1,\dfrac{1}{\lVert d_1^2d''\theta\rVert}\bigg)\\
\notag \ll &~(\log X)\max_{1\leq J\leq D_2^2N/D_1}\bigg|D_1\sum_{d_1\sim D_1}\sum_{d''\sim J}\uptau_3(d'')\min\bigg(\dfrac{X}{d_1^2d''} + 1,\dfrac{1}{\lVert d_1^2d''\theta\rVert}\bigg)\bigg|\\
\notag \ll &~\bigg\{D_1D_2^2N + X\bigg(\dfrac{1}{D_1} + \dfrac{qH}{X} + \dfrac{1}{qH}\bigg)^{1/4}\bigg\}(\log X)^9.
\end{align}
Hence, from the above estimate together with \eqref{meth2ii}, and recalling from \eqref{lemmapeq2} that $\lVert {\alpha}_1\rVert_2\lVert {\alpha}_2 \rVert_2 \ll X^{1/2}(\log X)^{h^2+1}$, we obtain
\begin{align}
\notag \lvert \mathcal{S}_{\textup{Type II}}\rvert \ll &~ X^{1/2}\bigg(MD_1 + D_1^2D_2^2 + D_1D_2^2N+X\bigg(\dfrac{1}{D_1} + \dfrac{qH}{X} + \dfrac{1}{qH}\bigg)^{1/4}\bigg)^{1/2}\\
\notag &\times (\log X)^{h^2+h_1^2/2+5}.
\end{align}
The above estimate on simplification completes the proof of the lemma.
\end{proof}

Let us now combine Lemma \ref{lemmatypeIIpointwise} and Lemma \ref{lemmatypeII} to obtain the following special case of Type II sums. In particular, we will use an optimization idea due to Mikawa \cite{Mik}.

\begin{cor}\label{Cor: D}
Let $D, M, N\geq 1$ be such that 
\[DM<X, \quad N\leq M \quad \text{and} \quad MN<X.\] Let $\theta=a/q+\beta$ with $(a, q)=1$ and $|\beta|<1/q^2$. Suppose ${\alpha}_1$ and ${\alpha}_2$ are two arithmetic functions supported in $[M, 2M]$ and $[N, 2N]$, respectively, and satisfy $|{\alpha}_1|, |{\alpha}_2|\leq \uptau_h\cdot  \log$ for some fixed integer $h\geq 1$. Furthermore, let $H=1+|\beta|X$ and $qH\in [1, X]$. Then for any integer $h_1\geq 1$, we have
\begin{align}
\notag \sum_{\substack{d\sim D}}\uptau_{h_1}(d)\cdot\max_{(c, d)=1}&\Bigg\lvert\sum_{\substack{mn <X\\m\sim M, n\sim N\\mn\equiv c\imod d}}{\alpha}_1(m){\alpha}_2(n)e(mn\theta)\Bigg\rvert\\
\notag &\ll X\bigg(\dfrac{DM}{X} + \dfrac{D^2}{X} + \dfrac{M^{1/9}}{X^{1/9}}+ \dfrac{(qH)^{1/9}}{X^{1/9}} + \dfrac{1}{(qH)^{1/9}}\bigg)^{1/2}(\log X)^{h^2+h_1^2/2 + 5}.
\end{align}
\end{cor}

\begin{proof}
Let $\Sigma_1$ be the sum we wish to estimate in the corollary. Then, by Lemma \ref{lemmatypeIIpointwise} and the fact that $N\leq M$, we have
\begin{align}
\label{cor6.3a}\Sigma_1
&\ll DX\bigg(\dfrac{M}{X} + \dfrac{qH}{X} + \dfrac{1}{qH}\bigg)^{1/2}(\log X)^{h^2 + h_1+1}.
\end{align}
Next, we apply Lemma \ref{lemmatypeII} with $D_1=D$ and $D_2=1$ along with the fact that $N\leq M$ to obtain
\begin{align}\label{cor6.3b}
\Sigma_1 \ll  X\bigg(\dfrac{DM}{X} + \dfrac{D^2}{X} + \dfrac{1}{D^{1/4}} + \dfrac{(qH)^{1/4}}{X^{1/4}} + \dfrac{1}{(qH)^{1/4}}\bigg)^{1/2}(\log X)^{h^2+h_1^2/2+5}.
\end{align}
From the inequalities \eqref{cor6.3a} and \eqref{cor6.3b}, we obtain
\begin{align}\label{cor6.3c}
\Sigma_1^2  \ll &~X^2\bigg\{\dfrac{DM}{X} + \dfrac{D^2}{X} +\dfrac{(qH)^{1/4}}{X^{1/4}} + \dfrac{1}{(qH)^{1/4}} + \min\bigg(\dfrac{1}{D^{1/4}}, \dfrac{D^2M}{X} + \dfrac{D^2qH}{X} + \dfrac{D^2}{qH}\bigg)\bigg\}\\
\notag &\times (\log X)^{2h^2+h_1^2+10}.
\end{align}
Next, we have
\begin{align}
\notag \min\bigg(\dfrac{1}{D^{1/4}}, \dfrac{D^2M}{X} + \dfrac{D^2qH}{X} + \dfrac{D^2}{qH}\bigg) &\leq \bigg(\dfrac{1}{D^{1/4}}\bigg)^{8/9}\cdot \bigg(\dfrac{D^2M}{X} + \dfrac{D^2qH}{X} + \dfrac{D^2}{qH}\bigg)^{1/9}\\
\notag &=\bigg(\dfrac{M}{X} + \dfrac{qH}{X} + \dfrac{1}{qH}\bigg)^{1/9}.
\end{align}
Finally, we substitute the above estimate in \eqref{cor6.3c} along with the fact that $qH\in [1, X]$ to complete the proof of the corollary.
\end{proof}

Now we combine Lemma \ref{lemmatypeII} and Corollary \ref{Cor: D} to deduce the following corollary.

\begin{cor}\label{Cor: D_1D_2ii}
Let $D_1, D_2, M, N\geq 1$ be such that 
\[D_1M<X, \quad N\leq M, \quad \text{and} \quad MN<X.\] Suppose that ${\alpha}_1$ and ${\alpha}_2$ are two arithmetic functions supported in $[M, 2M]$ and $[N, 2N]$, respectively, and satisfy $|{\alpha}_1|, |{\alpha}_2|\leq \uptau_h\cdot  \log$ for some fixed integer $h\geq 1$. Let $\theta=a/q+\beta$ with $(a, q)=1$ and $|\beta|<1/q^2$. Furthermore, let $H=1+|\beta|X$ and $qH\in [1, X]$. Set
\begin{align}
\notag \mathcal{S}:= \notag \mathop{\sum\sum}_{\substack{d_1\sim D_1\\d_2\sim D_2\\(d_1d_2, c)=1\\(d_1,d_2)=1}}&\Bigg\lvert\sum_{\substack{mn <X\\m\sim M, n\sim N\\mn\equiv c\imod {d_1d_2}}}{\alpha}_1(m){\alpha}_2(n)e(mn\theta)\Bigg\rvert 
\end{align}

Then the following estimates hold.
\begin{enumerate}[(a)]
\item If $D_1D_2^2N<X$, we have
\begin{align}
\notag \mathcal{S}\ll &~X\bigg(\dfrac{D_1M}{X} + \dfrac{(D_1D_2)^2}{X} + \dfrac{D_1D_2^2N}{X} + \dfrac{M^{1/9}}{X^{1/9}}+ \dfrac{(D_2M)^{1/5}}{X^{1/5}}+\dfrac{(qH)^{1/9}}{X^{1/9}} + \dfrac{1}{(qH)^{1/9}}\bigg)^{1/2}\\
\notag &\times (\log X)^{h^2+7}.
\end{align}
\item If $D_1D_2^{3/2}< X^{1/2}$, we have
\begin{align}
\notag \mathcal{S}\ll &~X\bigg(\dfrac{D_1M}{X} + \dfrac{(D_1D_2)^2}{X} + \dfrac{D_1D_2^{3/2}}{X^{1/2}} + \dfrac{M^{1/9}}{X^{1/9}}+ \dfrac{(D_2M)^{1/5}}{X^{1/5}}+\dfrac{(qH)^{1/9}}{X^{1/9}} + \dfrac{1}{(qH)^{1/9}}\bigg)^{1/2}\\
\notag &\times (\log X)^{h^2+7}.
\end{align}
\end{enumerate}

\end{cor}

\begin{proof}
 We apply Lemma \ref{lemmatypeII} with $h_1=1$ and $c'=c$ to obtain
\begin{align}\label{cor6.4i}
\mathcal{S} \ll X\bigg(\dfrac{D_1M}{X} + \dfrac{(D_1D_2)^2}{X} + \dfrac{D_1D_2^2N}{X} + \dfrac{1}{D_1^{1/4}} + \dfrac{(qH)^{1/4}}{X^{1/4}} + \dfrac{1}{(qH)^{1/4}}\bigg)^{1/2}(\log X)^{h^2+6}.
\end{align}
Next, we write $d=d_1d_2$, so that $d\in [D_1D_2, 4D_1D_2]$. We then apply Corollary \ref{Cor: D} with $D=D_1D_2$ and $h_1=2$ to obtain
\begin{align}\label{cor6.4ii}
\mathcal{S} \ll X\bigg(\dfrac{D_1D_2M}{X} + \dfrac{(D_1D_2)^2}{X} + \dfrac{M^{1/9}}{X^{1/9}} + \dfrac{(qH)^{1/9}}{X^{1/9}} + \dfrac{1}{(qH)^{1/9}}\bigg)^{1/2}(\log X)^{h^2+7}.
\end{align}
From the inequalities \eqref{cor6.4i} and \eqref{cor6.4ii}, we obtain
\begin{align}
\notag|\mathcal{S}|^2 \ll X^2\bigg\{\dfrac{D_1M}{X} + \dfrac{(D_1D_2)^2}{X}&+ \dfrac{M^{1/9}}{X^{1/9}} + \dfrac{(qH)^{1/9}}{X^{1/9}} + \dfrac{1}{(qH)^{1/9}} \\
\label{cor6.4iii}& + \min\bigg(\dfrac{D_1D_2^2N}{X} + \dfrac{1}{D_1^{1/4}}, \dfrac{D_1D_2M}{X}\bigg)\bigg\}(\log X)^{2h^2+14},
\end{align}
where we have used the fact that $qH \in [1, X]$. Now we optimize the right-hand side of the above expression to obtain
\begin{align}\notag 
 \min\bigg(\dfrac{D_1D_2^2N}{X} + \dfrac{1}{D_1^{1/4}}, \dfrac{D_1D_2M}{X}\bigg)  &\ll \dfrac{D_1D_2^2N}{X} + \bigg(\dfrac{1}{D_1^{1/4}}\bigg)^{4/5}\bigg(\dfrac{D_1D_2M}{X}\bigg)^{1/5}\\
 \notag &\ll \dfrac{D_1D_2^2N}{X} + \dfrac{(D_2M)^{1/5}}{X^{1/5}}.
\end{align} Substituting the above estimate in \eqref{cor6.4iii} completes the proof of the part (a) of the corollary. Next, we note that
\begin{align}\notag 
 \min\bigg(\dfrac{D_1D_2^2N}{X} + \dfrac{1}{D_1^{1/4}}, \dfrac{D_1D_2M}{X}\bigg)  &\ll \bigg(\dfrac{D_1D_2^2N}{X}\cdot \dfrac{D_1D_2M}{X}\bigg)^{1/2} + \bigg(\dfrac{1}{D_1^{1/4}}\bigg)^{4/5}\bigg(\dfrac{D_1D_2M}{X}\bigg)^{1/5}\\
 \notag &\ll \dfrac{D_1D_2^{3/2}}{X^{1/2}} + \dfrac{(D_2M)^{1/5}}{X^{1/5}}.
\end{align}
The above estimate together with \eqref{cor6.4iii} completes the proof of part (b).
\end{proof}

\section{Exponential sums over primes in arithmetic progressions}\label{sec: exponen}

\subsection{A general exponential sum estimate over primes in arithmetic progressions}
We consider a general exponential sum estimate. Our key aim is to reduce the exponential sum over primes in arithmetic progressions into estimating Type I and Type II sums via the Vaughan identity.

\begin{prop}[General exponential sum over primes in arithmetic progressions]\label{ExpGen}
Let $\delta>0$ be small and let $b$ be a fixed positive integer. Suppose that $\sigma$ is an arithmetic function such that $\sigma$ is supported in $[1, D]$ with $D\leq X^{1/2-\delta}$, $|\sigma|\leq \uptau$, and for each $d$ in the support of $\sigma$, $c_d$ is some reduced residue class modulo $d$.

Let $\theta =a/q+\beta$ with $(a, q)=1$ and $|\beta|<1/q^2$. Furthermore, let $H=1+|\beta|X$ and $qH\in [1, X]$. For any arithmetic functions $\alpha_1, \alpha_2, \alpha_3$ with $|\alpha_1|, |\alpha_2|, |\alpha_3|\leq \uptau_2\cdot \log$, suppose that the following two conditions holds.
\begin{enumerate}[(I)]
\item \label{TypeI} For $j\in \{0,1\}$, and for some constant $C_1>0$, we have
\begin{align}\notag 
\bigg|\sum_{\substack{d\leq D\\(d, b)=1}}{\sigma}(d)\sum_{\substack{mn<X\\ 1\leq m\leq X^{1/3}\\mn\equiv c_d\imod d}}{\alpha}_1(m)(\log n)^je(mn\theta)\bigg| \ll X\bigg(\dfrac{(qH)^{\delta/2}}{X^{\delta/2}} + \dfrac{1}{(qH)^{\delta/2}}\bigg)(\log X)^{C_1}.
\end{align}
\item \label{Type II} For $N\leq M$, $MN<X$, and for some constant $C_2>0$, we have
\begin{align}\notag 
\max_{\substack{D'\leq D\\ M, N\leq X^{2/3}}}\bigg|\sum_{\substack{d\sim D'\\(d, b)=1}}{\sigma}(d)\sum_{\substack{mn<X\\ m\sim M, n\sim N \\mn\equiv c_d\imod d}}{\alpha}_2(m){\alpha}_3(n)e(mn\theta)\bigg| \ll & X\bigg(\dfrac{(qH)^{\delta/2}}{X^{\delta/2}} + \dfrac{1}{(qH)^{\delta/2}}\bigg)(\log X)^{C_2}.
\end{align}
\end{enumerate}
Then, we have
\begin{align}\notag 
\bigg|\sum_{\substack{d\leq D\\(d, b)=1}}{\sigma}(d)\sum_{\substack{n<X\\n\equiv c_d\imod d\\(n, b)=1}}\Lambda(n)e(n\theta)\bigg| \ll_{b, \delta} X\bigg(\dfrac{(qH)^{\delta/2}}{X^{\delta/2}} + \dfrac{1}{(qH)^{\delta/2}}\bigg)(\log X)^{C_3},
\end{align}
where $C_3=\max\{C_1, C_2+3\}$.
\end{prop}

\begin{proof}
We may drop the condition $(n, b)=1$ in the sum. Indeed, the contribution of $(n, b)>1$ is
\begin{align}\notag
\ll \bigg|\sum_{\substack{d\leq D\\(d, b)=1}}{\sigma}(d)\sum_{\substack{n<X\\n\equiv c_d\imod d\\(n, b)>1}}\Lambda(n)e(n\theta)\bigg| & \ll D(\log D)(\log X)\uptau(b)\\
\notag & \ll_{b, \delta}X^{1/2-\delta}(\log X)^2 \ll_{b, \delta} X\cdot \dfrac{(qH)^{1/2+\delta}}{X^{1/2+\delta}}(\log X)^2,
\end{align}
which is negligible. Therefore, we can focus on bounding the following sum
\begin{align}\label{Gen Sigma}
\Sigma:=\sum_{\substack{d\leq D\\(d,b)=1}}{\sigma}(d)\sum_{\substack{n<X\\n\equiv c_d\imod d}}\Lambda(n)e(n\theta).
\end{align}
Let $U=X^{1/3}$. Then, by Vaughan's identity (see \cite[Lemma 23.1]{kou}), we have
\begin{align}\notag
\Lambda(n) &= \Lambda_{\leq U}(n) + (\mu_{\leq U}*\log)(n) - (\mathfrak{f}_{\leq U}*1)(n) - (\mathfrak{f}_{> U}*1)(n)+ (\mu_{>U}*\Lambda_{>U}*1)(n),
\end{align}
where $\mathfrak{f} = \mu_{\leq U}*\Lambda_{\leq U}$ and note that $|\mathfrak{f}|\leq \log $. This allows us to write the sum in \eqref{Gen Sigma} as 

\begin{align}
\notag \Sigma =& \sum_{\substack{d\leq D\\(d, b)=1}}\sigma(d)\sum_{\substack{n<X\\n\equiv c_d\imod d}}\bigg(\Lambda_{\leq U}(n) + (\mu_{\leq U}*\log)(n) - (\mathfrak{f}_{\leq U}*1)(n)\\
\notag & - (\mathfrak{f}_{> U}*1)(n)+ (\mu_{>U}*\Lambda_{>U}*1)(n)\bigg)e(n\theta)\\
\label{Gen Sigma i} =&~\Sigma_1 + \Sigma_2 - \Sigma_3 - \Sigma_4 + \Sigma_5,
\end{align}
say.

Since $\Lambda \leq \log$, we can bound the sum $\Sigma_1$ as
\begin{align}\label{Gen Sigma ii}
\Sigma_1 
 &\ll (\log X)\sum_{\substack{d\leq D\\(d, b)=1}}\uptau(d)\bigg(\dfrac{X^{1/3}}{d} + 1\bigg)\ll (X^{1/3} + D)(\log X)^3 \ll X^{1/2} \ll \dfrac{X (qH)^{\delta/2}}{X^{\delta/2}}.
\end{align}

Next, we estimate the sums $\Sigma_2$ and $\Sigma_3$ using condition \ref{TypeI} with $\alpha_1\in \{\mu_{\leq X^{1/3}}, \mathfrak{f}_{\leq X^{1/3}}\}$ to obtain
\begin{align}\label{Gen Sigma iii}
\Sigma_2, \Sigma_3 \ll X\bigg(\dfrac{(qH)^{\delta/2}}{X^{\delta/2}} + \dfrac{1}{(qH)^{\delta/2}}\bigg)(\log X)^{C_1}.
\end{align}

Now we estimate the sum $\Sigma_4$ given by
\begin{align}\notag 
\Sigma_4 = \sum_{\substack{d\leq D\\(d, b)=1}}\sigma(d)\sum_{\substack{mn<X\\mn\equiv c_d\imod d\\ X^{1/3} \leq m\leq X^{2/3}}}\mathfrak{f}(m)e(mn\theta).
\end{align}
By a dyadic decomposition of summation ranges, we find that
\begin{align}\notag 
\Sigma_4\ll (\log X)^3 \max_{1\leq D'\leq D}\max_{X^{1/3}\leq M\leq X^{2/3}}\max_{1\leq N\leq X^{2/3}}\bigg|\sum_{\substack{d\sim D'\\(d, b)=1}}\sigma(d)\sum_{\substack{mn<X\\mn\equiv c_d\imod d\\ m\sim M, n\sim N\\MN<X}}\mathfrak{f}(m)e(mn\theta)\bigg|.
\end{align}
Now we can apply condition \ref{Type II} with $\{\alpha_2, \alpha_3\}=\{\mathfrak{f}_{>X^{1/3}}, 1\}$ by considering whether $M$ or $N$ is longer or not. The key point is that both $M, N\leq X^{2/3}$. Therefore, we obtain
\begin{align}
\label{Gen Sigma iv}\Sigma_4 \ll X\bigg(\dfrac{(qH)^{\delta/2}}{X^{\delta/2}} + \dfrac{1}{(qH)^{\delta/2}}\bigg)(\log X)^{C_2+3}.
\end{align}

Similarly, by a dyadic decomposition of summation ranges in $\Sigma_5$, we have
\begin{align}\notag
\Sigma_5 \ll (\log X)^3 \max_{1\leq D'\leq D}\max_{X^{1/3}\leq M, N\leq X^{2/3}}\bigg|\sum_{\substack{d\sim D'\\(d, b)=1}}\sigma(d)\sum_{\substack{mn<X\\mn\equiv c_d\imod d\\ m\sim M, n\sim N\\MN<X}}\Lambda(m) (\mu_{> X^{1/3}}*1)(n)e(mn\theta)\bigg|.
\end{align}
Since both $M, N\in [X^{1/3}, X^{2/3}]$, without the loss of generality we can assume $N\leq M$ and apply condition \ref{Type II} with $\alpha_2=\Lambda_{>X^{1/3}}$ and $\alpha_3= \mu_{>X^{1/3}}*1$ to obtain
\begin{align}
\label{Gen Sigma v}\Sigma_5 \ll X\bigg(\dfrac{(qH)^{\delta/2}}{X^{\delta/2}} + \dfrac{1}{(qH)^{\delta/2}}\bigg)(\log X)^{C_2+3}.
\end{align}
Hence, substituting the estimates from \eqref{Gen Sigma ii}, \eqref{Gen Sigma iii}, \eqref{Gen Sigma iv}, \eqref{Gen Sigma v} in \eqref{Gen Sigma i} completes the proof of the proposition.
\end{proof}

\begin{rem}\label{Delta Remark}
Note that if $\delta>0$ small, $\delta_1\in \{\delta, \delta/2\}$, and $\delta_2, \delta_3\geq \delta_1$, then we have the following estimate
\begin{align}\label{Delta remark eq 1}
\dfrac{1}{X^{\delta_1}} + \dfrac{1}{X^{\delta_2}} + \dfrac{(qH)^{\delta_3}}{X^{\delta_3}} + \dfrac{1}{(qH)^{\delta_3}} \ll \dfrac{(qH)^{\delta_1}}{X^{\delta_1}} + \dfrac{1}{(qH)^{\delta_1}},
\end{align}
where $H=1+|\beta|X$ and $qH\in [1, X]$.

We will use the above estimate in several occasions in the paper.
\end{rem}

\subsection{Exponential sum estimates over primes in arithmetic progression}
We now employ Proposition \ref{ExpGen} to establish the following exponential sum estimate.

\begin{prop}[Exponential sum over primes in arithmetic progressions]\label{PropminorPrimeI}
Let $\delta>0$, let $b$ be a fixed positive integer, and let $D\leq X^{1/3-\delta}$. Let $\theta= a/q + \beta$ with $(a, q)=1$ and $\lvert \beta\rvert <1/q^2$. Furthermore, let $H=1+|\beta|X$ and $qH\in [1, X]$. Then for some constant $C_1$, we have 
\begin{equation}
\sum_{\substack{d\leq D\\(d,b)=1}}\max_{(c,d)=1}\bigg\lvert\sum_{\substack{n<X\\n \equiv c \imod d\\(n, b)=1}}\Lambda(n)e(n\theta)\bigg\rvert \ll_{b, \delta} X\bigg(\dfrac{(qH)^{\delta/2}}{X^{\delta/2}} + \dfrac{1}{(qH)^{\delta/2}}\bigg)(\log X)^{C_1}.
\end{equation}
\end{prop}

\begin{proof}
Without the loss of generality, we may assume that the maximum over $c$ is attained at $c_d$. Let $\lambda(d)\in \mathbb{C}$ be of absolute value $1$ whenever $c=c_d$ and $(d, bc_d)=1$ for $d\in [1, D]$, so that
\begin{align}
\notag \sum_{\substack{d\leq D\\(d,b)=1}}\max_{(c,d)=1}\bigg\lvert\sum_{\substack{n<X\\n \equiv c \imod d\\(n, b)=1}}\Lambda(n)e(n\theta)\bigg| =  \sum_{\substack{d\leq D\\(d, b)=1}}\lambda(d)\sum_{\substack{n<X\\n \equiv c_d \imod d\\(n, b)=1}}\Lambda(n)e(n\theta).
\end{align}
We may now use Proposition \ref{ExpGen} with $\sigma =\lambda$ to establish the required bound. Note that $|\lambda|\leq 1$ in this case. So, it is enough to estimate the Type I and Type II sums.

\noindent {\it Verifying condition \ref{TypeI}:} Recall that Type I sum in this case is of the following form,
\begin{align}\notag 
\Sigma_{\text{Type I}}:= \sum_{\substack{d\leq D\\(d,b)=1}}\max_{(c,d)=1}\bigg\lvert\sum_{\substack{n<X\\mn \equiv c \imod d\\ 1\leq m\leq X^{1/3}}}\alpha_1(m)(\log n)^je(mn\theta)\bigg|,
\end{align}
where $|\alpha_1| \leq \uptau_2\cdot \log$ and $j\in \{0, 1\}$. We apply Lemma \ref{lemmatype1} with $M=X^{1/3}$, $v=1$, $h_1=2$, $h_2=1$ and $h_3=1$ and Remark \ref{Delta Remark} to obtain
\begin{align}\notag 
\Sigma_{\text{Type I}}  &\ll X\bigg(\dfrac{(qH)^{\delta/2}}{X^{\delta/2}} + \dfrac{1}{(qH)^{\delta/2}}\bigg)(\log X)^8,
\end{align}
as desired.

\noindent {\it Verifying condition \ref{Type II}:} We wish to estimate the following Type II sum
\begin{align}\notag 
\Sigma_{\textup{Type II}}:=\max_{\substack{D'\leq D \\ M, N\leq X^{2/3}\\ N\leq M, NM<X}}\sum_{\substack{d\sim D'\\(d,b)=1}}\max_{(c, d)=1}\bigg|\sum_{\substack{mn<X\\ m\sim M, n\sim N \\mn\equiv c\imod d}}{\alpha}_2(m){\alpha}_3(n)e(mn\theta)\bigg|,
\end{align}
where $|\alpha_2|, |\alpha_3|\leq \uptau_2\cdot \log h$. Recalling that $D\leq X^{1/3-\delta}$, we may apply Corollary \ref{Cor: D} with $h=2$ and $h_1=1$ and Remark \ref{Delta Remark} to obtain
\begin{align}\notag 
\Sigma_{\textup{Type II}} &\ll X\bigg(\dfrac{(qH)^{\delta/2}}{X^{\delta/2}} + \dfrac{1}{(qH)^{\delta/2}}\bigg)(\log X)^{10}.
\end{align}
This completes the verification of condition \ref{Type II}, and hence the proof of the proposition.
\end{proof}

\subsection{Exponential sum over primes in arithmetic progressions with composite moduli}
We now establish the exponential sum over primes in arithmetic progressions with composite moduli, which is one of the key inputs to prove Theorem \ref{Main Theorem2}.

\begin{prop}[Exponential sum over primes with composite moduli]\label{PropminorPrimeII}
Let $\delta>0$ be small, and let $b$ be a fixed positive integer. Let $D_1\in [1, X^{1/3-\delta}]$ and $D_2\in [1, X^{1/9}]$. Let $\theta= a/q + \beta$ with $(a, q)=1$ and $|\beta|<1/q^2$. Furthermore, let $H=1+|\beta|X$ and $qH\in [1, X]$. Let $c$ be a fixed non-zero integer. Then for some constant $C_2>0$, we have
\begin{align}
\sum_{d_1\leq D_1}\sideset{}{^*}\sum_{\substack{d_2\leq D_2}}\bigg\lvert \sum_{\substack{n< X\\n \equiv c \imod {d_1d_2}\\(n, b)=1}}\Lambda(n)e(n\theta)\bigg\rvert \ll_{b, \delta}  X\bigg(\dfrac{(qH)^{\delta/2}}{X^{\delta/2}} + \dfrac{1}{(qH)^{\delta/2}}\bigg)(\log X)^{C_2},
 \end{align}
 where $*$ in the sum denotes the conditions $(d_1, d_2)=(d_1d_2, bc)=1$.
\end{prop}

\begin{proof}
We write
\begin{align}
\notag \sum_{d_1\leq D_1}\sideset{}{^*}\sum_{\substack{d_2\leq D_2}}\bigg\lvert \sum_{\substack{n< X\\n \equiv c \imod {d_1d_2}\\(n, b)=1}}\Lambda(n)e(n\theta)\bigg\rvert =\mathop{\sum\sum}_{\substack{d_1\leq D_1\\d_2\leq D_2}}\lambda(d_1, d_2)\sum_{\substack{n< X\\n \equiv c \imod {d_1d_2}\\(n, b)=1}}\Lambda(n)e(n\theta),
\end{align}
where $\lambda(d_1, d_2)$ is a complex number of absolute value $1$ whenever $( d_1d_2, bc)= (d_1, d_2)=1$ with $d_1\in [1, D_1]$ and $d_2\in [1, D_2]$. We now apply Proposition \ref{ExpGen} with \[{\sigma}(d)=\sum_{\substack{d_1d_2=d\\d_j\leq D_j\forall j}}\lambda(d_1, d_2),\] 
to establish the proposition. Note that $|{\sigma}| \leq \uptau$ in this case. So, it is enough to estimate the Type I and Type II sums. 

We can use Lemma \ref{lemmatype1} to estimate the Type I sums, which is similar to the proof of Proposition \ref{PropminorPrimeI}, so that condition \ref{TypeI} holds in Proposition \ref{ExpGen}. 

For Type II sums, we need to estimate the following sum
\begin{align}
\notag \Sigma_{\text{Type II}}:=\max_{\substack{D_1', D_2', M, N}}\mathop{\sum\sum}_{\substack{d_1\sim D_1'\\d_2\sim D_2'\\(d_1d_2, bc)=1\\(d_1,d_2)=1}}\bigg\lvert\sum_{\substack{mn<X\\mn\equiv c\imod {d_1d_2}\\m\sim M, n\sim N}}{\alpha}_1(m){\alpha}_2(n)e(mn\theta)\bigg\rvert,
\end{align}
where $|{\alpha}_1|, |{\alpha}_2|\leq \uptau_2 \cdot \log$ and the maximum is taken over those $D_1', D_2', M, N$ that satisfy
\begin{equation}\label{prop2eq1}
D_1'\in [1, D_1], \quad D_2'\in [1, D_2], \quad M, N\leq X^{2/3}, \quad  MN<X, \quad \text{and} \quad N\leq M.
\end{equation}

 We divide our analysis of the above Type II sum into two cases:

\noindent{\it Case 1:} Suppose that $M\leq X^{1/2}$. Then we write $d=d_1d_2$ so that $d_1d_2\in [D_1'D_2', 4D_1'D_2']$. We can now apply Corollary \ref{Cor: D} with $h=2$ and $h_1=2$ to obtain
\begin{align}\notag 
\Sigma_{\textup{Type II}} &\ll \max_{\substack{M, N\\ D\leq D_1D_2}}\sum_{\substack{d\sim D\\(d, bc)=1}}\uptau(d)\bigg\lvert\sum_{\substack{mn<X\\mn\equiv c\imod {d}\\m\sim M, n\sim N}}{\alpha}_1(m){\alpha}_2(n)e(mn\theta)\bigg\rvert\\
\notag &\ll \max_{\substack{M, N\\D\leq D_1D_2}}X\bigg(\dfrac{DM}{X} + \dfrac{D^2}{X} + \dfrac{M^{1/9}}{X^{1/9}} + \dfrac{(qH)^{1/9}}{X^{1/9}} + \dfrac{1}{(qH)^{1/9}}\bigg)^{1/2}(\log X)^{11}.
\end{align}
Note that by assumption $D_1\leq X^{1/3-\delta}$ and $D_2\leq X^{1/9}$. This implies that $D_1D_2\leq X^{4/9-\delta} \leq X^{1/2-\delta}$. Therefore,
\begin{align}\notag 
\Sigma_{\textup{Type II}}  &\ll X\bigg(\dfrac{(qH)^{\delta/2}}{X^{\delta/2}} + \dfrac{1}{(qH)^{\delta/2}}\bigg)(\log X)^{11}
\end{align}
by Remark \ref{Delta Remark}.

\noindent{\it Case 2:} Now we consider the case $M\geq X^{1/2}$. In this case, we have $N\leq X^{1/2}$. So,  applying Corollary \ref{Cor: D_1D_2ii} (b) with $h=2$, we obtain
\begin{align}\notag 
\Sigma_{\textup{Type II}} \ll &\max_{\substack{D_1', D_2', M, N}} X\bigg(\dfrac{D_1'M}{X} + \dfrac{(D_1'D_2')^2}{X} + \dfrac{D_1'(D_2')^{3/2}}{X^{1/2}} + \dfrac{M^{1/9}}{X^{1/9}} + \dfrac{(D_2'M)^{1/5}}{X^{1/5}} + \dfrac{(qH)^{1/9}}{X^{1/9}}\\
\notag &+ \dfrac{1}{(qH)^{1/9}}\bigg)^{1/2}(\log X)^{11}.
 \end{align}
Recall the relation \eqref{prop2eq1}, and note by assumption that $D_1\leq X^{1/3-\delta}$, and $D_2\leq X^{1/9}$, so that $D_1D_2\leq X^{4/9-\delta}$ and $D_1D_2^{3/2}\leq X^{1/2-\delta}$. Therefore, by Remark \ref{Delta Remark}, we have
\begin{align}
\notag \Sigma_{\textup{Type II}}&\ll X\bigg(\dfrac{(qH)^{\delta/2}}{X^{\delta/2}} + \dfrac{1}{(qH)^{\delta/2}}\bigg)(\log X)^{11}.
\end{align} 

The above two cases complete our analysis of Type II sum estimates. Hence, this completes the proof of the proposition.
\end{proof}

\subsection{Exponential sum over primes with a well-factorable function}
We now establish exponential sum over primes in arithmetic progressions weighted by a well-factorable function (see Definition \ref{def: well fac} for the notion of \emph{well-factorable}).

\begin{prop}[Well-factorable exponential sum estimate]\label{expwellfac}
Let $\delta>0$ and let $b$ be a fixed positive integer. Let $c$ be a fixed non-zero integer and let $\xi:\mathbb{N}\rightarrow\mathbb{R}$ be a well-factorable function of level $D\in [1, X^{1/2-\delta}]$ with $|\xi|\leq 1$. Let $\theta= a/q + \beta$ for some $(a, q)=1$ and $\lvert \beta\rvert <1/q^2$. Furthermore, let $H=1+|\beta|X$ and $qH\in [1, X]$. Then for some constant $C_3>0$, we have 
\begin{equation}
 \sum_{\substack{d\leq D\\(d,bc)=1}}\xi(d)\sum_{\substack{n<X\\n \equiv c \imod d\\(n, b)=1}}\Lambda(n)e(n\theta)\ll_{b, \delta}  X\bigg(\dfrac{(qH)^{\delta/2}}{X^{\delta/2}} + \dfrac{1}{(qH)^{\delta/2}}\bigg)(\log X)^{C_3}.
\end{equation}
\end{prop}

\begin{proof}
If $D\leq X^{1/3-\delta}$, we can apply Proposition \ref{PropminorPrimeI} and the fact that $|\xi|\leq 1$ to establish the required bound in the proposition with $C_3=C_1$. Therefore, we can assume that $D>X^{1/3-\delta}$ for the rest of the proof.

We will use Proposition \ref{ExpGen} with ${\sigma}(d)=\xi(d)$ for $d\in (X^{1/3-\delta}, X^{1/2-\delta}]$.
The calculations for the Type I sums are analogous to, as in the proof of Proposition \ref{PropminorPrimeI}. We can apply Lemma \ref{lemmatype1} to estimate the Type I sum, so that condition \ref{TypeI} holds in Proposition \ref{ExpGen}. The key difference is the estimate for the Type II sums. So, we will explain the Type II sum estimates in this case.
In order to do that, we must estimate the following Type II sum:
\begin{align}\notag 
\Sigma_{\text{well-fac,Type II}}:=\max_{\substack{D', M, N}}\bigg|\sum_{\substack{d\sim D'\\(d, bc)=1}}\xi(d)\sum_{\substack{mn<X\\mn\equiv c\imod {d_1d_2}\\m\sim M, n\sim N}}{\alpha}_1(m){\alpha}_2(n)e(mn\theta)\bigg|,
\end{align}
where $|{\alpha}_1|, |{\alpha}_2|\leq \uptau_2\cdot\log$ and the maximum is taken over those $D'$, $M$ and $N$ that satisfy
\[D'\in (X^{1/3-\delta}, X^{1/2-\delta}], \quad M, N\leq X^{2/3},\quad N\leq M, \quad  MN<X.\]
 As in the proof of Proposition \ref{PropminorPrimeII}, we divide the analysis of $\Sigma_{\text{well-fac, Type II}}$ into two cases:

\paragraph{\it Case 1:} Suppose that $M\leq X^{1/2}$. We apply Corollary \ref{Cor: D} with $h=2$, $h_1=1$, and the fact that $|\xi|\leq 1$ to obtain
\begin{align}\notag 
\Sigma_{\text{well-fac,Type II}} \ll X\bigg(\dfrac{(qH)^{\delta/2}}{X^{\delta/2}} + \dfrac{1}{(qH)^{\delta/2}}\bigg)(\log X)^{10}.
\end{align}

\paragraph{\it Case 2:} Suppose that $M\in [X^{1/2}, X^{2/3}]$. For any $d\sim D'$ in the support of $\xi$, we write
\[d=d_1d_2 \quad  \text{with}\quad (d_1, d_2)=1 \quad \text{for}\quad d_1\sim D_1\quad \text{and} \quad d_2\sim D_2,\]
 so that $D_1D_2\asymp D' \leq D\leq  X^{1/2-\delta}$. We take
 \[D_1=\dfrac{D'X^{1/2}}{M}.\]
Since $D'\in (X^{1/3-\delta}, X^{1/2-\delta}]$ and $M\in [X^{1/2}, X^{2/3}]$, we have $D_1\leq D'$ and $D'X^{1/2}\geq M$. Therefore,
\[\dfrac{D_1M}{X}\leq X^{-\delta},\quad \dfrac{D_1D_2^2N}{X}\leq X^{-\delta}, \quad \dfrac{D_2M}{X} \leq \dfrac{1}{X^{1/6}}.\]
Therefore, we can now apply Corollary \ref{Cor: D_1D_2ii} (a) with $h=2$ to obtain

\begin{align}
\notag &\Sigma_{\text{well-fac,Type II}}\\
\notag \ll &~(\log X)^2
\max_{\substack{D_1D_2 \asymp D'\\ M, N}}\mathop{\sum\sum}_{\substack{d_1\sim D_1\\d_2\sim D_2\\(d_1d_2, bc)=1\\(d_1,d_2)=1}}\bigg|\sum_{\substack{mn <X\\ m \sim M, n\sim N\\ mn \equiv c \imod {d_1d_2}}}{\alpha}_1(m){\alpha}_2(n)e(mn\theta)\bigg|\\
\notag \ll &~\max_{\substack{D_1D_2\asymp D'\\ M, N}}X\bigg(\dfrac{D_1M}{X} + \dfrac{(D_1D_2)^2}{X} + \dfrac{D_1D_2^2N}{X} + \dfrac{M^{1/9}}{X^{1/9}} + \dfrac{(D_2M)^{1/5}}{X^{1/5}} + \dfrac{(qH)^{1/9}}{X^{1/9}} + \dfrac{1}{(qH)^{1/9}}\bigg)^{1/2}\\
\notag &\times (\log X)^{13}\\
\notag \ll &~X\bigg(\dfrac{(qH)^{\delta/2}}{X^{\delta/2}} + \dfrac{1}{(qH)^{\delta/2}}\bigg)(\log X)^{13}.
\end{align}

The above two cases cover the entire range for the Type II sums. Therefore, condition \ref{Type II} holds in Proposition \ref{ExpGen}. Hence, this completes the proof of the proposition.
\end{proof}

\subsection{Exponential sum over primes with semi-linear sieve}

We will use Lemma \ref{wellfacsem} to estimate the exponential sum in the following proposition.

\begin{prop}[Semi-linear sieve exponential sum estimate]\label{expsemi}
Let $\varepsilon>0$ be small and let $\delta\in (0, 10^{-3}]$. Let $b$ be a fixed positive integer. Let $\lambda_{\textup{sem}}^-$ be a lower bound semi-linear sieve weights of level $D\in [2, X^{\frac{3}{7}(1-4\delta)-\varepsilon}]$, as given in Lemma \ref{fundsemi} and Lemma \ref{wellfacsem}. Let $\theta= a/q + \beta$ with $(a, q)=1$ and $\lvert \beta\rvert <1/q^2$. Furthermore, let $H=1+|\beta|X$ and $qH\in [1, X]$. Then, for some constant $C_4>0$, we have 
\begin{equation}\label{expsemieqM}
\sum_{\substack{d\leq D\\(d,2b)=1}}\lambda_{\textup{sem}}^-(d)\sum_{\substack{n<X\\n \equiv 1 \imod d\\ n\equiv 3\imod 8\\(n, b)=1}}\Lambda(n)e(n\theta)\ll_{b, \delta, \varepsilon}  X\bigg(\dfrac{(qH)^{\delta/2}}{X^{\delta/2}} + \dfrac{1}{(qH)^{\delta/2}}\bigg)(\log X)^{C_4}.
\end{equation}
\end{prop}
The above proposition is closely related to \cite[Theorem 1.5]{Ter}. In fact, we will borrow a few ideas from \cite{Ter} to establish the above proposition.

\begin{proof}
If $D\leq X^{1/10}$, the estimate in \eqref{expsemieqM} follows from Proposition \ref{PropminorPrimeI}. So, we may assume throughout the proof that $D\geq X^{1/10}$.

We now apply Proposition \ref{ExpGen} with ${\sigma}=\lambda_{\textup{sem}}^-$. In order to do that, we consider the following Type I and Type II sums:
\begin{align}
\notag \Sigma_{\textup{sem, Type I}}&:= \sum_{\substack{d\leq D\\(d, 2b)=1}}\lambda_{\textup{sem}}^-(d)\sum_{\substack{mn<X\\mn \equiv 1 \imod d\\mn\equiv 3\imod 8\\m\in [1, X^{1/3}]}}{\alpha}(m)\log^j(n)e(mn\theta),\\
\notag \Sigma_{\textup{sem, Type II}}&:= \max_{D', M, N}\bigg|\sum_{\substack{d\sim D'\\(d, 2b)=1}}\lambda_{\textup{sem}}^-(d)\sum_{\substack{mn<X\\mn \equiv 1 \imod d\\mn\equiv 3\imod 8\\m \sim M, n\sim N}}{\alpha}_1(m){\alpha}_2(n)e(mn\theta)\bigg|,
\end{align}
where $j\in \{0, 1\}$, $|{\alpha}|, |{\alpha}_1|, |{\alpha}_2|\leq \uptau_2\cdot \log$ and the maximum is over those $D', M, N$ that satisfy
\begin{equation}\label{semieq0}
D' \in [X^{1/10}, X^{\frac{3}{7}(1-4\delta)-\varepsilon}], \quad M, N\leq X^{2/3}, \quad N\leq M, \quad MN<X.
\end{equation}

First, we use Lemma \ref{lemmatype1} to estimate the type I sum with $h_1=2, h_2=1, h_3=1, M\leq X^{1/3}$ and $D\leq X^{3(1-4\delta)/7-\varepsilon}$ to obtain
\begin{align}
\notag \Sigma_{\textup{sem, Type I}} \ll X\bigg(\dfrac{(qH)^{\delta}}{X^{\delta}} + \dfrac{1}{(qH)^\delta}\bigg)(\log X)^6.
\end{align}
This implies that condition \ref{TypeI} holds in Proposition \ref{ExpGen}.

Next, by orthogonality of the Dirichlet characters $\chi_8$ modulo $8$, we have
\begin{align}
\notag |\Sigma_{\textup{sem, Type II}}|\leq \max_{D', M, N}\bigg|\sum_{\substack{d\sim D'\\(d, 2b)=1}}\lambda_{\textup{sem}}^-(d)\sum_{\substack{mn<X\\mn \equiv c \imod d\\m\sim M, n\sim N}}{\alpha}_1(m)\chi_8(m){\alpha}_2(n)\chi_8(n)e(mn\theta)\bigg|.
\end{align}

Next, we divide our analysis of the sum $\Sigma_{ \textup{sem, Type II}}$ into two cases. 

\noindent{\it Case 1:} Suppose that $M\leq X^{1/2}$. In this case, we use Corollary \ref{Cor: D} with $D=D'$, $c=1$, $h=2$, $h_1=1$, and the facts that $|\lambda_{\textup{sem}}^-|\leq 1$ and $D'\leq X^{3/7} \leq X^{1/2-2\delta}$, to obtain

\begin{align}\notag 
\Sigma_{ \textup{sem, Type II}}  & \ll  X\bigg(\dfrac{(qH)^\delta}{X^\delta} + \dfrac{1}{(qH)^\delta}\bigg)(\log X)^{10}.
\end{align}

\paragraph{\it Case 2:} Suppose that $M\in [X^{1/2}, X^{2/3}]$. The assumption on $M$ implies that $N\leq X^{1/2}$. We now consider two subcases.

\noindent{\it Case 2(a):} Suppose that $D' \in [X^{1/10}, X^{3(1-4\delta)/7-\varepsilon}]$ and $D' \leq X^{1-2\delta-\varepsilon^2}/M$. Recalling that $|\lambda_{\textup{sem}}^-|\leq 1$, and by Corollary \ref{Cor: D} with $D=D'$, $c=1$, $h=2$, $h_1=1$,  we obtain
\begin{align}
\notag \Sigma_{ \textup{sem, Type II}} & \ll \max_{D', M, N}X\bigg(\dfrac{D'M}{X}+ \dfrac{(D')^2}{X} +\dfrac{M^{1/9}}{X^{1/9}} + \dfrac{(qH)^{1/9}}{X^{1/9}} + \dfrac{1}{(qH)^{1/9}}\bigg)^{1/2}(\log X)^{10}.
\end{align}
By assumption, $D'M/X \leq X^{-2\delta-\varepsilon^2} \leq  X^{-2\delta}$, $D'\leq X^{3/7} \leq X^{1/2-2\delta}$, and $M\leq X^{2/3}$, so by Remark \ref{Delta Remark} we have
\begin{align}
 \notag \Sigma_{ \textup{sem, Type II}} &\ll  X\bigg(\dfrac{(qH)^\delta}{X^\delta} + \dfrac{1}{(qH)^\delta}\bigg)(\log X)^{10}.
\end{align}

\paragraph{\it Case 2(b):} Finally, we consider the case when $D' \in [X^{1/10}, X^{3(1-4\delta)/7-\varepsilon}]$ and $D' >X^{1-2\delta-\varepsilon^2}/M$. 
 
Note that the sifting parameter associated with $\lambda_{\textup{sem}}^-$ is $\leq X^{1/3-2\delta-2\varepsilon^2}$. We fix a parameter $D_0 \in [X^{1/3-2\delta-2\varepsilon^2}, X^{3(1-4\delta)/7-\varepsilon}]$ to be chosen shortly. Then any $d\sim D'$ in the support of $\lambda_{\textup{sem}}^-$ can be written as $d=d_1d_2$ with $d_1\in [X^{1/10}, D_0]$ and $d_1d_2^2\leq X^{1-4\delta-2\varepsilon^2}/D_0$. 
 
We take $D_0=X^{1-2\delta-\varepsilon^2}/M$. Note that since $M\in [X^{1/2}, X^{2/3}]$, this implies that $D_0\geq X^{1/3-2\delta-\varepsilon^2}$ and by assumption, $D_0=X^{1-2\delta-\varepsilon^2}/M < D'\leq X^{3(1-4\delta)/7 -\varepsilon}$, so Lemma \ref{wellfacsem} is applicable in this case. Next, we perform a dyadic decomposition of the range of the variables $d_1$ and $d_2$, so that
\begin{align}
\notag &d_1\sim D_1,\: d_2\sim D_2, \quad \text{where}\quad X^{1/10} \ll D_1\leq D_0, \: \quad D_1D_2^2 \leq \dfrac{X^{1-4\delta-2\varepsilon^2}}{D_0}, \quad D_1D_2\asymp D'.
 \end{align}
Therefore, we have 
\begin{align}\label{semieqn: i}
X^{1/10} \ll D_1 \leq \dfrac{X^{1-2\delta}}{M} \quad \text{and} \quad D_1D_2^2\leq \dfrac{X^{1-4\delta}}{D_0} \leq \dfrac{M}{X^{2\delta}}.
\end{align}
By Lemma \ref{lemmatypeII} with $h=2$ and $h_1=1$, we obtain
\begin{align}
\notag \Sigma_{\textup{sem, Type II}} \ll &~(\log X)^2 \max_{\substack{D_1D_2\asymp D'\\ M, N}}\mathop{\sum\sum}_{\substack{d_1\sim D_1\\d_2\sim D_2\\(d_1d_2, 2bc)=1\\(d_1,d_2)=1}}\bigg\lvert\sum_{\substack{mn <X\\ m \sim M\\ mn \equiv c \imod {d_1d_2}}}{\alpha}_1(m)\chi_8(m){\alpha}_2(n)\chi_8(n)e(mn\theta)\bigg\rvert\\
\notag \ll &~ \max_{\substack{D_1D_2\asymp D'\\ M, N}}X\bigg(\dfrac{D_1M}{X} + \dfrac{(D_1D_2)^2}{X} + \dfrac{D_1D_2^2N}{X} + \dfrac{1}{D_1^{1/4}}+ \dfrac{(qH)^{1/4}}{X^{1/4}} + \dfrac{1}{(qH)^{1/4}}\bigg)^{1/2}\\
\notag &\times (\log X)^{12}.
\end{align}
Using \eqref{semieqn: i}, recalling from \eqref{semieq0} that 
\[D_1D_2 \asymp D'\leq X^{3(1-4\delta)/7-\varepsilon}, \quad MN <X, \quad M\leq X^{2/3},\] and by Remark \ref{Delta Remark}, we have
\begin{align}
\notag \Sigma_{\textup{sem, Type II}}\ll  X\bigg(\dfrac{(qH)^{\delta}}{X^\delta} + \dfrac{1}{(qH)^{\delta}}\bigg)(\log X)^{12}.
\end{align}

The above cases cover the entire range for the Type II sums. Noting that $\delta>\delta/2$, we see that condition \ref{Type II} holds in Proposition \ref{ExpGen}. Hence, this completes the proof of the proposition.
\end{proof}

\begin{remark}
We note that our proof of Case 2(b) in Proposition \ref{expsemi} can be generalized to any well-factorable sieve weights of level $D$ as long as $D\leq X^{1/2-2\delta}$. The same idea will feature in the proof of Proposition \ref{explin}.
\end{remark}

\subsection{Exponential sum with linear sieve}
We will Lemma \ref{wellfacclin} to establish Proposition \ref{explin} given below.

\begin{prop}[Linear sieve exponential sum estimate]\label{explin}
Let $\varepsilon>0$ be small and let $\delta\in (0, 10^{-3}]$. Let $b$ be a fixed positive integer. Let $\lambda_{\textup{lin}}^+$ be an upper bound linear sieve  weights of level $D\in [2, X^{1/2-2\delta-\varepsilon}]$, as given in Lemma \ref{fundlin} and Lemma \ref{wellfacclin}. Let $L$ be a real number such that $L\in [X^{1/3-2\delta-\varepsilon}, X^{2/3+2\delta+\varepsilon}]$ and let $\mathfrak{h}$ be a bounded arithmetic real-valued function. Let $\theta= a/q + \beta$ with $(a, q)=1$ and $|\beta|<1/q^2$. Furthermore, let $H=1+|\beta|X$ and $qH\in [1, X]$. Then for some constant $C_5>0$, we have 
\begin{align}\label{linearexpestimate}
\notag \sum_{\substack{d\leq D\\(d,\: 2b)=1}}\lambda_{\textup{lin}}^+(d)\sum_{\substack{\ell\sim L\\(\ell,b)=1}}\mathfrak{h}(\ell)&\sum_{\substack{n<X/2\ell\\2\ell n+1\equiv 0\imod {d}\\ \ell n\equiv 1\imod 4\\(n,b)=1}}\Lambda(n)e\Big((2\ell n+1)\theta\Big)\\
\ll_{b, \delta, \varepsilon} &  X\bigg(\dfrac{(qH)^{\delta}}{X^\delta} + \dfrac{1}{(qH)^\delta}\bigg)(\log X)^{C_5}.
\end{align}
\end{prop}

\begin{proof}
Let $\Sigma_{\text{lin}}$ be the sum we wish to estimate. The proof is similar to the proof of Proposition \ref{expsemi}. 

We note that  $\mathfrak{h}(\ell)$ is supported on $[L, 2L)$ with $L\in [X^{1/3-2\delta-\varepsilon}, X^{2/3+2\delta+\varepsilon}]$. We can proceed in the same way as in the proof of Proposition \ref{expsemi}.

We begin with a dyadic decomposition of the range of $n$ variable, say $n\sim N$ in the sum $\Sigma_{\text{lin}}$. Note that since $L\in [X^{1/3-2\delta-\varepsilon}, X^{2/3+2\delta+\varepsilon}]$ and $n<X/2\ell$, we have that $N\leq X^{2/3+ 2\delta + \varepsilon}$. Moreover, we also have that $L\leq X^{2/3+2\delta + \varepsilon}$. 

We therefore define two new parameters $M'$ and $N'$, where 
\[M'=\max\{L, N\}, \quad N'=\min\{L, N\}, \quad \text{so that}  \quad  N'M'<X, \quad  N', M'\leq X^{2/3+2\delta+\varepsilon}.\]
We also perform a dyadic decomposition on the range of $d$ variable, say $d\sim D^\prime$, with $D^\prime \leq D$. Similarly to the proof of Proposition \ref{expsemi}, we introduce Dirichlet characters $\chi_4$ modulo $4$ to detect the congruence condition $\ell n\equiv 1\imod 4$. Therefore, we have
\begin{align}\notag
\Sigma_{\textup{lin}} \ll (\log X)^2\max_{D', M', N'}\bigg|\sum_{\substack{d\sim D^\prime\\(d, 2b)=1}}\lambda^+_{\textup{lin}}(d)\sum_{\substack{mn<X\\mn\equiv -1\imod d\\m\sim M', n\sim N'}}\alpha_1(m)\alpha_2(n)e(mn\theta)\bigg|,
\end{align}
where for $m\sim \{M', N'\}$,
\[\{\alpha_1(m), \alpha_2(m)\}=\{\mathfrak{h}(m)\chi_4(m)\cdot 1_{(m, b)=1}, \Lambda(m/2)\chi_4(m/2)\cdot 1_{(m/2, b)=1, \: 2|m}\},\]
and the maximum is over those $D', M', N'$ that satisfy
\begin{align}
\label{prop 5: eq1} D'\in [2, X^{1/2-2\delta-\varepsilon}], \quad M', N' \leq X^{2/3+2\delta+\varepsilon}, \quad N'\leq M', \quad M'N'<X.
\end{align}
Note that since $\mathfrak{h}$ is bounded and $\Lambda\leq \log$, we have $|\alpha_1|, |\alpha_2|\leq \log$. We also recall from Lemma \ref{fundlin} that $|\lambda_{\textup{lin}}^+|\leq 1$.  

Next, we divide our analysis of the above sum into three cases.

\noindent{\it Case 1:}  Suppose that $D'\leq X^{1/10}$. We can then apply Corollary \ref{Cor: D} with $D=D'$ and $h_1=h=1$ to obtain
\begin{align}\notag
\Sigma_{\textup{lin}} &\ll   X\bigg(\dfrac{(qH)^\delta}{X^\delta} + \dfrac{1}{(qH)^\delta}\bigg)(\log X)^9.
\end{align}

For the rest of the two cases, we can assume that $D'\geq X^{1/10}$.

\noindent{\it Case 2:} Suppose that $M'\leq X^{1/2}$ and $D'\in [X^{1/10}, X^{1/2-2\delta-\varepsilon}]$. In this case, we may apply Corollary \ref{Cor: D} with $D=D'$, $c=-1$, $h=h_1=1$ to obtain
\begin{align}
\notag \Sigma_{\textup{lin}} & \ll  X\bigg(\dfrac{(qH)^\delta}{X^\delta} + \dfrac{1}{(qH)^\delta}\bigg)(\log X)^9.
\end{align}

\noindent{\it Case 3:} Suppose that $M'\in [X^{1/2}, X^{2/3+2\delta+\varepsilon}]$.  The assumption on $M'$ implies that $N'\leq X^{1/2}$. We now consider two subcases.

\noindent{\it Case 3(a):} Suppose that $D^\prime \in [X^{1/10}, X^{1/2-2\delta-\varepsilon}]$ and $D^\prime \leq X^{1-2\delta-\varepsilon^2}/M'$. 
We apply Corollary \ref{Cor: D} with $D=D'$, $c=-1$, $h=h_1=1$ to obtain
\begin{align}
\notag \Sigma_{ \textup{lin}} \ll &~(\log X)^2\max_{D', M', N'} X\bigg(\dfrac{D'M'}{X} + \dfrac{(D')^2}{X} + \dfrac{(M')^{1/9}}{X^{1/9}} + \dfrac{(qH)^{1/9}}{X^{1/9}} + \dfrac{1}{(qH)^{1/9}}\bigg)^{1/2}(\log X)^7.
\end{align}
By assumption, $D^\prime M'/X \leq X^{-2\delta-\varepsilon^2} \leq X^{-2\delta}$, $M'\leq X^{2/3+2\delta+\varepsilon}$ and $D'\leq X^{1/2-2\delta-\varepsilon}$. Therefore, by Remark \ref{Delta Remark}, we see that
\begin{align}
\notag \Sigma_{ \textup{lin}}\ll  X\bigg(\dfrac{(qH)^\delta}{X^\delta} + \dfrac{1}{(qH)^\delta}\bigg)(\log X)^9.
\end{align}

\noindent{\it Case 3(b):} Finally, we consider the case when $D' \in [X^{1/10}, X^{1/2-2\delta-\varepsilon}]$ and $D' >X^{1-2\delta-\varepsilon^2}/M'$. 

If $d\sim D'$ we write $d=d_1d_2$, so that
 $d_1, d_2$ satisfy for every $D_0 \in [X^{1/5}, X^{1/2-2\delta-\varepsilon}]$, the inequalities $d_1\in [X^{1/10}, D_0]$ and $d_1d_2^2 \leq X^{1-4\delta-2\varepsilon^2}/D_0$. 

We take $D_0=X^{1-2\delta-\varepsilon^2}/M'$, which is in the range $[X^{1/5}, X^{1/2-2\delta-\varepsilon}]$ by the assumption on $D'$ and $M'$. This allows us to apply Lemma \ref{wellfacsem}. Next, we do a dyadic decomposition of the range of $d_1$ and $d_2$ variables so that 
\begin{align}
\notag & d_1\sim D_1, d_2\sim D_2, \quad \text{where}\quad X^{1/10} \ll D_1\leq D_0, \quad D_1D_2^2 \leq X^{1-4\delta-2\varepsilon^2}/D_0, \quad D_1D_2\asymp D'.
\end{align}
Therefore, we have 
\[X^{1/10} \ll D_1 \leq \dfrac{X^{1-2\delta-\varepsilon^2}}{M'} \leq \dfrac{X^{1-2\delta}}{M'} \quad \text{and}\quad D_1D_2^2\leq \dfrac{X^{1-4\delta-2\varepsilon^2}}{D_0} \leq \dfrac{ M'}{X^{2\delta}}.\]
Recalling from \eqref{prop 5: eq1} that $M'N'< X$ and $D'\leq X^{1/2-2\delta-\varepsilon}$, we can now use Lemma \ref{lemmatypeII} to obtain the desired estimate
\begin{align}
\notag \Sigma_{\textup{lin}} & \ll  \max_{\substack{D_1D_2\asymp D'\\M', N'}}X\bigg(\dfrac{D_1M'}{X} + \dfrac{(D_1D_2)^2}{X} + \dfrac{D_1D_2^2N'}{X} + \dfrac{1}{D_1^{1/4}} + \dfrac{(qH)^{1/4}}{X^{1/4}}+ \dfrac{1}{(qH)^{1/4}}\bigg)^{1/2}(\log X)^{10}\\
\notag & \ll  X\bigg(\dfrac{(qH)^\delta}{X^\delta} + \dfrac{1}{(qH)^\delta}\bigg)(\log X)^{10}.
\end{align}

The above three cases cover the entire range for the sum $\Sigma_{\textup{lin}}$ and hence, the proposition is established.
\end{proof}

\part{Circle method}\label{part: circle}
In this part of the paper, we establish Theorems \ref{Main Theorem}--\ref{Main Theorem3}, \ref{Hyp1} and \ref{Hyp2}. We will use the circle method and employ the exponential sums estimates from Part \ref{part: exp sum} to establish them.

\section{Proof of Theorems \ref{Main Theorem}--\ref{Main Theorem3}, \ref{Hyp1}, and \ref{Hyp2} }\label{mainproof}

\subsection{General Theorem}

In this section, we consider a general theorem for an arithmetic function $\mathfrak{f}$ satisfying some conditions (see Theorem \ref{Thm G}) to prove our main results.

\begin{thm}[General Theorem]\label{Thm G}
Let $\delta>0$ and let $b$ be an integer that is sufficiently large in terms of $\delta$. Let $k$ be a positive integer and set $X:=b^k$. Let $D$ be a real number such that $D\in [1, X^{1/2})$. Let $r\in \mathcal{A}\cap [0, b)$ with $(r, b)=1$ and let $s$ be a positive integer such that $(r-s, b)=1$. Let $\mathfrak{f}$ be an arithmetic function supported on integers co-prime to $b$ and $|\mathfrak{f}|\ll \log $. Suppose there exists an arithmetic function ${\sigma}$ such that ${\sigma}$ is supported on $[1, D]$, $|{\sigma}|\leq \uptau$, and
for each $d$ in the support of ${\sigma}$, $c_d$ is some reduced residue class modulo $d$. Furthermore, assume that the following three conditions hold.
\begin{enumerate}[(a)]
\item \label{Cond 1} (Partial sum estimate) For any $y\in [X^{3/4}, X]$, for any $A>0$ and for any integer $d \in [1, X)$, there exists a parameter $\lambda_d$ such that $|\lambda_d|\ll \log X$ and the relation
\begin{align}
\notag \sum_{\substack{n\leq y\\(n, d)=1}}\mathfrak{f}(n)= &y\lambda_{ d} + O_{A, b}\bigg(\dfrac{y}{(\log y)^A}\bigg)
\end{align}
holds.
\item  \label{Cond 2} (Equidistribution estimate in arithmetic progressions)  For any $A, C>0$, we have 
\begin{align}
\notag \sum_{\substack{d\leq D\\(d, b)=1}}\sum_{\substack{q\leq (\log X)^{C}\\q|X}}\max_{(c, d)=1}\max_{(m, q)=1}\max_{X^{3/4}\leq y \leq X}\Bigg|\sum_{\substack{n\leq y\\n\equiv c\imod d\\ n\equiv m\imod q }}\mathfrak{f}(n)-\dfrac{ y\lambda_{d}}{\varphi(dq)}\Bigg| \ll_{A, C, b} \dfrac{X}{(\log X)^A},
\end{align}
where $\lambda_d$ is as described in condition \ref{Cond 1}.

\item \label{Cond 3} (Exponential sum estimate)
Consider $\theta =a/q+\beta$ with $(a, q)=1$ and $|\beta|<1/q^2$. Furthermore, let $H=1+|\beta|X$ and $qH\in [1, X]$. Let $\omega$ be such that $\omega \in (0, 1)$ and $\alpha_b < \omega/2$ (where $\alpha_b$ is given by the relation \eqref{alphab}). Then there exists an absolute constant $C'>0$ such that

\begin{align}
\notag \sum_{\substack{d\leq D\\(d, b)=1}}{\sigma}(d)\sum_{\substack{n< X\\n\equiv c_d \imod d}}\mathfrak{f}(n)e(n\theta) \ll_{b, \delta, \omega
} X\bigg(\dfrac{(qH)^\omega}{X^\omega} + \dfrac{1}{(qH)^\omega}\bigg)(\log X)^{C'}.
\end{align}
\end{enumerate}
Then, for any $A>0$, we have
\begin{align}
\label{genthmeq1} \sum_{\substack{d\leq D\\(d, b)=1}}{\sigma}(d)\bigg(\sum_{\substack{n < X\\n\equiv c_d \imod d}}\mathfrak{f}(n)1_{\mathcal{A}_r}(n+s) - \dfrac{\lambda_{d}}{\varphi(d)}\dfrac{b}{\varphi(b)}\sum_{\substack{n<X}}1_{\mathcal{A}_r}(n)\bigg) \ll\dfrac{X^{\zeta}}{(\log X)^A},
\end{align}
where the implicit constant in Vinogradov's notation $\ll$ depends at most on $A$, $b$, $\delta$, and $\omega$.
\end{thm}

\begin{remark}
In Section \ref{digitsec} we will see that $\alpha_b$ given by \eqref{alphab} tends to $0$ as $b\rightarrow \infty$. So, our assumption that $\alpha_b<  \omega/2$ in condition \ref{Cond 3} of Theorem \ref{Thm G} is justified.
\end{remark}

Before embarking into the proof of Theorem \ref{Thm G}, we explain how to use it to deduce Theorems \ref{Main Theorem}-\ref{Main Theorem3}, \ref{Hyp1}, and \ref{Hyp2}.

\subsection{Proof of Theorems \ref{Main Theorem}--\ref{Main Theorem3}, \ref{Hyp1} and \ref{Hyp2}}

We begin with the proof of Theorem \ref{Main Theorem}.

\begin{proof}[Proof of Theorem \ref{Main Theorem}] We will show that for any $A>0$, 
\begin{align}
\label{abse} \sum_{\substack{d\leq X^{1/3-\delta} \\(d, b)=1}}\max_{(c, d)=1}\bigg|\sum_{\substack{n<X\\n\equiv c\imod d}}\Lambda(n)1_{\mathcal{A}_r}(n)-\dfrac{1}{\varphi(d)}\dfrac{b}{\varphi(b)}\sum_{n<X}1_{\mathcal{A}_r}(n)\bigg| \ll_{A, b, \delta} \dfrac{X^\zeta}{(\log X)^A}.
\end{align}
Without loss of generality, we can assume that $\max{(c, d)=1}$ is attained at some reduced residue class $c_d$ modulo $d$. Then,  in Theorem \ref{Thm G}, we take $\mathfrak{f}(n)=\Lambda(n)1_{(n, b)=1}$ for $n<X$, $D=X^{1/3-\delta}$, $s=0$ and ${\sigma}$ to be the corresponding sign of the expression inside the absolute value of the left-hand side of \eqref{abse} whenever $(d, bc_d)=1$. Clearly, $|{\sigma}|=1\leq \uptau$. 

Now we check the three conditions in Theorem \ref{Thm G}.

\noindent{Verifying condition \ref{Cond 1}:} The condition \ref{Cond 1} with $\lambda_d=1$ follows from the Prime Number Theorem \cite[Chapter 18]{dav} together with the fact that for any $y\geq 2$,
\begin{align}\label{thm1 eq1}
\sum_{\substack{n\leq y\\ (n, bd)>1}}\Lambda(n) \ll (\log bd)(\log y).
\end{align}

\noindent{Verifying condition \ref{Cond 2}:} In order to verify condition \ref{Cond 2}, we will show that, for any $A, C>0$, the relation
\begin{align}
\notag \sum_{\substack{d\leq D\\(d, b)=1}}\sum_{\substack{q\leq (\log X)^{C}\\q|X}}\max_{\substack{1\leq c<d\\(c, d)=1}}\max_{\substack{1\leq m<q\\(m, q)=1}}\max_{X^{3/4}\leq y \leq X}\Bigg|\sum_{\substack{n\leq y\\n\equiv c\imod d\\ n\equiv m\imod q\\(n, b)=1} }\Lambda(n)-\dfrac{ y}{\varphi(dq)}\Bigg| \ll_{A, C, b, \delta} \dfrac{X}{(\log X)^A}
\end{align}
holds. By \eqref{thm1 eq1}, we can drop the condition $(n, b)=1$ in the above sum with an admissible error of $\ll_{b}X^{1/3-\delta}(\log X)^{C+2}$. Therefore, it is enough to show that
\begin{align}
\notag \sum_{\substack{d\leq D\\(d, b)=1}}\sum_{\substack{q\leq (\log X)^{C}\\q|X}}\max_{\substack{1\leq c<d\\(c, d)=1}}\max_{\substack{1\leq m<q\\(m, q)=1}}\max_{X^{3/4}\leq y \leq X}\Bigg|\sum_{\substack{n\leq y\\n\equiv c\imod d\\ n\equiv m\imod q} }\Lambda(n)-\dfrac{ y}{\varphi(dq)}\Bigg| \ll_{A, C, b, \delta} \dfrac{X}{(\log X)^A}.
\end{align}

Since $q|X=b^k$ and $(d, b)=1$, we have that $(d, q)=1$. Without loss of generality, we can assume that the maximum over $(c, d)=1$ is attained at some reduced residue class modulo $d$, say, $c_d$ and the maximum over $(m, q)=1$ is attained at $m_q$, a reduced residue class modulo $q$. Then, by the Chinese Remainder Theorem, the system of congruences $n\equiv c_d\imod d$ and $n\equiv m_q\imod q$ has a unique solution modulo $dq$. Let us call this solution $u_{dq}$. Then, we have
\begin{align}\label{thm1 eq2}
\begin{aligned}
 \sum_{\substack{d\leq D\\(d, b)=1}}\sum_{\substack{q\leq (\log X)^{C}\\q|X}}&\max_{\substack{1\leq c<d\\(c, d)=1}}\max_{\substack{1\leq m<q\\(m, q)=1}}\max_{X^{3/4}\leq y \leq X}\Bigg|\sum_{\substack{n\leq y\\n\equiv c\imod d\\ n\equiv m\imod q\\(n, b)=1} }\Lambda(n)-\dfrac{ y}{\varphi(dq)}\Bigg|\\
 &\ll \sum_{\substack{d'\leq D(\log X)^C}}\uptau(d')\max_{X^{3/4}\leq y \leq X}\Bigg|\sum_{\substack{n\leq y\\n\equiv u_{d'}\imod {d'}} }\Lambda(n)-\dfrac{ y}{\varphi(d')}\Bigg|.
 \end{aligned}
\end{align}
Note that
\[\bigg|\sum_{\substack{n\leq y\\n\equiv u_{d'}\imod {d'}} }\Lambda(n)-\dfrac{ y}{\varphi(d')}\bigg|\ll\dfrac{y(\log y)}{d'}\] 
So, by the Cauchy-Schwarz inequality and by the Bombieri-Vinogradov Theorem \cite[Chapter 28]{dav}, the sum in \eqref{thm1 eq2} is
\begin{align}\notag 
&\ll \bigg(X(\log X)\sum_{d'\leq D(\log X)^C}\dfrac{\uptau(d')^2}{d'}\bigg)^{1/2}\bigg(\max_{X^{3/4}\leq y \leq X}\Bigg|\sum_{\substack{n\leq y\\n\equiv u_{d'}\imod {d'}} }\Lambda(n)-\dfrac{ y}{\varphi(d')}\Bigg|\bigg)^{1/2} \ll \dfrac{X}{(\log X)^A}.
\end{align}
This completes the verification of condition \ref{Cond 2}.

\noindent{\it Verifying condition \ref{Cond 3}:} Condition \ref{Cond 3} holds with $\omega=\delta/2$ and $C'=C_1$ by Proposition \ref{PropminorPrimeI}. Since $b$ is large in terms of $\delta$, we have $\alpha_b < \delta/4$. 

Thus, the estimate \eqref{genthmeq1} in Theorem \ref{Thm G} holds for $\Lambda(n)1_{(n, b)=1}$ for $n<X$. We can finally replace $\Lambda(n)1_{(n, b)=1}$ by $\Lambda(n)$ for $n\in [1, X)$ by noting that 
\[\sum_{d\leq D}\sum_{\substack{n<X\\ n\equiv c_d\imod d\\(n, b)>1}}\Lambda(n)1_{\mathcal{A}_r}(n) \ll D(\log b)(\log D) \ll_{b, \delta} X^{1/3-\delta}(\log X)\] to complete the proof of Theorem \ref{Main Theorem}.

\end{proof}

The proofs of Theorems \ref{Main Theorem2}, \ref{Main Theorem3} and \ref{Hyp1} are similar to the above proof of Theorem \ref{Main Theorem}. We will only briefly explain the key changes in the set-up.

\begin{proof}[Proof of Theorem \ref{Main Theorem2}]
We apply Theorem \ref{Thm G} with $\mathfrak{f}(n)=\Lambda(n)1_{(n, b)=1}$ for $n<X$, $s=0, c_d=c$ (a fixed reduced residue class), \[{\sigma}(d) =\sum_{\substack{d=d_1d_2\\ d_j \leq D_j \forall j}}\lambda(d_1, d_2),\]
 where $\lambda(d_1, d_2)$ is a complex number of absolute value $1$, and $D=D_1D_2$ with $D_1\leq X^{1/3-\delta}$ and $D_2\leq X^{1/9}$. Note that $|{\sigma}|\leq \uptau$ in this case. We may now check three conditions of Theorem \ref{Thm G}.
 \begin{enumerate}[(i)]
\item It is evident that by the Prime Number Theorem \cite[Chapter 18]{dav}, condition \ref{Cond 1} holds with $\lambda_{d}=1$ for any $d\in [1, X)$.
\item Condition \ref{Cond 2} follows from the Bombieri-Vinogradov Theorem \cite[Chapter 28]{dav} and the Cauchy-Schwarz inequality. 
\item Proposition \ref{PropminorPrimeII} implies condition \ref{Cond 3} with $\omega=\delta/2$ and $C'=C_2$.
\end{enumerate}
 
As noted above in the proof of Theorem \ref{Main Theorem}, we can remove the co-primality condition $(n, b)=1$ with an admissible error $\ll_{b, \delta} X^{{4/9-\delta}}(\log X)^2$. This establishes Theorem \ref{Main Theorem2}.
\end{proof}
 
\begin{proof}[Proof of Theorem \ref{Main Theorem3}]
In order to prove Theorem \ref{Main Theorem3}, we take $\mathfrak{f}(n)=\Lambda(n)1_{(n, b)=1}$ for $n<X$, $s=0, c_d=c$ (a fixed reduced residue class), ${\sigma} =\xi$ and $D=X^{1/2-\delta}$ in Theorem \ref{Thm G}. In particular, 
\begin{enumerate}[(i)]
\item condition \ref{Cond 1} follows from the Prime Number Theorem \cite[Chapter 18]{dav} with $\lambda_d=1$ for any $d\in [1, X)$,
\item condition \ref{Cond 2} follows from the Bombieri-Vinogradov Theorem \cite[Chapter 28]{dav}.,
\item Proposition \ref{expwellfac} to check condition \ref{Cond 3} with $\omega=\delta/2$ and $C'=C_3$.
\end{enumerate}

In this case also, we can extend it to $\Lambda(n)$ with an admissible error $\ll_{b, \delta}X^{1/2-\delta}(\log X)$ to deduce Theorem \ref{Main Theorem3}.
\end{proof}

\begin{proof}[Proof of Theorem \ref{Hyp1}] Theorem \ref{Hyp1} follows from Theorem \ref{Thm G} by taking 
\[\mathfrak{f}(n)=\Lambda(n)1_{n\equiv 3\imod 8}1_{(n, b)=1} \quad \text{for $n\in [1, X)$},\]$s=0, c_d=1, {\sigma} =\lambda_{\textup{sem}}^-$ and $D=X^{3(1-4\delta)/7-\varepsilon}$. Clearly, 
\begin{enumerate}[]
\item condition \ref{Cond 1} follows from the Prime Number Theorem in arithmetic progressions \cite[Chapters 20, 22]{dav} with $\lambda_{d}=1/4$,  \item condition \ref{Cond 2} follows from the Bombieri-Vinogradov Theorem, \cite[Chapter 28]{dav},
\item Proposition \ref{expsemi} implies condition \ref{Cond 3}.
\end{enumerate}

 Finally, we can replace $\Lambda(n)1_{n\equiv 3\imod 8}1_{(n, b)=1}$ by $\Lambda(n)$ with an admissible error $\ll_{b,\delta} X^{3(1-4\delta)/7}(\log X)$ to complete the proof of Theorem \ref{Hyp1}.
 \end{proof}

\begin{proof}[Proof of Theorem \ref{Hyp2}] Finally, we apply Theorem \ref{Thm G} to deduce Theorem \ref{Hyp2} by taking\[\mathfrak{f}(n)=(\mathfrak{h}*\Lambda)(n/2)1_{\substack{(n, b)=1,\: n\equiv 2\imod 8}} \quad \text{for $n\in [1, X)$},\] 
where $\mathfrak{h}$ is supported on $[L, 2L]$ with $L\in [X^{1/3-2\delta-\varepsilon}, X^{2/3+2\delta+\varepsilon}]$. Furthermore, we take $s=1, c_d=-1, {\sigma}=\lambda_{\textup{lin}}^+$ and $D=X^{1/2-2\delta-\varepsilon}$ in Theorem \ref{Thm G}. Now we check three conditions of Theorem \ref{Thm G}.
\begin{enumerate}[(i)]
\item By the Prime Number Theorem in arithmetic progressions \cite[Chapters 20, 22]{dav}, condition \ref{Cond 1} holds with $\lambda_d=\sum_{\ell\sim L,\: (\ell, 2bd)=1}\mathfrak{h}(\ell)/\ell$ for any $d\in [1, X)$. Since $\mathfrak{h}$ is bounded, we have $|\lambda_d|\ll \log L\ll\log X$.
 \item Arguing as in the proof of Theorem \ref{Main Theorem}, condition \ref{Cond 2} follows from the Bombieri-Vinogradov Theorem for the Dirichlet convolution and the Cauchy-Schwarz inequality. In particular, we apply \cite[Theorem 9.17]{opera}  with $\alpha = \mathfrak{h}$ and $\beta=\Lambda$. Note that since $\mathfrak{h}$ is supported on $[L, 2L)$ with $L\in [X^{1/3-2\delta-\varepsilon}, X^{2/3+2\delta+\varepsilon}]$, for $\ell n<y$ we have that $n<y/L$. Moreover, by Siegel-Walfisz theorem, $\Lambda$ satisfies the Siegel-Walfisz condition, and if $y\in [X^{3/4}, X]$, then we have $\ell, n <y/(\log y)^B$ for some $B$ large. Therefore, the Bombieri-Vinogradov type estimate holds for the above function $(\mathfrak{h}*\Lambda)\big(n/2\big)1_{\substack{n\equiv 2\imod 8}}$ for $n\leq y$. The Bombieri-Vinogradov type estimate together with the Cauchy-Schwarz inequality implies condition \ref{Cond 2}.
\item We can apply Proposition \ref{explin} to check condition \ref{Cond 3} with $\omega=\delta$ and $C'=C_5$ to complete the proof of Theorem \ref{Hyp2}.
\end{enumerate}
This completes the proof of Theorem \ref{Hyp2}.
\end{proof}

\subsection{Proof outline of Theorem \ref{Thm G}}
We give a brief outline of the proof of Theorem \ref{Thm G} following the set-up from Section \ref{Setup}. 

By Fourier inversion (see relation \eqref{fourinv}), we have
\begin{align}
\sum_{\substack{n < X\\n\equiv c_d \imod d}}\mathfrak{f}(n)1_{\mathcal{A}_r}(n+s) = \dfrac{1}{X}\sum_{0\leq t<X}\widehat{1}_{\mathcal{A}_r}\bigg(\dfrac{t}{X}\bigg)\widehat{\mathfrak{f}}_{d, c_d}\bigg(\dfrac{-t}{X}\bigg)e\bigg(\dfrac{-st}{X}\bigg),
\end{align}
where for any $(c, d)=(d, b)=1$ and for any real number $\theta\in [0, 1)$,
\begin{align}
\label{Fourf}\widehat{\mathfrak{f}}_{d, c}(\theta): = \sum_{\substack{n<X\\n\equiv c\imod d}}\mathfrak{f}(n)e(n\theta).
\end{align}

\begin{remark}
 Since $|\mathfrak{f}|\ll \log$, we have for any real number $\theta\in [0, 1)$ and for $d<X$, 
\begin{align}\label{Triv: for f}
\big|\widehat{\mathfrak{f}}_{d, c}(\theta)\big| \leq \sum_{\substack{n<X\\n\equiv c\imod d}}|\mathfrak{f}(n)| \ll \dfrac{X(\log X)}{d}.
\end{align}

\end{remark}

The strategy to prove Theorem \ref{Thm G} roughly goes as follows:
\begin{enumerate}[(a)]
\item As outline in Section \ref{Setup}, we dissect $t/X$ into so-called major arcs and minor arcs.
\item The major arcs contribution is estimated in Proposition \ref{PropMajArc} by employing conditions \ref{Cond 1} and \ref{Cond 2} of Theorem \ref{Thm G}.
\item The minor arcs contribution is estimated in Proposition \ref{PropMinorArc1} by using Lemma \ref{lemd2} (hybrid bound) and condition \ref{Cond 3} of Theorem \ref{Thm G}.
\item  Finally, in Section \ref{sec: proof of gen} we combine Proposition \ref{PropMajArc} (major arcs estimate) and Proposition \ref{PropMinorArc1} (minor arcs estimate) to deduce Theorem \ref{Thm G}.
\end{enumerate}

\section{Fourier estimates for the digit function}\label{digitsec}
In this section, we collect the key properties of  $\widehat{1}_{\mathcal{A}_r}$ from Maynard \cite{May2}.
For the purpose of this section, we introduce the following notation for brevity. For any integer $j\in [1, k]$ and for any real number $\theta\in [0, 1)$, we set
\begin{align}
\widehat{1}_{\mathcal{N} \cap [0,\: b^j)}(\theta):=\sum_{n<b^j}1_{\mathcal{N}}(n)e(n\theta),
\end{align}
where $\mathcal{N}=\mathcal{A}$ or $\mathcal{A}_r$. In particular, $\widehat{1}_{\mathcal{A}_r} = \widehat{1}_{\mathcal{A}_r\cap [0, b^k)}$.

We begin with the $L^1$ bound in the following lemma.
\begin{lemma}[$L^1$ bound] \label{lemd1} There exists a constant $C_b \in [1/\log b, 1+3/\log b]$ such that
\begin{align}
\notag \sup_{\vartheta\in \mathbb{R}}\sum_{0\leq t <b^k}\bigg|\widehat{1}_{\mathcal{A}_r}\bigg(\frac{t}{b^k} + \vartheta\bigg)\bigg| \ll_{b} (C_b b
\log b)^k.
\end{align}
\end{lemma}

\begin{proof}
 We write $n=\sum_{j=0}^{k-1}n_jb^j$ with $n_0=r$, so that for any real number $\theta\in [0, 1)$,
\begin{align}
\label{factorMay} \widehat{1}_{\mathcal{A}_r}(\theta)=e(r\theta)\widehat{1}_{\mathcal{A}\cap[0, b^{k-1})}(b\theta).
\end{align}
The above factorization allows us to express our sum as
\begin{align}\notag
\sup_{\vartheta\in \mathbb{R}}\sum_{0\leq t <b^k}\bigg|\widehat{1}_{\mathcal{A}_r}\bigg(\dfrac{t}{b^k}+ \vartheta\bigg)\bigg| \leq b\cdot \sup_{\vartheta\in \mathbb{R}}\sum_{0\leq t <b^{k-1}}\bigg|\widehat{1}_{\mathcal{A}}\bigg(\dfrac{t}{b^{k-1}} + b\vartheta\bigg)\bigg| \ll b(C_bb\log b)^{k-1} \ll_{b} (C_bb\log b)^k,
\end{align}
where we have used \cite[Lemma 5.1]{May2} with $b^k$ replace by $b^{k-1}$ to the sum over $t$.
\end{proof}

Next, we have the following large-sieve type estimate for the Fourier transform of the set $\mathcal{A}_r$.
\begin{lemma}[Large-sieve type estimate]\label{lemmlargsieve}
Let $Q\geq 1$. Then, we have
\begin{align}
\label{largsieveeq}\sup_{\vartheta\in \mathbb{R}}\sum_{q\sim Q}\sum_{\substack{0<a <q\\(a, q)=1}}\sup_{\lvert \varepsilon\rvert < 1/2Q^2}\bigg|\widehat{1}_{\mathcal{A}_r}\bigg(\frac{a}{q}+\varepsilon +\vartheta\bigg)\bigg| \ll_{b} (Q^2+b^k)(C_b \log b)^k,
\end{align}
where $C_b$ is the constant as in Lemma \ref{lemd1}.
\end{lemma}

\begin{proof}
Note that $a/q+\varepsilon$ with $(a, q)=1, q\sim Q$ and $|\varepsilon|<1/2Q^2$ are well-spaced by $\gg 1/Q^2$ in the interval $[0, 1]$. Therefore, by the Gallagher-Sobolev type inequality (see \cite{gal}), we have
\begin{align}\notag
\sup_{\vartheta\in \mathbb{R}}\sum_{q\sim Q}\sum_{\substack{0<a <q\\(a, q)=1}}\sup_{\lvert \varepsilon\rvert < 1/2Q^2}\bigg|\widehat{1}_{\mathcal{A}_r}\bigg(\frac{a}{q}+\varepsilon + \vartheta\bigg)\bigg| \ll Q^2\int_0^1 |\widehat{1}_{\mathcal{A}_r}(u)|\textup{d}u + \int_0^1 \bigg|\dfrac{\textup{d}\widehat{1}_{\mathcal{A}_r}(u)}{\textup{d}u}\bigg|\textup{d}u.
\end{align}
By the relation \eqref{factorMay} and arguing similarly as in the proof of \cite[Lemma 5.2]{May2}, we may estimate the above sum as
\begin{align}
\notag \sup_{\vartheta\in \mathbb{R}}\sum_{q\sim Q}\sum_{\substack{0<a <q\\(a, q)=1}}\sup_{\lvert \varepsilon\rvert < 1/2Q^2}\bigg|\widehat{1}_{\mathcal{A}_r}\bigg(\frac{a}{q}+\varepsilon + \vartheta\bigg)\bigg|& \ll Q^2\int_{0}^1|\widehat{1}_{\mathcal{A}\cap [0, b^{k-1})}(bu)|\textup{d}u + \int_0^1 \bigg|\dfrac{\textup{d}\widehat{1}_{\mathcal{A}\cap [0, b^{k-1})}(bu)}{\textup{d}u}\bigg|\textup{d}u.\\
\notag &\ll_{b} (Q^2 + b^k)(C_b\log b)^k,
\end{align}
as desired.
\end{proof}

We also have the following hybrid bound for the Fourier transform of the set $\mathcal{A}_r$.

\begin{lemma}[Hybrid estimate]\label{lemd2} Let $Q, B \geq 1$. Then, we have 
\begin{align}
\notag \sum_{q\sim Q}\sum_{\substack{1\leq a<q\\(a, q)=1}}\sum_{\substack{|\eta|<B\\b^k a/q + \eta \in \mathbb{Z}}}\bigg|\widehat{1}_{\mathcal{A}_r}\bigg(\frac{a}{q} + \frac{\eta}{b^k}\bigg)\bigg| \ll_{b} (b-1)^k(Q^2B)^{\alpha_b} + Q^2B(C_b\log b)^k,
\end{align}
where $C_b$ is the constant described in Lemma \ref{lemd1}, and
\begin{align}
\label{alphab} \alpha_b = \dfrac{\log \bigg(C_b\dfrac{b\log b}{b-1}\bigg)}{\log b}.
\end{align}
\end{lemma}
\begin{remark}
We note that $\alpha_b$ tends to $0$ as $b\rightarrow \infty$. Therefore, $\alpha_b$ will be small if we take $b$ large enough, which is a crucial point in our entire Fourier analytic set-up.
\end{remark}
\begin{proof}
The proof follows from the relation \eqref{factorMay} in combination with the arguments of \cite[Lemma 5.3]{May2}. 
\end{proof}

We end this section with the $L^\infty$ bound for $\widehat{1}_{\mathcal{A}_r}$.
\begin{lemma}[$L^\infty$ bound]\label{lemd3}
Let $q<b^{k/3}$ be of the form $q=q_1q_2$ with $(q_1, b)=1$ and $q_1\neq 1$, and let $|\varepsilon| < 1/2b^{2k/3}$. Then, for any integer $a$ with $(a, q)=1$, we have
\begin{align}
\notag \bigg|\widehat{1}_{\mathcal{A}_r}\bigg(\frac{a}{q}+\varepsilon\bigg)\bigg| \ll_{b} (b-1)^k\exp\bigg(-c_b\dfrac{k}{\log q}\bigg),
\end{align}
for some constant $c_b>0$ depending only on $b$.
\end{lemma}
\begin{proof}
The proof follows from the relation \eqref{factorMay} in conjunction with the argument of \cite[Lemma 5.4]{May2}.
\end{proof}

\section{Major arcs}\label{sec: maj}
We devote this section to establishing the major arcs estimate. Throughout, $\widehat{1}_{\mathcal{A}_r}$ denotes the Fourier transform of the set $\mathcal{A}_r$ given by \eqref{Arfourier} and $\widehat{\mathfrak{f}}_{d, c}$ is given by \eqref{Fourf}. 
\begin{prop}[Major arcs estimate for Theorem \ref{Thm G}]\label{PropMajArc}
Let $C\geq 1$ be a large real number. Assume the setting of Theorem \ref{Thm G} and recall that $s$ is a positive integer such that $(r-s, b)=1$. Then we have
\begin{align}
\notag \sum_{\substack{d\leq D\\(d, b)=1}}\max_{(c, d)=1}\bigg\lvert\dfrac{1}{X}\sum_{\substack{0\leq t<X\\t\in \mathfrak{M}}}\widehat{1}_{\mathcal{A}_r}\bigg(\dfrac{t}{X}\bigg)\widehat{\mathfrak{f}}_{d, c}\bigg(\dfrac{-t}{X}\bigg)e\bigg(\dfrac{-st}{X}\bigg)-\dfrac{\lambda_{d}}{\varphi(d)}\dfrac{b}{\varphi(b)}\sum_{\substack{n< X}}1_{\mathcal{A}_r}(n)\bigg| \ll\dfrac{X^{\zeta}}{(\log X)^{5C+5}},
\end{align}
where $\mathfrak{M}$ is given by the relation \eqref{maj}, and the implicit constant in $\ll$ depends at most on $b$, $C$ and $\delta$.

\end{prop}

Recall from the relation \eqref{maj} that $\mathfrak{M} = \mathfrak{M}_1\cup \mathfrak{M}_2\cup \mathfrak{M}_3$. In order to prove Proposition \ref{PropMajArc}, we will estimate separately the contribution coming from $\mathfrak{M}_1, \mathfrak{M_2},$ and $\mathfrak{M}_3$ in Lemmas \ref{Majorarc1}, \ref{Majorarc2}, and \ref{Majorarc3}, respectively. We begin with estimating the contribution of $\mathfrak{M}_1$ in the following lemma.

\begin{lemma}\label{Majorarc1} Let $C\geq 1$, $D\in [1, X)$, and recall the set $\mathfrak{M}_1$ is given by
\begin{align}
\notag \mathfrak{M}_1 =&\bigg\{t \in [0, X)\cap \mathbb{Z}: \bigg\lvert\dfrac{t}{X}-\dfrac{a}{q}\bigg\rvert \leq \dfrac{(\log X)^C}{X}\, \text{for some}\, (a, q)=1, 1\leq a<q \leq (\log X)^C, q\nmid X\bigg\}.
\end{align}
Assume the setting of Theorem \ref{Thm G}. Then we have
\begin{align}\label{Majorarceq1}
\dfrac{1}{X}\sum_{\substack{d\leq D\\(d, b)=1}}\max_{(c,d)=1}\bigg\lvert
\sum_{\substack{0\leq t<X\\t\in \mathfrak{M}_1}}\widehat{1}_{\mathcal{A}_r}\bigg(\dfrac{t}{X}\bigg)\widehat{\mathfrak{f}}_{d, c}\bigg(\dfrac{-t}{X}\bigg)e\bigg(\dfrac{-st}{X}\bigg)\bigg\rvert \ll \dfrac{X^{\zeta}}{(\log X)^{5C+5}},
\end{align}
where the implicit constant in $\ll$ depends at most on $b$, $C$ and $\delta$.
\end{lemma}

\begin{proof}
If $t\in \mathfrak{M}_1$, we use Lemma \ref{lemd3} to obtain
\begin{align}
\notag \bigg\lvert \widehat{1}_{\mathcal{A}_r}\bigg(\dfrac{t}{X}\bigg)\bigg\rvert \ll_{b, C} \dfrac{X^{\zeta}}{(\log X)^{8C+7}}.
\end{align}
We note that the cardinality of the set $\mathfrak{M}_1$ is at most $\ll (\log X)^{3C}$. Therefore, by relation \eqref{Triv: for f},
\begin{align}
\notag \dfrac{1}{X}\sum_{\substack{d \leq D\\(d,b)=1}}\max_{(c,d)=1}\bigg\lvert
\sum_{\substack{0\leq t<X\\t\in \mathfrak{M}_1}}\widehat{1}_{\mathcal{A}_r}\bigg(\dfrac{t}{X}\bigg)\widehat{\mathfrak{f}}_{d, c}\bigg(\dfrac{-t}{X}\bigg)\bigg\rvert & \ll_{b, C} \dfrac{1}{X}\cdot (\log X)^{3C}\cdot \dfrac{X^{\zeta}}{(\log X)^{8C+7}}\cdot X(\log X)\cdot\sum_{d\leq D}\dfrac{1}{d}\\
\notag & \ll_{b, C} \dfrac{X^{\zeta}}{(\log X)^{5C+5}}.
\end{align}
This completes the proof of the lemma.
\end{proof}

Now we estimate the contribution coming from $\mathfrak{M}_2$.

\begin{lemma}\label{Majorarc2}
Let $C\geq 1$. Recall that the set $\mathfrak{M}_2$ is given by
\begin{align}
\notag \mathfrak{M}_2 =\bigg\{&t \in [0, X)\cap \mathbb{Z}: \dfrac{t}{X}=\dfrac{a}{q}+ \dfrac{\eta}{X}\: \text{for some}\, (a, q)=1, 0\leq a <q\leq (\log X)^C,\\
\notag & \: q\geq 1, q\vert X,\: 0<\lvert \eta \rvert \leq (\log X)^C\bigg\}.
\end{align}
Assume the setting of Theorem \ref{Thm G}. Then we have
\begin{align}\label{Majorarceq2}
\dfrac{1}{X}\sum_{\substack{d\leq D\\(d, b)=1}}\max_{(c,d)=1}\bigg\lvert
\sum_{\substack{0\leq t<X\\t\in \mathfrak{M}_2}}\widehat{1}_{\mathcal{A}_r}\bigg(\dfrac{t}{X}\bigg)\widehat{\mathfrak{f}}_{d,c}\bigg(\dfrac{-t}{X}\bigg)e\bigg(\dfrac{-st}{X}\bigg)\bigg\rvert \ll\dfrac{X^{\zeta}}{(\log X)^{5C+5}},
\end{align}
where the implicit constant in $\ll$ depends at most on $b$, $C$ and $\delta$.
\end{lemma}

\begin{proof}
 We call the left-hand side of \eqref{Majorarceq2} as $\Sigma_{\textup{Major}}$ and simplify the sum as
\begin{align}\notag
\Sigma_{\textup{Major}} & \leq \dfrac{1}{X}\sum_{\substack{d\leq D\\(d,b)=1}}\max_{(c,d)=1}\sum_{\substack{q\leq (\log X)^C\\q|X}}\sum_{\substack{a=0\\(a,q)=1}}^q\sum_{0<|\eta|\leq (\log X)^C}\bigg|\widehat{1}_{\mathcal{A}_r}\bigg(\dfrac{a}{q}+\dfrac{\eta}{X}\bigg)\widehat{\mathfrak{f}}_{d,c}\bigg(\dfrac{-a}{q}-\dfrac{\eta}{X}\bigg)\bigg|.
\end{align}
In the right-hand side of the above expression, since $q|X$, we have that $\eta$ is an integer in this case. Next, we use the trivial bound $|\widehat{1}_{\mathcal{A}_r}(a/q+\eta/X)|\leq X^{\zeta}$ in the above estimate to obtain
\begin{align}\label{inner0}
\Sigma_{\textup{Major}} \leq \dfrac{X^{\zeta}}{X}\sum_{\substack{d\leq D\\(d,b)=1}}\max_{(c,d)=1}\sum_{\substack{q\leq (\log X)^C\\q|X}}\sum_{\substack{a=0\\(a,q)=1}}^q\sum_{0<|\eta|\leq (\log X)^C}\bigg|\widehat{\mathfrak{f}}_{d,c}\bigg(\dfrac{-a}{q}-\dfrac{\eta}{X}\bigg)\bigg|.
\end{align}
 Therefore, in order to establish the lemma it is enough to show that
 \begin{align}\notag
\Sigma_{\textup{Major}}':=\sum_{\substack{d\leq D\\(d,b)=1}}\max_{(c,d)=1}\sum_{\substack{q\leq (\log X)^C\\q|X}}\sum_{\substack{a=0\\(a, q)=1}}^q\sum_{0<|\eta|\leq (\log X)^C}\bigg|\widehat{\mathfrak{f}}_{d,c}\bigg(\dfrac{-a}{q}-\dfrac{\eta}{X}\bigg)\bigg| \ll_{C}\dfrac{X}{(\log X)^{5C+5}}.
 \end{align}
 We have
\begin{align}
\notag \widehat{\mathfrak{f}}_{d,c}\bigg(\dfrac{-a}{q}-\dfrac{\eta}{X}\bigg)=\sum_{\substack{m=1}}^qe\bigg(\dfrac{-ma}{q}\bigg)\sum_{\substack{n<X\\n\equiv c\imod {d}\\ n\equiv m\imod q}}\mathfrak{f}(n)e\bigg(\dfrac{-n \eta }{X}\bigg).
\end{align}
We note that since $q|X = b^k$ and $(n, b)=1$ (as $\mathfrak{f}(n)$ is supported on integers $n$ such that $(n , b)=1$), the congruence $n\equiv m\imod q$ implies that we either have $(m, q)=1$ or the above sum is empty.
Furthermore, since $(d, b)=1$ and $q|X=b^k$, we have $(d, q)=1$. Therefore, by the Chinese Remainder Theorem, the system of congruences 
\begin{align}
\notag n\equiv c\imod d \quad \text{and}\quad  n\equiv m\imod q
\end{align}
has a unique solution modulo $dq$. This allows us to write
\begin{align}
\label{inner}\widehat{\mathfrak{f}}_{d, c}\bigg(\dfrac{-a}{q}-\dfrac{\eta}{X}\bigg)=\sum_{\substack{m=1\\(m,q)=1}}^qe\bigg(\dfrac{-ma}{q}\bigg)\sum_{\substack{n<X\\n\equiv c \imod{d}\\n\equiv m\imod q}}\mathfrak{f}(n)e\bigg(\dfrac{-n \eta }{X}\bigg).
\end{align}
Next, if $(c, d)=(m, q)=1$, then for any $y\geq 2$, we denote
\begin{align}
\notag \Psi_{\mathfrak{f}}(y; d, c; q, m):=\sum_{\substack{n<y\\n\equiv c\imod {d}\\n\equiv m\imod q}}\mathfrak{f}(n).
\end{align} 
For any $y\in [X^{3/4}, X]$, we denote 
\[\Delta_{\mathfrak{f}}(y; d, c; q, m):=\Psi_{\mathfrak{f}}(y; d, c; q, m)-\frac{ y\lambda_{d}}{\varphi(dq)}.\]
Using partial summation and the inequality \eqref{Triv: for f}, we have
\begin{align}
\notag \sum_{\substack{n<X\\n\equiv c \imod{d}\\n\equiv m\imod q}}\mathfrak{f}(n)e\bigg(\dfrac{-n \eta }{X}\bigg)=& \int_{X^{3/4}}^X e\bigg(\dfrac{-y\eta}{X}\bigg)\text{d}\Delta_{\mathfrak{f}}(y; d, c; q, m)
+ \dfrac{\lambda_{d}}{\varphi(dq)}\int_{X^{3/4}}^X e\bigg(\dfrac{-y\eta}{X}\bigg)\text{d}y\\
\notag &+ O\bigg(\log X + \dfrac{X^{3/4}\log X}{dq}\bigg)\\
\label{inner01}  =& \int_{X^{3/4}}^X e\bigg(\dfrac{-y\eta}{X}\bigg)\text{d}\Delta_{\mathfrak{f}}(y; d, c; q, m) + O\bigg(\dfrac{\lambda_{d}X^{3/4}}{\varphi(dq)} + \log X + \dfrac{X^{3/4}\log X}{dq}\bigg),
\end{align}
where we have used the fact that $\eta$ is an integer, so that $\int_{1}^X e(-y\eta/X)\text{d}y=O(1)$. Next, using integration by parts, we have
\begin{align}
\label{inner02}\int_{X^{3/4}}^X e\bigg(\dfrac{-y\eta}{X}\bigg)\text{d}\Delta_{\mathfrak{f}}(y; d, c; q, m) &\ll (1+|\eta|)\max_{X^{3/4}<y\leq X}|\Delta_{\mathfrak{f}}(y; d, c; q, m)|.
\end{align}
Using the estimate from \eqref{inner02} in \eqref{inner01} along with the facts that $|\eta| \leq (\log X)^C$ and that $\lambda_d\ll \log X$, we obtain
\begin{align}
\label{inner03} \sum_{\substack{n<X\\n\equiv c \imod{d}\\n\equiv m\imod q}}\mathfrak{f}(n)e\bigg(\dfrac{-n \eta }{X}\bigg) 
 &\ll (\log X)^C\max_{X^{3/4}< y\leq X}\big|\Delta_{\mathfrak{f}}(y; d, c; q, m)\big|+ \dfrac{X^{3/4} (\log X)^2}{dq}.
\end{align}
Therefore, using the inequalities \eqref{inner} and \eqref{inner03}, we have
\begin{align}
\notag \Sigma_{\textup{Major}}' :=&\sum_{\substack{d\leq D\\(d,b)=1}}\max_{(c,d)=1}\sum_{\substack{q\leq (\log X)^C\\q|X}}\sum_{\substack{a=0\\(a,q)=1}}^q\sum_{0<|\eta|\leq (\log X)^C}\bigg|\widehat{\mathfrak{f}}_{d,c}\bigg(\dfrac{-a}{q}-\dfrac{\eta}{X}\bigg)\bigg|\\
\notag \ll &~(\log X)^{4C}\sum_{\substack{d\leq D\\(d,b)=1}}\sum_{\substack{q\leq (\log X)^C\\q|X}}\max_{(c, d)=1} \max_{(m, q)=1}\max_{X^{3/4}<y\leq X}|\Delta_{\mathfrak{f}}(y; d, c; q, m)| + X^{3/4}(\log X)^{3C+3}.
\end{align}
We can now apply \ref{Cond 2} with $A=9C+5$ to obtain
\begin{align}\notag
\Sigma_{\textup{Major}}' \ll (\log X)^{4C} \cdot \dfrac{X}{(\log X)^{9C+5}} + X^{3/4}(\log X)^{3C+3} \ll \dfrac{X}{(\log X)^{5C+5}},
\end{align}
as desired.
\end{proof}

Finally, we end this section by analyzing the set $\mathfrak{M}_3$ (given below), which yields the expected main term in Proposition \ref{PropMajArc}.  

\begin{lemma}\label{Majorarc3}
Let $C\geq 1$. Recall that the set $\mathfrak{M}_3$ is given by
\begin{align}
\notag \mathfrak{M}_3 =&\bigg\{t \in [0, X)\cap \mathbb{Z}: \dfrac{t}{X}=\dfrac{a}{q}\: \text{for some}~~(a, q)=1, 0\leq a <q\leq (\log X)^C, q\geq 1, \: q\lvert X\bigg\}.
\end{align}
Assume the setting of Theorem \ref{Thm G} and recall that $s$ is a positive integer such that $(r-s, b)=1$. Then, we have
\begin{align}
\label{majorarceqn3} \sum_{\substack{d\leq D\\(d, b)=1}}\max_{(c,d)=1}\bigg\lvert\dfrac{1}{X}
\sum_{\substack{0\leq t<X\\t\in \mathfrak{M}_3}}\widehat{1}_{\mathcal{A}_r}\bigg(\dfrac{t}{X}\bigg)\widehat{\mathfrak{f}}_{d, c}\bigg(\dfrac{-t}{X}\bigg)e\bigg(\dfrac{-st}{X}\bigg)-\dfrac{ \lambda_{d}}{\varphi(d)}\dfrac{b}{\varphi(b)}\sum_{\substack{n< X}}1_{\mathcal{A}_r}(n)\bigg\rvert  \ll \dfrac{X^{\zeta}}{(\log X)^{5C+5}},
\end{align}
where the implicit constant in $\ll$ depends at most on $b$, $C$ and $\delta$.
\end{lemma}

\begin{proof} We begin with the following observation that for $k$ large enough
\[``q\leq (\log X)^C,\: q|X=b^k" \quad \text{is equivalent to} \quad ``q\leq (\log X)^C,\: \text{for every}\: p|q, \: \text{we have} \: p|b".\]
Therefore, we have
\begin{align}
\label{Maintermeq0}\dfrac{1}{X}\sum_{\substack{0\leq t<X\\t\in \mathfrak{M}_3}}\widehat{1}_{\mathcal{A}_r}\bigg(\dfrac{t}{X}\bigg)\widehat{\mathfrak{f}}_{d, c}\bigg(\dfrac{-t}{X}\bigg)e\bigg(\dfrac{-st}{X}\bigg)=\dfrac{1}{X}\sum_{\substack{q \leq (\log X)^C\\p|q\Rightarrow p|b}}\sum_{\substack{0\leq a<q\\(a,q)=1}}\widehat{1}_{\mathcal{A}_r}\bigg(\dfrac{a}{q}\bigg)\widehat{\mathfrak{f}}_{d,c}\bigg(\dfrac{-a}{q}\bigg)e\bigg(\dfrac{-sa}{q}\bigg).
\end{align}
If $(a, q)=1$, we have
\begin{align}
\label{Maintermeq00} \widehat{\mathfrak{f}}_{d, c}\bigg(\dfrac{-a}{q}\bigg)=\sum_{\substack{m=1}}^q e\bigg(\dfrac{-ma}{q}\bigg)\sum_{\substack{n<X\\n\equiv c\imod{d}\\ n\equiv m\imod q}}\mathfrak{f}(n).
\end{align}
Arguing similarly as in Lemma \ref{Majorarc2}, we note that $(d, q)=(m, q)=1$. For brevity, let
\begin{align}
\Delta_{\mathfrak{f}}(X; dq):=\max_{(c, d)=1}\max_{(m, q )=1}\Bigg\lvert\sum_{\substack{n<X\\n\equiv c\imod d\\n\equiv m\imod q}}\mathfrak{f}(n)-\dfrac{ \lambda_{d} X}{\varphi(dq)}\Bigg\rvert.
\end{align}
Then, we have
\begin{align}
\notag  \widehat{\mathfrak{f}}_{d, c}\bigg(\dfrac{-a}{q}\bigg) &=\dfrac{\lambda_{d} X}{\varphi(dq)}\sum_{\substack{m=1\\(m,q)=1}}^q e\bigg(\dfrac{-ma}{q}\bigg)+O\Big(\varphi(q)\Delta_{\mathfrak{f}}(X;dq)\Big)\\
\label{maintermeq1} &=\dfrac{ \lambda_{d}\mu(q) X}{\varphi(dq)}+O\Big(\varphi(q)\Delta_{\mathfrak{f}}(X;dq)\Big),
\end{align}
where we have used we have the expression for the Ramanujan sum (see \cite[p. 149]{dav}): if $(a, q)=1$, then
\begin{align}
\label{Ramanujan}\sum_{\substack{m=1\\(m,q)=1}}^q e\bigg(\dfrac{-ma}{q}\bigg)=\mu(q).
\end{align}

Note that $|\widehat{1}_{\mathcal{A}_r}(a/q)| \leq X^{\zeta}$, so that the contribution of the big-Oh term from relation \eqref{maintermeq1} to the expression in  \eqref{majorarceqn3} is
\begin{align}
\notag \ll \dfrac{1}{X}\sum_{\substack{d \leq D\\(d,b)=1}}\sum_{\substack{q\leq (\log X)^C\\q|X}}\sum_{0\leq a <q}\bigg \lvert \widehat{1}_{\mathcal{A}_r}\bigg(\dfrac{a}{q}\bigg)\bigg \rvert\varphi(q)\Delta_{\mathfrak{f}}\Big(X; dq\Big)\ll &~\dfrac{X^{\zeta}(\log X)^{2C}}{X}\sum_{\substack{d\leq D\\(d,b)=1}}\sum_{\substack{q\leq (\log X)^C\\q|X}}\Delta_{\mathfrak{f}}\Big(X;dq\Big).
\end{align}
We apply condition \ref{Cond 2} with $A=7C+5$, so that the above sum is
\begin{align}
\notag & \ll_{b, C}  \dfrac{X^{\zeta}(\log X)^{2C}}{X}\cdot \dfrac{X}{(\log X)^{7C+5}} \ll_{b, C} \dfrac{X^\zeta}{(\log X)^{5C+5}},
\end{align}
which is admissible.

We are therefore left with showing that
\[\dfrac{\lambda_d}{X}\sum_{\substack{q \leq (\log X)^C\\p|q\Rightarrow p|b}}\sum_{\substack{0\leq a<q\\(a,q)=1}}\widehat{1}_{\mathcal{A}_r}\bigg(\dfrac{a}{q}\bigg)\dfrac{X\mu(q)}{\varphi(dq)}e\bigg(\dfrac{-sa}{q}\bigg) = \dfrac{\lambda_d}{\varphi(d)}\dfrac{b}{\varphi(b)}\sum_{\substack{n<X}}1_{\mathcal{A}_r}(n).\]

It is evident from the relation \eqref{maintermeq1} that $q$ is supported on square-free integers. Therefore, we have $q|b$ since for every prime $p|q$ implies $p|b$. Then, by the relation \eqref{Arfourier} and the fact that $q|b$, we have
\begin{align}
\label{maintermeq2}\widehat{1}_{\mathcal{A}_r}\bigg(\dfrac{a}{q}\bigg) =e\bigg(\dfrac{ar}{q}\bigg)\sum_{\substack{n< X}}1_{\mathcal{A}_r}(n).
\end{align}

We also note that $(d, q)=1$ as $(d, b)=1$ and $q|b$. Furthermore, since by our assumption $(r-s, b)=1$, we have $(r-s, q)=1$. Therefore, using  \eqref{maintermeq1}, \eqref{Ramanujan}, and \eqref{maintermeq2} allows us to estimate the main term as
\begin{align}
\notag \dfrac{1}{X}\sum_{\substack{q \leq (\log X)^C\\p|q\Rightarrow p|b}}\sum_{\substack{0\leq a<q\\(a,q)=1}}&\widehat{1}_{\mathcal{A}_r}\bigg(\dfrac{a}{q}\bigg)\dfrac{\lambda_{d} X\mu(q)}{\varphi(dq)}e\bigg(\dfrac{-sa}{q}\bigg)\\
\notag &= \dfrac{\lambda_{d}}{\varphi(d)}\sum_{\substack{q \leq (\log X)^C\\p|q\Rightarrow p|b}}\sum_{\substack{0\leq a<q\\(a,q)=1}}\widehat{1}_{\mathcal{A}_r}\bigg(\dfrac{a}{q}\bigg)\dfrac{\mu(q)}{\varphi(q)}e\bigg(\dfrac{-sa}{q}\bigg)\\
\notag &=\dfrac{\lambda_{d}}{\varphi(d)}\sum_{\substack{n<X}}1_{\mathcal{A}_r}(n)\sum_{\substack{q|b}}\sum_{\substack{0\leq a<q\\(a,q)=1}}e\bigg(\dfrac{a(r-s)}{q}\bigg)\dfrac{\mu(q)}{\varphi(q)}\\
\notag &=\dfrac{\lambda_{d}}{\varphi(d)}\sum_{\substack{n<X}}1_{\mathcal{A}_r}(n)\sum_{q|b}\dfrac{\mu^2(q)}{\varphi(q)}= \dfrac{ \lambda_{d}}{\varphi(d)}\dfrac{b}{\varphi(b)}\sum_{\substack{n<X}}1_{\mathcal{A}_r}(n)
\end{align}
as desired.
\end{proof}

We can now combine Lemma \ref{Majorarc1}, Lemma \ref{Majorarc2} and Lemma \ref{Majorarc3} along with the fact that $\mathfrak{M}=\mathfrak{M}_1\cup\mathfrak{M}_2\cup\mathfrak{M}_3$ (see the relation \eqref{maj}) to complete the proof of Proposition \ref{PropMajArc}.

\section{Minor arcs}\label{minar1}
In this section, we will establish the minor arcs estimate for Theorem \ref{Thm G} by combining condition \ref{Cond 3} of Theorem \ref{Thm G} and Lemma \ref{lemd2}.

\begin{prop}[Minor arcs estimate for Theorem \ref{Thm G}]\label{PropMinorArc1}
Let $Q, B\geq 1$ with $QB \ll X^{1/2}$. Assume the set-up of Theorem \ref{Thm G}. Then we have
\begin{align}
\notag  \sum_{q\sim Q}\sum_{\substack{a=1\\(a,q)=1}}^q\sum_{\substack{B< |\eta|+1\leq 2B\\X a/q+\eta \in \mathbb{Z}}}\bigg|\widehat{1}_{\mathcal{A}_r}\bigg(\frac{a}{q}+\frac{\eta}{X}\bigg)\sum_{\substack{d\leq D\\(d, b)=1}}{\sigma}(d)\widehat{\mathfrak{f}}_{d, c_d}\bigg(-\bigg(\frac{a}{q}+\frac{\eta}{X}\bigg)\bigg)\bigg|\\
\notag \ll X^{1+\zeta}\bigg(\dfrac{1}{(Q^2B)^{\omega/2-\alpha_b}} + \dfrac{X^{\alpha_b}}{X^{\omega/2}}\bigg)(\log X)^{C'},
\end{align}
where $\widehat{1}_{\mathcal{A}_r}, \widehat{\mathfrak{f}}_{d, c_d}$, and $\alpha_b$ are given by \eqref{Arfourier}, \eqref{Fourf}, and \eqref{alphab}, respectively. Furthermore, $\omega$ and $C'$ are as in condition \ref{Cond 3} of Theorem \ref{Thm G}, and  the implicit constant in $\ll$ depends at most on $b$, $\delta$, and $\omega$.

\end{prop}

\begin{proof}
We use ideas of Maynard from the proof of \cite[Lemma 6.1]{May2}.
 For $X=b^k$, we have
\begin{align}
\widehat{\mathfrak{f}}_{d, c_d}\bigg(\dfrac{-a}{q}-\dfrac{\eta}{X}\bigg)=\sum_{\substack{n<X\\n\equiv c_d\imod d}}\mathfrak{f}(n)e\bigg(-n\bigg(\dfrac{a}{q}+\dfrac{\eta}{X}\bigg)\bigg).
\end{align}

We use condition \ref{Cond 3} with $\theta=a/q+\beta$ where $\beta=\eta/X$. Hence, for $q\sim Q$ and $(1+|\eta|)\sim B$ with $B\geq 1$, we note that $q (1 + |\beta|X)\asymp QB$ to obtain
\begin{align}
\label{minor2}\sup_{\substack{q\sim Q\\(a,q)=1\\ (|\eta|+1) \sim B}}\sum_{\substack{d\leq D\\(d, b)=1}}{\sigma}(d)\widehat{\mathfrak{f}}_{d, c_d}\bigg(\dfrac{-a}{q}-\dfrac{\eta}{X}\bigg) \ll X\bigg(\dfrac{(QB)^{\omega}}{X^\omega}+\dfrac{1}{(QB)^{\omega}}\bigg)(\log X)^{C'}.
\end{align}
By assumption $QB\ll X^{1/2}$, so that $Q^2B \ll  X$. Therefore, Lemma \ref{lemd2} implies that
\begin{align}
\label{minor3} \sum_{q\sim Q}\sum_{\substack{1\leq a <q\\(a,q)=1}}\sum_{\substack{(|\eta|+1)\sim B\\X a/q + \eta \in \mathbb{Z}}}\bigg|\widehat{1}_{\mathcal{A}_r}\bigg(\frac{a}{q} + \frac{\eta}{X}\bigg)\bigg| \ll X^{\zeta}(Q^2B)^{\alpha_b},
\end{align}
assuming that $b$ is large enough, so that $\alpha_b < 1$ and $\alpha_b \rightarrow 0$ as $b \rightarrow \infty$.
Putting the estimates from \eqref{minor2} and \eqref{minor3} together, we have
\begin{align}
\notag &\sum_{q\sim Q}\sum_{\substack{a=1\\(a,q)=1}}^q\sum_{\substack{(|\eta|+1)\sim B\\X a/q+\eta \in \mathbb{Z}}}\bigg|\widehat{1}_{\mathcal{A}_r}\bigg(\frac{a}{q}+\frac{\eta}{X}\bigg)\sum_{\substack{d\leq D\\(d, b)=1}}{\sigma}(d)\widehat{\mathfrak{f}}_{d, c_d}\bigg(\dfrac{-a}{q}-\dfrac{\eta}{X}\bigg)\bigg|\\
\notag &\ll X^{1+\zeta}\bigg(\dfrac{(QB)^{\omega}(Q^2B)^{\alpha
_b}}{X^\omega}+\dfrac{(Q^2B)^{\alpha_b}}{(QB)^{\omega}}\bigg)(\log X)^{C'}.
\end{align}
By assumption $QB\ll X^{1/2}$, and by the fact that $(QB)^{\omega}> (Q^2B)^{\omega/2}$ for $B> 1$, the above estimate is
\begin{align} 
\notag \ll X^{1+\zeta}\bigg(\dfrac{X^{\alpha_b }}{X^{\omega/2}}+ \dfrac{1}{(Q^2B)^{\omega/2-\alpha_b}}\bigg)(\log X)^{C'}.
\end{align}
This establishes the desired result.
\end{proof}

\section{Proof of Theorem \ref{Thm G} }\label{sec: proof of gen}

We are now ready to give the proof of Theorem \ref{Thm G} by combining Proposition \ref{PropMajArc} (major arcs estimate) and Proposition \ref{PropMinorArc1} (minor arcs estimate).

\begin{proof}[Proof of Theorem \ref{Thm G}]
By Fourier inversion (relation \eqref{fourinv}), we have
\begin{align}
\label{Thm7eq1} \sum_{\substack{n < X\\n\equiv c_d \imod d}}\mathfrak{f}(n)1_{\mathcal{A}_r}(n+s) = \dfrac{1}{X}\sum_{0\leq t<X}\widehat{1}_{\mathcal{A}_r}\bigg(\dfrac{t}{X}\bigg)\widehat{\mathfrak{f}}_{d, c_d}\bigg(\dfrac{-t}{X}\bigg)e\bigg(\dfrac{-st}{X}\bigg),
\end{align}
where $\widehat{1}_{\mathcal{A}_r}$ and $\widehat{\mathfrak{f}}_{d, c_d}$ are given by \eqref{Arfourier} and \eqref{Fourf}, respectively.

We consider the parameter $C>0$ to be chosen later. Then we dissect the fractions $t/X$ with $t\in [0, X)\cap \mathbb{Z}$ into two sets: major arcs $\mathfrak{M}$ and minor arcs $\mathfrak{m}$ (see Section \ref{Setup} for definition of these two sets).

We may now use Proposition \ref{PropMajArc} to estimate the major arcs $\mathfrak{M}$ contribution. We will show that
\begin{align}
\label{Thm7eq2} \sum_{\substack{d\leq D\\(d, b)=1}}&{\sigma}(d)\bigg(\dfrac{1}{X}\sum_{\substack{0\leq t<X\\t\in \mathfrak{M}}}\widehat{1}_{\mathcal{A}_r}\bigg(\dfrac{t}{X}\bigg)\widehat{\mathfrak{f}}_{d, c_d}\bigg(\dfrac{-t}{X}\bigg)e\bigg(\dfrac{-st}{X}\bigg)-\dfrac{ \lambda_{d}}{\varphi(d)}\dfrac{b}{\varphi(b)}\sum_{\substack{n< X}}1_{\mathcal{A}_r}(n)\bigg)\ll \dfrac{X^\zeta}{(\log X)^C}.
\end{align}
For brevity, let us denote
\begin{align}
\notag \mathcal{E}(d):=\dfrac{1}{X}\sum_{\substack{0\leq t<X\\t\in \mathfrak{M}}}\widehat{1}_{\mathcal{A}_r}\bigg(\dfrac{t}{X}\bigg)\widehat{\mathfrak{f}}_{d, c_d}\bigg(\dfrac{-t}{X}\bigg)e\bigg(\dfrac{-st}{X}\bigg)-\dfrac{ \lambda_{d}}{\varphi(d)}\dfrac{b}{\varphi(b)}\sum_{\substack{n< X}}1_{\mathcal{A}_r}(n).
\end{align}
We note that
\[\quad \#\mathfrak{M} \ll (\log X)^{3C},\quad  \big|\widehat{1}_{\mathcal{A}_r}(t/X)\big| \ll X^\zeta,\quad \text{and}\quad  \big|\widehat{\mathfrak{f}}_{d, c_d}(-t/X)\big| \ll X(\log X)/d\]
so that trivially, we have
\[|\mathcal{E}(d)| \ll \dfrac{X^\zeta(\log X)^{C+1}}{d}.\]
Next, we apply the Cauchy-Schwarz inequality and use the assumption that $|\sigma|\leq \uptau$ to obtain
\begin{align}
\notag \sum_{\substack{d\leq D\\(d, b)=1}}{\sigma}(d)\mathcal{E}(d)  &\ll \bigg(X^\zeta(\log X)^{3C+1}\sum_{d\leq D}\dfrac{\uptau(d)^2}{d}\bigg)^{1/2}\bigg(\sum_{\substack{d\leq D\\(d, b)=1}}\big|\mathcal{E}(d)\big|\bigg)^{1/2} \ll \dfrac{X^{\zeta}}{(\log X)^C},
\end{align}
using Proposition \ref{PropMajArc}. This completes our analysis of the major arcs.

We now use Proposition \ref{PropMinorArc1} for the remaining cases, that is, the minor arcs. We apply Dirichlet's approximation theorem to find reduced fractions $a/q$ with $1\leq q\leq X^{1/2}$ such that
\[ \bigg|\dfrac{t}{X}-\dfrac{a}{q}\bigg|\leq \dfrac{1}{qX^{1/2}}.\]
 Hence,  we have
 \[\dfrac{t}{X}=\dfrac{a}{q} + \dfrac{\eta}{X},\]
where $\max \{q, |\eta|\}\geq (\log X)^C$ and $q|\eta|\leq X^{1/2}$. Next, we perform a dyadic decomposition over $q\sim Q$ and $|\eta|+1\sim B$, so that $QB\ll X^{1/2}$. Also note that we have $\max\{Q, B\} \gg (\log X)^C$ in this case. Therefore, the contribution of minor arcs is
\begin{align}
\notag&\ll \dfrac{(\log X)^2}{X}\sum_{q\sim Q}\sum_{\substack{a=1\\(a,q)=1}}^q\sum_{\substack{B<|\eta|+1\leq 2B\\X\frac{a}{q} + \eta \in \mathbb{Z}}}\bigg|\widehat{1}_{\mathcal{A}_r}\bigg(\dfrac{a}{q}+\dfrac{\eta}{X}\bigg)\sum_{\substack{d\leq D\\(d, b)=1}}{\sigma}(d)\widehat{\mathfrak{f}}_{d, c_d}\bigg(\dfrac{-a}{q}-\dfrac{\eta}{X}\bigg)\bigg|\\
\notag &\ll X^{\zeta} \bigg(\dfrac{(\log X)^{C'}}{(\log X)^{C(\omega/2-\alpha_b)}} + \dfrac{X^{\alpha_b}(\log X)^{C'}}{X^{\omega/2}}\bigg)(\log X)^2\ll \dfrac{X^{\zeta}}{(\log X)^{A}},
\end{align}
where we have chosen $C=(A + C' +2)/(\omega/2-\alpha_b)$. Note that $\alpha_b$ goes to $0$ as $b\rightarrow \infty$ (see the relation \eqref{alphab}). In particular, since the base $b$ is sufficiently large, we have $\omega/2> \alpha_b$. Along with \eqref{Thm7eq2}, this completes the proof of Theorem \ref{Thm G}. 
\end{proof}

\bibliographystyle{alpha}

\end{document}